\documentclass[12pt]{amsart}
\def\cal{\mathcal}

\usepackage{amsmath,amssymb,amsthm,latexsym}
\usepackage{amsfonts,array,amscd}

\input epsf

\def\epsfigbracket[#1=#2]{\epsffile{ps/#2}}
\def\epsfig#1{\epsfigbracket[#1]}

\def\mod{{\rm mod}}

\def\Conj{{\cal S}}
\def\labeled{{labeled}}

\def\sigmab{{\sigma}}

\def\Ch{{\widehat{\nabla}}}
\def\Hh{{\widehat{H}}}
\def\Fh{{\widetilde{F}}}
\def\ob{{\rm ob}}

\def\w{{\rm w}}

\def\A{{\cal A}}
\def\Ab{{\bar{\A}}}
\def\L{{\cal L}}

\def\gl{{\mathfrak{gl}}}

\def\glz{{\gl_0}}

\def\osp{{\mathfrak{osp}}}
\def\Wgl{{\overline{W}_\gl}}

\def\id{{\rm id}}

\def\mod{{\rm mod}}

\def\Ker{{\rm Ker}}

\def\Conw{{\cal C}}
\def\Conjrep{{\hat{\Conj}}}

\def\Tor{{\rm Tor}}
\def\implies{\Rightarrow}
\def\Ring{R}
\def\hRing{{\widehat{\Ring}}}

\def\barz{0}
\def\baro{1}

\def\Conquo{\Conw^\#}

\def\sgn{{\rm sgn}}
\def\obG{{chord diagram}}
\def\ob{{\rm ob}}
\def\glz{{\gl_0}}
\def\glo{\glz}
\def\oK{{\overline{K}}}
\def\overK{\oK}
\def\overL{{\overline{L}}}
\def\cor{{\Z}}
\def\cori{{\Z[z]}}
\def\gen{{{\cal G}_\Sigma}}
\def\genb{{\overline{\cal G}_\Sigma}}
\def\pe{{P}}

\def\What{{\widehat{W}}}

\def\picChord{{\begin{picture}(1.1,0.6)(-0.55,-0.1)
\thicklines
\put(-0.5,0.5){\vector(0,-1){1}}
\put(0.5,0.5){\vector(0,-1){1}}
\thinlines
\put(-0.5,0){\line(1,0){1}}
\end{picture}}}

\def\picX{{\begin{picture}(1.1,0.6)(-0.55,-0.1)
\thicklines
\put(-0.5,0.5){\vector(1,-1){1}}
\put(0.5,0.5){\vector(-1,-1){1}}
\end{picture}}}

\def\picXs{{\begin{picture}(1,0.6)(-0.5,-0.1)
\setlength{\unitlength}{14pt}
\thicklines
\put(-0.5,0.5){\vector(1,-1){1}}
\put(0.5,0.5){\vector(-1,-1){1}}
\end{picture}}}

\def\picXps{{\begin{picture}(1,0.6)(-0.5,-0.1)
\setlength{\unitlength}{14pt}
\thicklines
\put(0.09,-0.09){\vector(1,-1){0.40}}
\qbezier[50](0.09,-0.09)(0.3,-0.3)(0.5,-0.5)
\qbezier[50](-0.5,0.5)(-0.3,0.3)(-0.1,0.1)
\put(0.5,0.5){\vector(-1,-1){1}}
\end{picture}}}

\def\picXms{{\begin{picture}(1,0.6)(-0.5,-0.1)
\setlength{\unitlength}{14pt}
\thicklines
\qbezier[50](0.5,0.5)(0.3,0.3)(0.1,0.1)
\put(-0.5,0.5){\vector(1,-1){1}}
\put(-0.1,-0.1){\vector(-1,-1){0.4}}
\qbezier[50](-0.1,-0.1)(-0.3,-0.3)(-0.5,-0.5)
\end{picture}}}

\def\picIIs{{\begin{picture}(1,0.6)(-0.5,-0.1)
\setlength{\unitlength}{14pt}
\thicklines
\qbezier[50](-0.5,0.5)(0,0)(-0.5,-0.5)
\qbezier[50](0.5,0.5)(0,0)(0.5,-0.5)
\put(-0.45,-0.45){\vector(-1,-1){0.05}}
\put(0.45,-0.45){\vector(1,-1){0.05}}
\end{picture}}}

\def\picBo{{\begin{picture}(1.5,0.6)(-0.6,-0.1)
\setlength{\unitlength}{15pt}
\thicklines
\put(-0.5,0.5){\vector(0,-1){1}}
\put(0.5,0.5){\vector(0,-1){1}}
\qbezier[50](0.4,0.05)(0.5,0.15)(0.6,0.25)
\qbezier[50](0.4,0.25)(0.5,0.15)(0.6,0.05)
\put(0.7,-0.2){{\makebox(0.3,0.6){\small $1$}}}
\end{picture}}}

\def\picBt{{\begin{picture}(1.5,0.6)(-0.6,-0.1)
\setlength{\unitlength}{15pt}
\thicklines
\put(-0.5,0.5){\vector(0,-1){1}}
\put(0.5,0.5){\vector(0,-1){1}}
\qbezier[50](-0.4,0.05)(-0.5,0.15)(-0.6,0.25)
\qbezier[50](-0.4,0.25)(-0.5,0.15)(-0.6,0.05)
\put(-0.3,-0.2){{\makebox(0.3,0.6){\small $1$}}}
\qbezier[50](0.4,0.05)(0.5,0.15)(0.6,0.25)
\qbezier[50](0.4,0.25)(0.5,0.15)(0.6,0.05)
\put(0.7,-0.2){{\makebox(0.3,0.6){\small $2$}}}
\end{picture}}}

\def\picuXpl{{\begin{picture}(1,0.6)(-0.5,-0.1)
\setlength{\unitlength}{17pt}
\thicklines
\put(0.09,-0.09){\line(1,-1){0.40}}
\qbezier[50](0.09,-0.09)(0.3,-0.3)(0.5,-0.5)
\qbezier[50](-0.5,0.5)(-0.3,0.3)(-0.1,0.1)
\put(0.5,0.5){\line(-1,-1){1}}
\end{picture}}}

\def\picuXml{{\begin{picture}(1,0.6)(-0.5,-0.1)
\setlength{\unitlength}{17pt}
\thicklines
\qbezier[50](0.5,0.5)(0.3,0.3)(0.1,0.1)
\put(-0.5,0.5){\line(1,-1){1}}
\put(-0.1,-0.1){\line(-1,-1){0.4}}
\qbezier[50](-0.1,-0.1)(-0.3,-0.3)(-0.5,-0.5)
\end{picture}}}

\def\picuIIl{{\begin{picture}(1,0.6)(-0.5,-0.1)
\setlength{\unitlength}{17pt}
\thicklines
\qbezier(-0.5,0.5)(0,0)(-0.5,-0.5)
\qbezier(0.5,0.5)(0,0)(0.5,-0.5)
\end{picture}}}

\def\picuInvl{{\begin{picture}(1,0.6)(-0.5,-0.1)
\setlength{\unitlength}{17pt}
\thicklines
\qbezier(-0.5,0.5)(0,0)(0.5,0.5)
\qbezier(-0.5,-0.5)(0,0)(0.5,-0.5)
\end{picture}}}

\def\picuInvmed{{\begin{picture}(1,0.6)(-0.5,-0.1)
\setlength{\unitlength}{14pt}
\thicklines
\qbezier(-0.5,0.5)(0,0)(0.5,0.5)
\qbezier(-0.5,-0.5)(0,0)(0.5,-0.5)
\end{picture}}}

\def\picuXp{{\begin{picture}(1,0.6)(-0.5,-0.1)
\setlength{\unitlength}{10pt}
\thicklines
\put(0.09,-0.09){\line(1,-1){0.40}}
\qbezier[50](0.09,-0.09)(0.3,-0.3)(0.5,-0.5)
\qbezier[50](-0.5,0.5)(-0.3,0.3)(-0.1,0.1)
\put(0.5,0.5){\line(-1,-1){1}}
\end{picture}}}

\def\picuXm{{\begin{picture}(1,0.6)(-0.5,-0.1)
\setlength{\unitlength}{10pt}
\thicklines
\qbezier[50](0.5,0.5)(0.3,0.3)(0.1,0.1)
\put(-0.5,0.5){\line(1,-1){1}}
\put(-0.1,-0.1){\line(-1,-1){0.4}}
\qbezier[50](-0.1,-0.1)(-0.3,-0.3)(-0.5,-0.5)
\end{picture}}}

\def\picuII{{\begin{picture}(1,0.6)(-0.5,-0.1)
\setlength{\unitlength}{10pt}
\thicklines
\qbezier(-0.5,0.5)(0,0)(-0.5,-0.5)
\qbezier(0.5,0.5)(0,0)(0.5,-0.5)
\end{picture}}}

\def\picuInv{{\begin{picture}(1,0.6)(-0.5,-0.1)
\setlength{\unitlength}{10pt}
\thicklines
\qbezier(-0.5,0.5)(0,0)(0.5,0.5)
\qbezier(-0.5,-0.5)(0,0)(0.5,-0.5)
\end{picture}}}

\def\picKt{{\begin{picture}(0.9,0.7)(-0.75,-0.1)
\setlength{\unitlength}{10pt}
\thicklines
\put(-0.8,1){\line(0,-1){0.64}}
\put(-0.8,-0.2){\line(0,-1){0.8}}
\qbezier[50](-0.63,-0.27)(-0.53,-0.45)(-0.25,-0.45)
\qbezier[50](-0.25,-0.45)(0.1,-0.45)(0.1,0.0)
\qbezier[50](-0.25,0.45)(0.1,0.45)(0.1,0.0)
\qbezier[50](-0.8,-0.2)(-0.8,0.45)(-0.25,0.45)
\end{picture}}}

\def\picKI{{\begin{picture}(0.4,0.7)(-0.2,-0.1)
\setlength{\unitlength}{10pt}
\thicklines
\put(0,1){\line(0,-1){2}}
\end{picture}}}

\def\picKts{{\begin{picture}(0.6,1)(-0.4,-0.167)
\setlength{\unitlength}{5pt}
\thicklines
\put(-0.8,1){\line(0,-1){0.64}}
\put(-0.8,-0.2){\line(0,-1){0.8}}
\qbezier[50](-0.63,-0.27)(-0.53,-0.45)(-0.25,-0.45)
\qbezier[50](-0.25,-0.45)(0.1,-0.45)(0.1,0.0)
\qbezier[50](-0.25,0.45)(0.1,0.45)(0.1,0.0)
\qbezier[50](-0.8,-0.2)(-0.8,0.45)(-0.25,0.45)
\end{picture}}}

\def\picKIs{{\begin{picture}(0.4,1)(-0.2,-0.167)
\setlength{\unitlength}{5pt}
\thicklines
\put(0,1){\line(0,-1){2}}
\end{picture}}}

\def\picux{{\begin{picture}(0.5,0.3)(-0.25,-0.1)
\thicklines
\qbezier(-0.2,0.2)(0,0)(0.2,-0.2)
\qbezier(0.2,0.2)(0,0)(-0.2,-0.2)
\end{picture}}}

\def\picu2{{\begin{picture}(0.5,0.3)(-0.25,-0.1)
\thicklines
\qbezier(-0.2,0.2)(0,0)(-0.2,-0.2)
\qbezier(0.2,0.2)(0,0)(0.2,-0.2)
\end{picture}}}

\def\picuinf{{\begin{picture}(0.5,0.3)(-0.25,-0.1)
\thicklines
\qbezier(-0.2,0.2)(0,0)(0.2,0.2)
\qbezier(-0.2,-0.2)(0,0)(0.2,-0.2)
\end{picture}}}

\def\DottedCircle{
\bezier{4}(0.966,-0.259)(1.04,0)(0.966,0.259)
\bezier{4}(0.966,0.259)(0.897,0.518)(0.707,0.707)
\bezier{4}(0.707,0.707)(0.518,0.897)(0.259,0.966)
\bezier{4}(0.259,0.966)(0,1.04)(-0.259,0.966)
\bezier{4}(-0.259,0.966)(-0.518,0.897)(-0.707,0.707)
\bezier{4}(-0.707,0.707)(-0.897,0.518)(-0.966,0.259)
\bezier{4}(-0.966,0.259)(-1.04,0)(-0.966,-0.259)
\bezier{4}(-0.966,-0.259)(-0.897,-0.518)(-0.707,-0.707)
\bezier{4}(-0.707,-0.707)(-0.518,-0.897)(-0.259,-0.966)
\bezier{4}(-0.259,-0.966)(0,-1.04)(0.259,-0.966)
\bezier{4}(0.259,-0.966)(0.518,-0.897)(0.707,-0.707)
\bezier{4}(0.707,-0.707)(0.897,-0.518)(0.966,-0.259)
}
\def\Endpoint[#1]{
}
\def\Arc[#1]{
\thicklines                     
\ifcase#1
\bezier{25}(0.966,-0.259)(1.04,0)(0.966,0.259)
\or
\bezier{25}(0.966,0.259)(0.897,0.518)(0.707,0.707)
\or
\bezier{25}(0.707,0.707)(0.518,0.897)(0.259,0.966)
\or
\bezier{25}(0.259,0.966)(0,1.04)(-0.259,0.966)
\or
\bezier{25}(-0.259,0.966)(-0.518,0.897)(-0.707,0.707)
\or
\bezier{25}(-0.707,0.707)(-0.897,0.518)(-0.966,0.259)
\or
\bezier{25}(-0.966,0.259)(-1.04,0)(-0.966,-0.259)
\or
\bezier{25}(-0.966,-0.259)(-0.897,-0.518)(-0.707,-0.707)
\or
\bezier{25}(-0.707,-0.707)(-0.518,-0.897)(-0.259,-0.966)
\or
\bezier{25}(-0.259,-0.966)(0,-1.04)(0.259,-0.966)
\or
\bezier{25}(0.259,-0.966)(0.518,-0.897)(0.707,-0.707)
\or
\bezier{25}(0.707,-0.707)(0.897,-0.518)(0.966,-0.259)
\fi}
\def\DottedArc[#1]{
\ifcase#1
\bezier{4}(0.966,-0.259)(1.04,0)(0.966,0.259)
\or
\bezier{4}(0.966,0.259)(0.897,0.518)(0.707,0.707)
\or
\bezier{4}(0.707,0.707)(0.518,0.897)(0.259,0.966)
\or
\bezier{4}(0.259,0.966)(0,1.04)(-0.259,0.966)
\or
\bezier{4}(-0.259,0.966)(-0.518,0.897)(-0.707,0.707)
\or
\bezier{4}(-0.707,0.707)(-0.897,0.518)(-0.966,0.259)
\or
\bezier{4}(-0.966,0.259)(-1.04,0)(-0.966,-0.259)
\or
\bezier{4}(-0.966,-0.259)(-0.897,-0.518)(-0.707,-0.707)
\or
\bezier{4}(-0.707,-0.707)(-0.518,-0.897)(-0.259,-0.966)
\or
\bezier{4}(-0.259,-0.966)(0,-1.04)(0.259,-0.966)
\or
\bezier{4}(0.259,-0.966)(0.518,-0.897)(0.707,-0.707)
\or
\bezier{4}(0.707,-0.707)(0.897,-0.518)(0.966,-0.259)
\fi}
\def\Chord[#1,#2]{
\thinlines
\ifnum#1>#2\Chord[#2,#1]
\else\ifnum#1<#2
\ifcase#1
\ifcase#2
\or\qbezier(1,0)(0.516,0.138)(0.866,0.5)
\or\qbezier(1,0)(0.45,0.26)(0.5,0.866)
\or\qbezier(1,0)(0.327,0.327)(0,1)
\or\qbezier(1,0)(0.179,0.311)(-0.5,0.866)
\or\qbezier(1,0)(0.0536,0.2)(-0.866,0.5)
\or\put(1, 0){\line(-2, 0){2}}
\or\qbezier(1,0)(0.0536,-0.2)(-0.866,-0.5)
\or\qbezier(1,0)(0.179,-0.311)(-0.5,-0.866)
\or\qbezier(1,0)(0.327,-0.327)(0,-1)
\or\qbezier(1,0)(0.45,-0.26)(0.5,-0.866)
\or\qbezier(1,0)(0.516,-0.138)(0.866,-0.5)
\fi
\or\ifcase#2\or
\or\qbezier(0.866,0.5)(0.378,0.378)(0.5,0.866)
\or\qbezier(0.866,0.5)(0.26,0.45)(0,1)
\or\qbezier(0.866,0.5)(0.12,0.446)(-0.5,0.866)
\or\qbezier(0.866,0.5)(0,0.359)(-0.866,0.5)
\or\qbezier(0.866,0.5)(-0.0536,0.2)(-1,0)
\or\put(0.866, 0.5){\line(-5, -3){1.73}}
\or\qbezier(0.866,0.5)(0.146,-0.146)(-0.5,-0.866)
\or\qbezier(0.866,0.5)(0.311,-0.179)(0,-1)
\or\qbezier(0.866,0.5)(0.446,-0.12)(0.5,-0.866)
\or\qbezier(0.866,0.5)(0.52,0)(0.866,-0.5)
\fi
\or\ifcase#2\or\or
\or\qbezier(0.5,0.866)(0.138,0.516)(0,1)
\or\qbezier(0.5,0.866)(0,0.52)(-0.5,0.866)
\or\qbezier(0.5,0.866)(-0.12,0.446)(-0.866,0.5)
\or\qbezier(0.5,0.866)(-0.179,0.311)(-1,0)
\or\qbezier(0.5,0.866)(-0.146,0.146)(-0.866,-0.5)
\or\put(0.5, 0.866){\line(-3, -5){1}}
\or\qbezier(0.5,0.866)(0.2,-0.0536)(0,-1)
\or\qbezier(0.5,0.866)(0.359,0)(0.5,-0.866)
\or\qbezier(0.5,0.866)(0.446,0.12)(0.866,-0.5)
\fi
\or\ifcase#2\or\or\or
\or\qbezier(0,1.)(-0.138,0.516)(-0.5,0.866)
\or\qbezier(0,1.)(-0.26,0.45)(-0.866,0.5)
\or\qbezier(0,1.)(-0.327,0.327)(-1,0)
\or\qbezier(0,1.)(-0.311,0.179)(-0.866,-0.5)
\or\qbezier(0,1.)(-0.2,0.0536)(-0.5,-0.866)
\or\put(0, 1){\line(0, -2){2}}
\or\qbezier(0,1.)(0.2,0.0536)(0.5,-0.866)
\or\qbezier(0,1.)(0.311,0.179)(0.866,-0.5)
\fi
\or\ifcase#2\or\or\or\or
\or\qbezier(-0.5,0.866)(-0.378,0.378)(-0.866,0.5)
\or\qbezier(-0.5,0.866)(-0.45,0.26)(-1,0)
\or\qbezier(-0.5,0.866)(-0.446,0.12)(-0.866,-0.5)
\or\qbezier(-0.5,0.866)(-0.359,0)(-0.5,-0.866)
\or\qbezier(-0.5,0.866)(-0.2,-0.0536)(0,-1)
\or\put(-0.5, 0.866){\line(3, -5){1}}
\or\qbezier(-0.5,0.866)(0.146,0.146)(0.866,-0.5)
\fi
\or\ifcase#2\or\or\or\or\or
\or\qbezier(-0.866,0.5)(-0.516,0.138)(-1,0)
\or\qbezier(-0.866,0.5)(-0.52,0)(-0.866,-0.5)
\or\qbezier(-0.866,0.5)(-0.446,-0.12)(-0.5,-0.866)
\or\qbezier(-0.866,0.5)(-0.311,-0.179)(0,-1)
\or\qbezier(-0.866,0.5)(-0.146,-0.146)(0.5,-0.866)
\or\put(-0.866, 0.5){\line(5, -3){1.73}}
\fi
\or\ifcase#2\or\or\or\or\or\or
\or\qbezier(-1,0)(-0.516,-0.138)(-0.866,-0.5)
\or\qbezier(-1,0)(-0.45,-0.26)(-0.5,-0.866)
\or\qbezier(-1,0)(-0.327,-0.327)(0,-1)
\or\qbezier(-1,0)(-0.179,-0.311)(0.5,-0.866)
\or\qbezier(-1,0)(-0.0536,-0.2)(0.866,-0.5)
\fi
\or\ifcase#2\or\or\or\or\or\or\or
\or\qbezier(-0.866,-0.5)(-0.378,-0.378)(-0.5,-0.866)
\or\qbezier(-0.866,-0.5)(-0.26,-0.45)(0,-1)
\or\qbezier(-0.866,-0.5)(-0.12,-0.446)(0.5,-0.866)
\or\qbezier(-0.866,-0.5)(0,-0.359)(0.866,-0.5)
\fi
\or\ifcase#2\or\or\or\or\or\or\or\or
\or\qbezier(-0.5,-0.866)(-0.138,-0.516)(0,-1)
\or\qbezier(-0.5,-0.866)(0,-0.52)(0.5,-0.866)
\or\qbezier(-0.5,-0.866)(0.12,-0.446)(0.866,-0.5)
\fi
\or\ifcase#2\or\or\or\or\or\or\or\or\or
\or\qbezier(0,-1.)(0.138,-0.516)(0.5,-0.866)
\or\qbezier(0,-1.)(0.26,-0.45)(0.866,-0.5)
\fi
\or\ifcase#2\or\or\or\or\or\or\or\or\or\or
\or\qbezier(0.5,-0.866)(0.378,-0.378)(0.866,-0.5)
\fi\fi\fi\fi}
\def\FullChord[#1,#2]{
\Endpoint[#1]
\Endpoint[#2]
\Arc[#1]
\Arc[#2]
\Chord[#1,#2]
}
\def\EndChord[#1,#2]{
\Endpoint[#1]
\Endpoint[#2]
\Chord[#1,#2]
}
\def\Picture#1{
\begin{picture}(2,1)(-1,-0.167) 
#1
\end{picture}
}
\def\DottedChordDiagram[#1,#2]{
\Picture{\DottedCircle \FullChord[#1,#2]}
}
\newtheorem{lemma}{Lemma}
\newtheorem{prop}[lemma]{Proposition}

\newtheorem{theorem}[lemma]{Theorem}
\newtheorem{coro}[lemma]{Corollary}
\newtheorem{defi}[lemma]{Definition}

\newcommand{\N}{{\ensuremath{\mathbb N}}}
\newcommand{\Z}{{\ensuremath{\mathbb Z}}}
\newcommand{\Q}{{\ensuremath{\mathbb Q}}}
\newcommand{\R}{{\ensuremath{\mathbb R}}}

\begin{document}
{
\title{Skein modules of links in cylinders 
over surfaces}
\author{Jens Lieberum}
\address{\hskip-\parindent
        Jens Lieberum\\
        Mathematisches Institut\\
        Universit\"at Basel\\
        Rheinsprung 21\\
        CH-4051 Basel}
\email{lieberum@msri.org}
\date{}

\thanks{
I would like to thank C.\ Riedtmann and C.\ Kassel for helpful remarks. 
I thank the German Academic Exchange Service and the
Schweizerischer Nationalfonds for financial support.
Research at MSRI is supported in part by NSF grant DMS-9701755.
}
\maketitle

\begin{abstract}
{\parindent0cm We define the Conway skein module~$\Conw(M)$ of 
ordered based links in a~$3$\nobreakdash-ma\-ni\-fold~$M$. This 
module gives rise to~$\Conw(M)$-valued invariants of usual links 
in~$M$.
We determine a 
basis of the~$\Z[z]$-module~$\Conw(\Sigma\times[0,1])/
\Tor(\Conw(\Sigma\times[0,1]))$ where $\Sigma$ is the real 
projective plane or a surface with boundary.
For cylinders over the M\"obius strip or the projective plane
we derive special properties of the Conway skein module,
among them a refinement
of a theorem of Hartley and Kawauchi about the Conway polynomial of strongly 
positive amphicheiral
knots in~$S^3$.
In addition, we determine the Homfly and
Kauffman skein modules of~$\Sigma\times [0,1]$ where~$\Sigma$ is an oriented surface with 
boundary.
\bigskip



}
\end{abstract}

\section*{Introduction}
\parindent0cm
In 1969 J.\ H.\ Conway found a normalized version~$\tilde{\nabla}$ 
of the Alexander polynomial of a link.  
It satisfies a so-called skein relation on certain triples of links differing only by
local modifications, 
namely 

\begin{equation}\label{e:skein}
\tilde{\nabla}(\picXps)-\tilde{\nabla}(\picXms)=z\tilde{\nabla}(\picIIs)
\end{equation}

\noindent (see \cite{Ale},\cite{Con}).  
The Conway polynomial also satisfies~$\tilde{\nabla}(O_n)=\delta_{1n}$, where~$O_n$ 
is the trivial link with~$n$ components and~$\delta_{ij}$ is~$1$ if~$i=j$ 
and~$0$ otherwise.
The skein relation together with~$\tilde{\nabla}(O_1)=1$ 
uniquely determine
the isotopy invariant~$\tilde{\nabla}(L)$ for every link~$L$.
More precisely, the so-called skein module~$\Conw(S^3)$ generated 
over~$\Z[z^{\pm 1}]$ by isotopy classes of links
modulo the skein relation is free with basis~$O_1$. 
The Jones, Homfly and Kauffman polynomials of links in~$S^3$ (\cite{Jon}, 
\cite{HOM}, \cite{Ka2})
can be characterized by similar skein relations.
In~1988 J.\ Przytycki and V.\ Turaev
generalized the Homfly and Kauffman polynomials to links in a solid torus
and determined the corresponding skein modules.
Since then,
various versions of skein modules related to these polynomials
have been studied (\cite{Pr2}), among them 
the Homfly skein module of links in a handlebody (\cite{Pr1}).

In~1990 Vassiliev invariants were introduced (\cite{Vas}). They can be characterized as follows:
every link invariant~$v$ with values in an abelian group
can recursively 
be extended to an invariant of singular links by $v(\picXs)=v(\picXps)-v(\picXms)$.
Vassiliev invariants~$v$ of degree~$n$ of links are those invariants that
vanish on singular links with~$n+1$ double points.
If the $3$-manifold is not orientable, then there exists no preferred choice of a sign in
$v(\picXs)=\pm v(\picXps)\mp v(\picXms)$, but it still makes sense to speak of 
Vassiliev invariants
of degree~$n$. Let $I=[0,1]$.
In 1996, Andersen, Mattes and Reshetikhin found a universal Vassiliev 
invariant~$Z_{\Sigma\times I}$
of links in a cylinder over an orientable surface~$\Sigma$ with boundary (\cite{AMR}).
We use~$Z_{\Sigma\times I}$ to 
determine the Homfly and Kauffman skein modules of~$\Sigma\times I$. 
In comparison with~\cite{Pr1} our choice of a basis of the Homfly skein module 
has some technical advantages.
The result seems to be new for the Kauffman skein module. 
 

In \cite{Lie} the invariant~$Z_{\Sigma\times I}$ was extended to a universal Vassiliev
invariant of
links in cylinders over non-orientable
surfaces with boundary and over the real projective plane~$P^2$.
The main idea that lead to this article was to extract an 
explicitly computable link invariant
from~$Z_{\Sigma\times I}$ in the non-orientable case.
We define a skein module~$\Conw(M)$ of a $3$-manifold~$M$ using a version of 
the skein relation in Equation~(\ref{e:skein}) for ordered based links.
In Theorem~\ref{t:conwaypoly} we establish an isomorphism 
of~$\cori$-modules

\begin{equation}\label{e:conwaymod}
\Conquo(\Sigma\times I)\cong \cori\otimes_\Z\Conw_0^\#(\Sigma\times I),
\end{equation}

\noindent where $\Sigma=P^2$ or $\Sigma$ is a non-orientable surface with boundary,
$\Conquo(M)$ is the quotient of 
the $\cori$-module~$\Conw(M)$ by its 
torsion submodule, and~$\Conw_0^\#(M)$ is 
the quotient of the~$\cor$-module~$\Conw(M)/(z\Conw(M))$
by its torsion submodule.
We determine the structure of~$\Conw_0^\#(M)$ for an arbitrary $3$-manifold. It is 
isomorphic to a tensor product of a polynomial algebra and an exterior algebra.
A basis of~$\Conquo(\Sigma\times I)$ is given by certain descending links.
We describe explicitly how to compute a Conway polynomial 
inducing the isomorphism in Equation~(\ref{e:conwaymod}) and how to obtain invariants of
links (without order or basepoints) from the Conway skein module. 
When~$\Sigma$ is the M\"obius strip~$X$ or~$\Sigma=P^2$ 
we show that the isomorphism in Equation~(\ref{e:conwaymod}) is an isomorphism of 
involutive algebras (Theorem~\ref{t:specase}).
For strongly 
positive amphicheiral knots in~$S^3$ a theorem of Hartley and Kawauchi (\cite{KaH}) says that
the Conway polynomial is a square. 
We prove a refinement of this theorem by showing that the
square root of this polynomial can be computed directly by using the Conway skein relation 
for ordered based links in~$P^2\times I$ (Theorem~\ref{t:Csquare}).

The paper is organized as follows. In Section~\ref{s:defc} we 
define~$\Conw(M)$ and prove some of its general properties. After 
this section we concentrate on the case~$M=\Sigma\times I$. In 
Section~\ref{s:surf} we introduce decomposed surfaces and make a 
choice of representatives of conjugacy classes in the fundamental 
group of~$\Sigma$. 
This choice is important in the definition of a
descending link.
Only a careful choice leads to a basis of~$\Conquo(\Sigma\times I)$ 
over~$\Z[z]$.
In Section~\ref{s:mr} we state the
Theorems~\ref{t:conwaypoly}, \ref{t:specase}, and~\ref{t:Csquare} mentioned above.
Sections~\ref{s:genc} to~\ref{s:squareproof} are devoted to the proof of 
these theorems. In Section~\ref{s:genc} we show that descending 
knots and so-called cabled descending knots generate~$\Conw(\Sigma\times 
I)$ as an algebra.  By making computations in~$\Conquo(X\times I)$ we show in 
Section~\ref{s:gencx} that cabled descending knots vanish 
in~$\Conquo(\Sigma\times I)$. In Section~\ref{s:wC} we introduce 
labeled ordered based chord diagrams and prove that a map 
$W_\glo^\ob$ defined on these diagrams is compatible with certain 
relations. In Section~\ref{s:eC} we compose the universal 
Vassiliev invariant~$Z_{\Sigma\times I}$ adapted to ordered based 
links with~$W_\glo^\ob$ to prove Theorems~\ref{t:conwaypoly} and~\ref{t:specase}.
In Section~\ref{s:squareproof} 
we prove Theorem~\ref{t:Csquare} by using a map between chord diagrams
induced by the covering~$S^2\times I\longrightarrow P^2\times I$.
In an appendix we determine the structure of the Homfly and 
Kauffman skein modules of
$\Sigma\times I$ for oriented surfaces~$\Sigma$ with boundary.

\parindent0cm
\section{The Conway skein module}\label{s:defc}

Throughout this paper we work in the PL-category.
Let~$M$ be a connected $3$-manifold 
equipped with an oriented subset~$U$ homeomorphic to a closed $3$-ball~$D^3$. 
We denote the fundamental group of~$M$ 
with respect to a basepoint in~$U$ by~$\pi_1(M)$. The first Stiefel-Whitney class  
induces a homomorphism $\sigmab:\pi_1(M)\longrightarrow \Z/2$ called {\em orientation character}
of~$M$. 
Explicitly,
if a local orientation at the basepoint stays
the same when we push it around~$w\in\pi_1(M)$,
then we have $\sigmab(w)=0$, and~$\sigmab(w)=1$ otherwise.

Ordered based links in~$M$ are 
links in~$M$ with a linear order of their components and
with a basepoint in~$U$ on each component. 
In particular, every component of an ordered based link meets~$U$.
The basepoints are not allowed to leave~$U$ 
during an isotopy (resp.\ homotopy) of ordered based links. 
The homotopy class of a based knot~$K$ can canonically be considered as an
element~$\alpha(K)\in\pi_1(M)$. 
We denote~$\sigma(\alpha(K))$ simply 
by~$\sigma(K)$.
Notice that~$\sigma(K)$ does not depend on the basepoint of~$K$.
We will regard the empty link~$\emptyset$ as an ordered based link because this
will simplify slightly the statement of our results.

In Figure~\ref{f:hpict} parts of ordered based links are shown. The numbers~$1$ and~$2$
near the basepoints in this figure indicate that the corresponding component is first
or second in the order of this link.

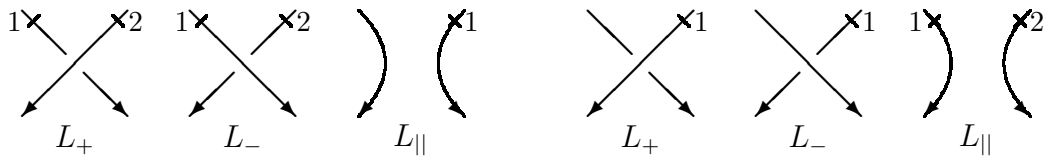
\begin{figure}[!h]
$$
\begin{picture}(2,1)(-1,-0.167)
\thicklines
\put(1,1){\vector(-1,-1){2}}
\put(-1,1){\line(1,-1){0.83}}
\put(0.17,-0.17){\vector(1,-1){0.83}}
\qbezier[20](0.7,0.9)(0.8,0.8)(0.9,0.7)
\qbezier[20](-0.7,0.9)(-0.8,0.8)(-0.9,0.7)
\put(-1.3,0.5){\makebox(0.3,0.6){\small $1$}}
\put(1.0,0.5){\makebox(0.3,0.6){\small $2$}}
\put(-1,-1.8){\makebox(2,0.7){$L_+$}}
\end{picture}
\qquad
\Picture{
\thicklines
\put(1,1){\line(-1,-1){0.83}}
\put(-0.17,-0.17){\vector(-1,-1){0.83}}
\put(-1,1){\vector(1,-1){2}}
\qbezier[20](0.7,0.9)(0.8,0.8)(0.9,0.7)
\qbezier[20](-0.7,0.9)(-0.8,0.8)(-0.9,0.7)
\put(-1.3,0.5){\makebox(0.3,0.6){\small $1$}}
\put(1.0,0.5){\makebox(0.3,0.6){\small $2$}}
\put(-1,-1.8){\makebox(2,0.7){$L_-$}}
}
\qquad
\begin{picture}(2,1)(-1,-0.167)
\thicklines
\qbezier[90](1,1)(0,0)(1,-1)
\qbezier[90](-1,1)(0,0)(-1,-1)
\put(-0.9,-0.9){\vector(-1,-1){0.1}}
\put(0.9,-0.9){\vector(1,-1){0.1}}
\qbezier[20](0.71,0.9)(0.81,0.8)(0.91,0.7)
\qbezier[20](0.91,0.9)(0.81,0.8)(0.71,0.7)
\put(1.0,0.5){\makebox(0.3,0.6){\small $1$}}
\put(-1,-1.8){\makebox(2,0.7){$L_{\mid\mid}$}}
\end{picture}
\qquad\qquad
\Picture{
\thicklines
\put(1,1){\vector(-1,-1){2}}
\put(-1,1){\line(1,-1){0.83}}
\put(0.17,-0.17){\vector(1,-1){0.83}}
\qbezier[20](0.7,0.9)(0.8,0.8)(0.9,0.7)
\put(1.0,0.5){\makebox(0.3,0.6){\small $1$}}
\put(-1,-1.8){\makebox(2,0.7){$L_+$}}
}
\qquad
\Picture{
\thicklines
\put(1,1){\line(-1,-1){0.83}}
\put(-0.17,-0.17){\vector(-1,-1){0.83}}
\put(-1,1){\vector(1,-1){2}}
\qbezier[20](0.7,0.9)(0.8,0.8)(0.9,0.7)
\put(1.0,0.5){\makebox(0.3,0.6){\small $1$}}
\put(-1,-1.8){\makebox(2,0.7){$L_-$}}
}
\qquad
\Picture{
\thicklines
\qbezier[90](1,1)(0,0)(1,-1)
\qbezier[90](-1,1)(0,0)(-1,-1)
\put(-0.9,-0.9){\vector(-1,-1){0.1}}
\put(0.9,-0.9){\vector(1,-1){0.1}}
\qbezier[20](0.71,0.9)(0.81,0.8)(0.91,0.7)
\qbezier[20](0.91,0.9)(0.81,0.8)(0.71,0.7)
\qbezier[20](-0.71,0.9)(-0.81,0.8)(-0.91,0.7)
\qbezier[20](-0.91,0.9)(-0.81,0.8)(-0.71,0.7)
\put(-1.3,0.5){\makebox(0.3,0.6){\small $1$}}
\put(1.0,0.5){\makebox(0.3,0.6){\small $2$}}
\put(-1,-1.8){\makebox(2,0.7){$L_{\mid\mid}$}}
}
$$\nopagebreak\vspace*{15pt}
\caption{Skein triples of ordered based links}\label{f:hpict}
\end{figure}

Ordered based links~$(L_+,L_-,L_{\vert\vert})$
form a {\em skein triple} if 
they differ inside of a ball $D^3\subset U$ as 
one of the two triples~$(L_+,L_-,L_{\vert\vert})$ in Figure~\ref{f:hpict}, they coincide on the
outside of this ball, and the order and basepoints of the remaining components also coincide.
We say that~$L_+$ and~$L_-$ are related by changing a crossing and~$L_{\vert\vert}$ is
obtained from~$L_+$ or~$L_-$ by splicing a crossing.

\begin{defi}\label{d:csm}
Let $\Conw(M)$ be the $\cori$-module generated by isotopy classes of 
ordered based links and the Relations $(Ord)$, $(Bas)$ and $(Skein)$ below.
\end{defi}


$(Ord): L=(-1)^{\sigmab(K_1)\sigmab(K_2)} L'$, 
where $K_1, K_2$ are two components
of an arbitrary ordered based link~$L$ that are neighbors in
the order of the components of~$L$. The ordered based link $L'$ is  
equal to~$L$ except that the order of~$K_1$ and~$K_2$ is interchanged.

\smallskip

$(Bas): L=(-1)^{\sigmab(a)\sigmab(b)}L'$, where~$L$ and~$L'$ are ordered based links that are
equal
except for one basepoint on a component~$K$ of~$L$
and~$K'$ of~$L'$. The element 
$a\in\pi_1(\Sigma\times I)$ is given by the curve described by the basepoint of~$K$ 
when it is pushed along~$K$ 
until it coincides with the
basepoint of~$K'$, and~$b=\alpha(K')a^{-1}$. 

\smallskip

$(Skein): L_+-L_-=z L_{\vert\vert}$, where $(L_+,L_-,L_{\vert\vert})$ 
is a skein triple of ordered based links.

\medskip

Let the homomorphism of rings~$\varphi:\cori\longrightarrow\cori$ be given 
by~$\varphi(p)(z)=p(-z)$.
We turn~$\Conw$ into a functor in the sense of the following proposition.

\begin{prop}\label{p:emb}
Let the $3$-manifolds $M_i$ ($i=1,2$) be equipped with the oriented subsets $U_i\cong D^3$, 
and let~$i:M_1\longrightarrow M_2$ be an embedding with~$i(U_1)\subseteq U_2$.

(1) If the orientations of~$i(U_1)$ and~$U_2$ coincide, 
then~$i$ induces a morphism of~$\Z[z]$-modules $i_*:\Conw(M_1)\longrightarrow\Conw(M_2)$.

(2) If the orientations of~$i(U_1)$ and~$U_2$ are opposite, 
then~$i$ induces a morphism of~$\Z$-modules $i_*:\Conw(M_1)\longrightarrow\Conw(M_2)$
satisfying~$i_*(pm)=\varphi(p)i_*(m)$ for all~$p\in\cori$ and~$m\in\Conw(M_1)$.
\end{prop}
\begin{proof}[Sketch of proof]
(1) is obvious.

(2) The only defining relation of~$\Conw(M_1)$ that depends on the orientation
of~$U_1$ is the skein relation. Reversing the orientation of~$U_1$ interchanges the 
parts~$L_+$ and~$L_-$ of a link which can be compensated by replacing~$z$ by~$-z$.
This way we obtain Part~(2) of the proposition by using Part~(1).
\end{proof}

\medskip

The existence of a suitable homeomorphism from~$M$ to itself implies the following corollary
which explains why we do not refer explicitly to~$U$ in the name of~$\Conw(M)$.

\begin{coro}\label{c:Unotimp}
If the $3$-manifold $M$ is oriented (resp.\ non-orientable), then 
for different choices of~$U_i\subset M_i=M$ with~$U_i\cong D^3$ 
equipped with the induced (resp.\ an arbitrary) orientation the~$\cori$-modules~$\Conw(M_i)$
are isomorphic.
\end{coro}

Let~$c$ be an arbitrary crossing of an ordered based link. We apply an isotopy such 
that~$c$ is contained in~$U$. Then we use Relations~$(Ord)$ and~$(Bas)$ 
such that~$c$ becomes a crossing involved in a skein triple.
The skein relation then implies that the
$\cor$-module $\Conw_0(M)=\Conw(M)/(z\Conw(M))$ is generated by homotopy
classes of ordered based links.
We turn a disjoint union of ordered based links~$L=L_1\cup L_2$ into an ordered 
based link 
by extending the order of the components of~$L_i$ to~$L$ by~$K_1<K_2$ for based 
knots~$K_i\subset L_i$.
This operation is not well-defined on~$\Conw(M)$, but
it turns~$\Conw_0(M)$ into an associative ring 
with~$1$\nobreakdash-ele\-ment~$\emptyset$. 
Let $T=\Tor(\Conw_0(M))$ be the
torsion submodule of the $\Z$-module $\Conw_0(M)$ which is a two-sided ideal 
of~$\Conw_0(M)$. Define $\Conw_0^\#(M)=\Conw_0(M)/T$.
Let us describe the structure of~$\Conw_0^\#(M)$.
Denote the conjugacy class of an element~$w\in\pi_1(M)$
by~$[w]$. Let

\begin{eqnarray}
& & \Conj_1=\{[w]\,\vert\,w\in\pi_1(M), \sigmab(w)=\baro\},\\
& & \Conj_0^*=\left\{[w]\,\vert\, w\in\pi_1(M),  
\sigmab(w)=\barz\mbox{ and }
N(w)\subseteq\Ker\,\sigmab\right\},
\end{eqnarray}

where $N(w)=\{v\in\pi_1(M)\,\vert\, vw=wv\}$ is the normalizer of~$w$. 
For a set~$X$, let~$V_X$ be the free $\cor$-module with
basis~$\{t_x\;\vert\; x\in X\}$. Denote the symmetric algebra on a free 
$\cor$-module~$V$ 
by~$S(V)$ and the exterior algebra by~$\Lambda(V)$. 
For $w\in\pi_1(M)$ we denote the corresponding based knot in~$\Conw_0^\#(M)$ by~$\overK_w$.

\begin{prop}\label{p:ringiso}
For every choice~$\Conjrep_0^*$ of representatives of the conjugacy 
classes in~$\Conj_0^*$ exists an isomorphism of algebras

$$
\psi:S\left(V_{\Conjrep_0^*}\right)\otimes
\Lambda\left(V_{\Conj_1}\right)\longrightarrow \Conw_0^\#(M).
$$

given by
mapping $t_v\otimes 1$ to $\overK_v$ ($v\in\Conjrep_0^*$) and 
$1\otimes t_{[w]}$ to $\overK_v$ for some~$v\in [w]\in\Conj_1$.
\end{prop}
\begin{proof}
The $\Z$-module $\Conw_0(M)$ can be described by  
non-commutative generators $k_w$ ($w\in\pi_1(M)$)
and the Relations 

\begin{eqnarray*}
(Ord') & : & k_vk_w=(-1)^{\sigmab(v)\sigmab(w)}k_wk_v\\
(Bas') & : & k_{vw}=(-1)^{\sigmab(v)\sigmab(w)}k_{wv}
\end{eqnarray*}

for all $v,w\in\pi_1(M)$.
By the Relation $(Bas')$ we have $k_a=k_b$ for all $a,b\in [a]\in\Conj_1$.
The Relation~$(Ord')$ then implies that~$\psi$ is well-defined.
If~$\sigmab(w)=0$ and there exists~$v\in\pi_1(M)$ with~$vw=wv$ 
and~$\sigmab(v)=1$, then 

\begin{equation}
k_w=k_{vwv^{-1}}=(-1)^{\sigmab(v)\sigmab(wv^{-1})}k_{wv^{-1}v}=-k_w
\end{equation}

implying $k_w=0\in\Conw_0^\#(M)$. 
If~$\sigma(w)=1$, then~$k_w^2=-k_w^2$ implies~$k_w^2=0\in\Conw_0^\#(M)$.
Therefore~$\psi$ is surjective.
Since we know a presentation of~$\Conw_0(M)$ by generators and relations,
it is easy to verify that an inverse map of~$\psi$ is 
well-defined. 
\end{proof}

\medskip

Now we fix a choice of a linearly 
ordered set~$\Conjrep^*$ consisting of representatives of~$\Conj_0^*\cup\Conj_1$.
Then by Proposition~\ref{p:ringiso} we can equip~$\Conw_0^\#(M)$ with the basis 

$$
\overK_{w_1}\ldots \overK_{w_n}\quad\mbox{ where }\quad w_i\in\Conjrep^*,
w_i\leq w_{i+1}, (w_i=w_{i+1}\implies \sigmab(w_i)=0).
$$

We define the linear map $\theta:\Conw_0^\#(M)\longrightarrow \cor$\label{d:theta}
to be equal to~$1$ on these basis elements. For non-orientable~$M$ 
the definition of~$\theta$ depends on the representatives chosen 
for~$\Conj_0^*$ and on the order chosen on~$\Conj_1$.

The module~$\Conw(M)$ gives rise to invariants of links in~$M$ (without order 
or basepoints) as follows. 
Given a link~$L$, we choose a
basepoint~$b_i$ on each component and pull all basepoints into~$U$ along disjoint
paths~$\beta_i$. 
Then we choose an arbitrary 
order of the components of~$L$ and obtain an ordered based link~$L_b$.
We denote the image of~$L_b\in\Conw(M)$ under the~$\cor$-linear 
projection~$\Conw(M)\longrightarrow
\Conw_0^\#(M)$ by~$\overline{L}_b$.
Define~$v_i(L)$ and~$C(L)$ by

\begin{eqnarray}
\Conw(M)/z^{i+1}\Conw(M)\;\ \ni & v_i(L) & =\ \;\theta\left(\overline{L}_b\right)\overL_b,\\
\{\{a,b\}\mid a,b\in\Conw(M)\}\;\ \ni &
C(L) & =\ \; 
\left\{\begin{array}{ll}
\{\theta\left(\overline{L}_b\right)\overL_b\} & \mbox{if $\theta\left(\overline{L}_b\right)
\not=0$,}\\
\{\pm \overL_b\} & \mbox{if $\theta\left(\overline{L}_b\right)=0$}.\label{e:defCL}
\end{array}
\right.
\end{eqnarray}

\begin{prop}\label{p:linkinv}
The maps $C$ and~$v_i$ ($i\geq 0$) are isotopy invariants of links. In addition,
the maps $v_i$ are Vassiliev invariants of degree~$i$. 
\end{prop}
\begin{proof}
Let $L$ be a link and let $L_a$, $L_b$ be ordered based links 
obtained from $L$ by
choosing basepoints $a_i$ and $b_i$ on the components of
$L$ that are pulled into $U$ along paths $\alpha_i$ and $\beta_i$ respectively.
We first assume that all the paths $\alpha_i$ and $\beta_i$ are disjoint. 
Then we can pull the points $b_i$ on $L_a$ into $U$ along $\beta_i$ and the points 
$a_i$ on $L_b$ into $U$ along $\alpha_i$.
This shows that we can pass from~$L_a$ to~$L_b$
by isotopies of ordered based links and by application of the Relations~$(Ord)$ and
$(Bas)$. The Relations~$(Ord)$ and $(Bas)$ 
only influence signs and $\theta\left(\overline{L}_b\right)$ changes signs simultaneously 
with~$\overL_b$.
If the paths $\alpha_i$ and $\beta_i$ are not disjoint,
then we choose basepoints $c_i$ of $L$ and paths $\gamma_i$ disjoint to
$\alpha_i$ and $\beta_i$ and apply the argument above to the pairs $(L_a,L_c)$ and $(L_c,L_b)$.
This shows that the maps~$v_i$ and~$C$ are well-defined isotopy invariants of links.
It follows from the skein relation by the same arguments as for the 
usual Conway polynomial of links
that the maps $v_i$ are Vassiliev invariants of degree~$i$ (see \cite{Vog}).
\end{proof}

\medskip

\section{Decomposed surfaces}\label{s:surf}

Let us now consider the special case~$M=\Sigma\times I$, where~$\Sigma$ is a 
connected compact surface and~$I=[0,1]$ is an oriented interval. 
First assume that~$\partial \Sigma\not=\emptyset$. We choose an oriented 
subset~$B_0\cong D^2$ 
of~$\Sigma$ such that~$B_0\cap\partial\Sigma\not=\emptyset$.
Define~$\Conw(\Sigma\times I)$ with 
respect to~$U=B_0\times I\cong D^3$.
Fix the choice of a point~$u\in B_0\cap \partial \Sigma$.
For the following definition we need some notation: we denote the fundamental group of~$\Sigma$
with respect to a basepoint in~$B_0$ by~$\pi_1(\Sigma)$.
Let~$p_\Sigma$ and~$p_I$ be the projections from~$\Sigma\times I$ to~$\Sigma$ and~$I$ 
respectively.
For~$x\in\Sigma\times I$ we call the value~$p_I(x)$ the height of~$x$.

\begin{defi}
(1) A based knot $K\subset\Sigma\times I$ 
is called {\em descending with respect to its basepoint}, if 
the height of~$K$ is descending when we travel along~$K$ starting at the basepoint~$A$ 
and following the orientation of~$K$ until we reach a point~$A'$ 
with~$p_\Sigma(A)=p_\Sigma(A')$ from which~$K$ leads back to~$A$ by increasing 
the height and by keeping the projection to~$\Sigma$ constant.

(2) A based knot~$K\subset\Sigma\times I$ is called a {\em descending knot}, 
if~$K$ is descending with respect to its basepoint~$A$
and there exists a neighborhood~$V\cong D^2$ of~$u$
in~$B_0$ such that~$(V\times I)\cap K$ 
consists of an interval containing~$A$.
\end{defi}

The important property of descending knots is that for~$w\in\pi_1(\Sigma)$ 
there exists up to isotopy exactly one
descending knot~$K_w$ such that~$\alpha(K_w)=w$.
If~$v,w\in\pi_1(\Sigma)$ are conjugate then one can pass from~$K_v$ to~$K_w$ by 
isotopies and crossing
changes and by moving the basepoint. 
In this section we will choose a set~$\Conjrep$ of representatives of conjugacy classes 
in~$\pi_1(\Sigma)$ with good properties (compare the example preceding Equation~\ref{e:nogen}
for a choice with bad properties).

\begin{defi}
A {\em decomposed surface} is a triple $(\Sigma,(B_0,B_1,\ldots,B_k),u)$ 
where~$\Sigma$ is a surface, $I^2\cong B_0\subset\Sigma$ is oriented, 
$u\in B_0\cap\partial\Sigma$, and
$I^2\cong B_i\subset\Sigma$ such that $\bigcup_{i=0}^k B_i=\Sigma$ and
Equations~(\ref{e:decsur1}) and~(\ref{e:decsur2}) are satisfied.

\begin{eqnarray}
& & B_0\cap B_i\cong I\times\{0,1\}\mbox{ for }i\geq 1,\label{e:decsur1}\\ 
& & B_i\cap B_j=\emptyset\mbox{ for }i\not=j, i,j\geq 1.\label{e:decsur2}
\end{eqnarray}
\end{defi}

We represent decomposed surfaces graphically as shown in Figure~\ref{f:decsur}
by an example. 
By convention the orientation of~$B_0$ is counterclockwise in this figure.
We will also represent~$\Sigma\times I$ as in 
Figure~\ref{f:decsur} where we assume
that on~$B_0$ the interval~$I$ directs towards the reader.

\begin{figure}[!h]
\centering
\setbox1=\hbox{\input{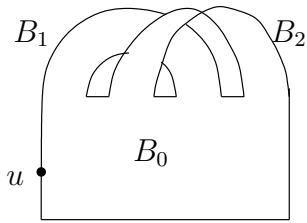}}
$\vcenter{\box1}$
\caption{A picture of a decomposed surface}\label{f:decsur}
\end{figure}

The decomposition of~$\Sigma$
determines generators~$x_i$ ($i=1,\ldots,k$)
of~$\pi_1(\Sigma)$
that pass through the band~$B_i$ exactly once in the clockwise sense in our pictures
and that do not meet the other bands~$B_j$ ($j\geq 1, i\not=j$).
For~$w\in\pi_1(\Sigma)$ we denote by~$\ell(w)$ 
the length of the unique reduced expression
of~$w$ in the generators~$x_i$ of the free group~$\pi_1(\Sigma)$
and we denote the generator at position~$i$ 
($1\leq i\leq \ell(w)$) in this expression
by~$w(i)$.
Traveling along~$\partial B_0$ in the clockwise sense starting at~$u$
we number the~$2k$ intervals~$B_i\cap B_0$ ($i=1,\ldots,k$) from~$1$ to~$2k$. We denote these
intervals by~$I_\nu$ ($\nu\in\{1,\ldots,2k\}$). 
This numbering determines maps
$s,d:\pi_1(\Sigma)\longrightarrow\{1,\ldots,2k\}$
where~$I_{s(w)}$ (resp.~$I_{d(w)}$) is the interval 
that a generic 
representative of~$w$ intersects first (resp.~last).
For example, we have~$s(x_1^{-1}x_2)=3$ and~$d(x_1^{-1}x_2)=4$ for the decomposed surface
in Figure~\ref{f:decsur}.
We define an order on elements of~$\pi_1(\Sigma)$ of fixed length~$r$ such that 
for~$v,w$ with~$s(v)<s(w)$ we have~$v<w$. For the remaining
pairs~$(v,w)\in\pi_1(\Sigma)^2$ with $v\not=w$ we choose~$m>0$ such that
$t:=v(1)\ldots v(m-1)=w(1)\ldots w(m-1)$ and $x:=v(m)\not=w(m)=:y$
and say that~$v<w$ if

\begin{equation} 
\sigmab(t)+\rho(s(y),s(x))+\rho(d(t),s(x))+\rho(d(t),s(y))\equiv 0\;\mod\;2,
\end{equation}

where $\rho(a,b)$ is~$0$ for~$a>b$ and~$1$ for~$a\leq b$.
Notice that we have
$s(x)\not=s(y)$, $d(t)\not=s(x)$, $d(t)\not=s(y)$. 
See Figure~\ref{f:order} for an example.

\begin{figure}[!h]
\centering
\setbox1=\hbox{\input{tex/conorder}}
$\vcenter{\box1}$
\caption{$x_1x_2^{-1}x_1^{-1}=v>w=x_1x_2^{-1}x_1$}\label{f:order}
\end{figure}

We can visualize the relation
defined above as follows:
the condition $s(v)<s(w)$ says that if
the first interval where~$v$ enters a band lies to the left of the first interval
where~$w$ enters a band, then we have~$v<w$.
For the second condition, represent~$v$ by a generic loop on~$\Sigma$.
We start drawing~$w$ by entering the first band to the right of~$v$ 
and then by running close and parallel to~$v$
through the first~$m-1$ bands. 
Then~$w$ leaves the $m-1$-st band to the left of~$v$ iff~$\sigma(t)=0$.
Three configurations of $(s(x),s(y),d(t))\in\N^3$ and the possible ways how
we can continue drawing~$w$ for~$\sigma(t)=0$ are shown in Figure~\ref{f:draww}.

\begin{figure}[!h]
\centering
\setbox1=\hbox{\input{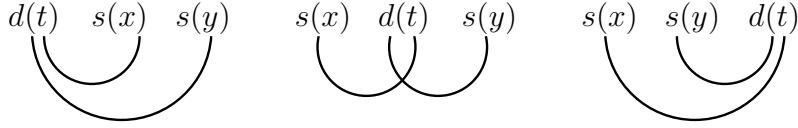}}
$\vcenter{\box1}$
\caption{Configurations of $(s(x),s(y),d(t))$ and intersections}\label{f:draww}
\end{figure}

The remaining possibilities in Figure~\ref{f:draww} for~$\sigma(t)=0$ 
are given by interchanging~$x$ and~$y$. The case~$\sigma(t)=1$ 
is treated in a similar way.
We see that 

$$
\sigma(t)+\rho(s(y),s(x))+\rho(d(t),s(x))+\rho(d(t),s(y)\equiv 1\;\mod\;2
$$ 

iff the
part of~$w$ corresponding to~$y$ going from the interval~$d(t)$ to~$s(y)$ 
must intersect the part of~$v$ corresponding to~$x$ going from~$d(t)$ to $s(x)$.
It is easy to verify that the defined relation is an order relation.

An important property of $v,w\in\pi_1(\Sigma)$ with~$v>w$ is the following.
Assume that~$v$ and~$w$ are represented by 
generic loops on~$\Sigma$ corresponding to reduced words, we have $v(1)=w(1)$, 
$v$ enters the first band at the interval $s(v)$ 
to the left of~$w$, and~$v$ and~$w$ have no intersections 
while they are traversing bands~$B_i$~($i\geq 1$). Choose the smallest number~$n$ such 
that the part
of~$v$ going from the interval~$d(v(n-1))$ to~$s(v(n))$ 
intersects the part of~$w$ going from~$d(w(n-1))$ to~$s(w(n))$ 
in an odd number of points.
Then we necessarily have~$v(i)=w(i)$ for all $i<n$. 

Now 
we can choose representatives~$\hat{w}$ 
of conjugacy classes of elements~$w\in\pi_1(\Sigma)$ as follows:

\begin{eqnarray}
\dot{w} & = & \min\{v\;\vert\; \mbox{$v$ is cyclically reduced, $[v]=[w]$ 
or $[v^{-1}]=[w]$}\},\label{e:dotw}\\
\hat{w} & = & \label{e:hatw}
\left\{\begin{array}{ll}
\dot{w} & \mbox{if $[\dot{w}]=[w]$,}\\
\dot{w}^{-1} & \mbox{if $[\dot{w}^{-1}]=[w]$.}
\end{array}
\right.
\end{eqnarray}

For a decomposed surface we 
choose~$\Conjrep:=\{\hat{w}\;\vert\;w\in\pi_1(\Sigma)\}$. 
For $\Sigma=P^2$ 
we obviously have to choose~$\Conjrep=\pi_1(P^2,*)\cong \Z/2$.
By a descending knot~$K$ we shall mean in this case
a based knot~$K$ with constant height~$p_I(K)$. 
The isotopy class of a descending knot~$K$ is again uniquely determined 
by~$\alpha(K)\in\pi_1(P^2,*)$.

For $i=0,1$ denote by $f_i$  an
orientation preserving 
homeomorphism~$I\cong [i/2, (i+1)/2]$.
By Proposition~\ref{p:emb} we obtain 
isomorphisms of $\Z[z]$-modules

$$
(\id\times f_i)_*:\Conw(\Sigma\times I)
\longrightarrow\Conw(\Sigma\times [i/2,(i+1)/2]).
$$
 
The identity 
$
\Sigma\times[1/2,1]\cup
\Sigma\times [0,1/2]=\Sigma\times I
$ 
induces a $\Z[z]$-bilinear product 

$$
\Conw(\Sigma\times I)\times\Conw(\Sigma\times I)\cong
\Conw(\Sigma\times [1/2,1])\times\Conw(\Sigma\times [0,1/2])
\longrightarrow \Conw(\Sigma\times I),
$$

where we extend the
order of the components of ordered based links~$L_1$, $L_2$ to their product~$L_1L_2$ 
by $K_1<K_2$ for $K_1\subseteq L_1$, $K_2\subseteq L_2$.
The projection $\Conw(\Sigma\times I)\longrightarrow
\Conw_0(\Sigma\times I)$ becomes a homomorphism of rings.

Since the multiplication in $\Conw(\Sigma\times I)$ is in general not commutative
we will use a standard order for products of descending knots.
For this purpose we fix an arbitrary choice of an order on the set~$\Conjrep$.
For example, for a decomposed surface we 
can take the lexicographical order on the alphabet~$x_1,x_1^{-1},x_2,\ldots$ using reduced
expressions.

\begin{defi}\label{d:desclink}
An {$\Conjrep$-descending link} is an ordered based link isotopic to a 
product of descending knots~$K_1\ldots K_n$ with~$\alpha(K_i)\in\Conjrep$ 
and~$\alpha(K_1)\leq \ldots \leq \alpha(K_n)$.
\end{defi}

\section{Main results about the Conway skein module}\label{s:mr}

Let~$\Tor(\Conw(M))=\{m\in\Conw(M)\;\vert\; \exists p\in\cori:pm=0\}$ 
be the torsion submodule of~$\Conw(M)$ and
denote the quotient~$\Conw(M)/\Tor(\Conw(M))$ by~$\Conquo(M)$.
Let~$\Ring(M)$ be the quotient of~$\Conw_0^\#(M)$ by the two-sided ideal generated by
the elements~$\overK_e\overK_w$, where~$e$ is 
the neutral element in~$\pi_1(M)$ and~$w\in\pi_1(M)$ is arbitrary.
For non-orientable~$M$ we have~$\Ring(M)=\Conw_0^\#(M)$ 
because a link~$L$ one of whose components is a trivial knot satisfies~$2L=0\in \Conw(M)$
by the Relation~$(Bas)$ and isotopies. 

We choose~$\Conjrep$ as in Section~\ref{s:surf}.
Recall the definition of~$\theta$ from page~\pageref{d:theta}.
With the notation from above we can state the following theorem.

\begin{theorem}\label{t:conwaypoly}
Let $\Sigma=P^2$ or let $\Sigma$ be a compact connected surface with boundary.
Then there exists a unique $\cori$-linear map \nopagebreak

$$
\nabla:\Conw(\Sigma\times I)\longrightarrow
\cori\otimes_\cor \Ring(\Sigma\times I)
$$\nopagebreak

satisfying $\nabla(L)=1\otimes \overL$ for all $\Conjrep$-descending links~$L$
with~$\theta\left(\overL\right)\not=0$.
The map~$\nabla$ induces an isomorphism
$\Conquo(\Sigma\times I)\cong\cori\otimes_\Z\Ring(\Sigma\times I)$.
\end{theorem}

Theorem~\ref{t:conwaypoly} will be proven in Section~\ref{s:eC}.
In Sections~\ref{s:genc} and~\ref{s:gencx} 
we prove the uniqueness of~$\nabla$. This part of the
proof contains a constructive algorithm for the computation of~$\nabla(L)$. 
In Sections~\ref{s:wC} and~\ref{s:eC} we prove the existence of~$\nabla$.
For this part of the proof we use a universal Vassiliev invariant~$Z_{\Sigma\times I}$.

As a first example, 
observe that for $\Sigma=I^2$ we have $\nabla(L)=\tilde{\nabla}(L)\overK_e$
where $\tilde{\nabla}(L)$ is the usual Conway polynomial of~$L$ (see Equation~(\ref{e:skein})).
 




For a link~$L\subset \Sigma\times I$ we denote the mirror image under reflection in
$\Sigma\times\{1/2\}$ 
by~$L^*$. For a link~$L$ with ordered components, we use the same order for the components 
of~$L^*$ as for their preimages.
For~$i=0,1$ denote by~$\vert L\vert_i$
the number of components~$K$ of~$L$ 
satisfying~$\sigmab(K)=i$. 
Then we have the following theorem about special properties of the Conway skein
module for certain surfaces.

\begin{theorem}\label{t:specase}
(1) For $\Sigma=I^2$, $S^1\times I$, $P^2$
or the M\"obius strip~$X$ the 
map~$\nabla$ satisfies~$\nabla(L)=1\otimes\overL$ 
and~$\nabla(L_1 L_2)=\nabla(L_1)\nabla(L_2)$ for all~$\Conjrep$-descending 
links~$L, L_1, L_2$.

(2) For $\Sigma=P^2$ or $\Sigma=X$ we have
$L=(-1)^{\vert L\vert_0} L^*\in\Conquo(\Sigma\times I).$
\end{theorem}

Theorem~\ref{t:specase} will be proven in Section~\ref{s:eC}.
Notice that we do not require~$\theta(L)\not=0$ in Part~(1) of the theorem.
In the case where~$\Sigma=I^2$ we have the well-known 
symmetry property $L^*=(-1)^{{\vert L\vert}+1} L$ similar but different to the formula in
Part~(2) of Theorem~\ref{t:specase}
whereas for other surfaces no formula of this type is valid. 

%

Let~$L$ be an ordered based link in~$P^2\times I$ represented by a link in 
$X\times I$ as in Figure~\ref{f:Conway}, 
where $X=P^2\setminus D^2$ is the M\"obius strip.
Let~$L'$ be a link
in~$I^2\times I$ 
as in Figure~\ref{f:Conway}, where the parts of the two diagrams 
consisting of a box \labeled{}~$T$ are identical.

\begin{figure}[!h]
\centering
\setbox1=\hbox{\input{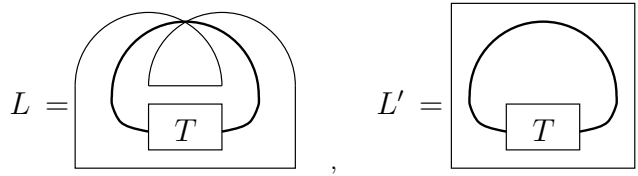}}
$\vcenter{\box1}$
\caption{A link in $P^2\times I$ and a link in $I^2\times I$}
\label{f:Conway}
\end{figure}

If $L$ is a knot, then we have~$\nabla(L)=\tilde{\nabla}(L')\overK_s$
where~$s$ denotes the unique non-trivial
element of~$\pi_1(P^2,*)$. 
In general, we can compose~$\nabla$ with
the invariant~$C(L)$ (see Equation~(\ref{e:defCL})) and obtain

\begin{equation}
\nabla(C(L))=\left\{
\begin{array}{ll}
\{\tilde{\nabla}(L')\overK_s\} & \mbox{if $L$ is a knot,}\\
\{\pm\tilde{\nabla}(L')\overK_s\} & \mbox{otherwise.}
\end{array}
\right.
\end{equation}

The example above is related to a simple special case of the following theorem.

\begin{theorem}\label{t:Csquare}
Let $K$ be a knot in~$P^2\times I$ with~$\nabla(K)=f(z)\overK_s$. 
Then $$\nabla(p^{-1}(K))=f(z)^2\overK_e,$$
where $p:S^2\times I\longrightarrow P^2\times I$ denotes the $2$-fold covering map.
\end{theorem}

Theorem~\ref{t:Csquare} will be proven in Section~\ref{s:squareproof} by using relations
between Vassiliev invariants of~$K$ and~$p^{-1}(K)$.

Recall that a knot~$K$ in $S^3$ is called strongly amphicheiral, if there exists an
orientation-reversing involution~$i$ of~$S^3$ with~$i(K)=K$. 
One says that~$K$ is strongly positive amphicheiral, if there exists an involution~$i$ 
as above that preserves the orientation of~$K$.
We have the following corollary.

\begin{coro}\label{c:Csquare}
Let $K$ be a strongly positive amphicheiral knot in~$S^3$. Then~$\tilde{\nabla}(K)=f(z)^2$
for some $f\in\Z[z]$.
\end{coro}
\begin{proof}
Strongly positive amphicheiral knots in~$S^3$ are in one-to-one correspondence with 
knots~$K$ in~$P^2\times I$ with~$\sigma(K)=1$.
To see this, recall that 
the set of fixed points of an orientation-reversing involution~$i$ of~$S^3$ is 
either~$S^0$ or~$S^2$ (see~\cite{Smi}). For strongly positive
amphicheiral knots~$K$ it is sufficient to consider
involutions~$i$ with fixed point set~$S^0\subset S^3$ and~$K\cap S^0=\emptyset$. 
By a theorem of Livesay (see~\cite{Liv}, \cite{Rub}), we then have
$(S^3\setminus S^0)/ i\cong P^2\times ]0,1[$.
The corollary follows by applying Theorem~\ref{t:Csquare} to the
image of~$K$ in~$(S^3\setminus S^0)/ i$.
\end{proof}

Corollary~\ref{c:Csquare} was proven directly by Kawauchi and Hartley
using the Blanchfield pairing 
(see~\cite{KaH}, \cite{Bla}, \cite{Kaw}).
Theorem~\ref{t:Csquare} contains the stronger statement 
that~$f(z)$ itself can be calculated using
a skein relation. 
It would be interesting
to have a direct skein-theoretical proof of Theorem~\ref{t:Csquare}.

\section{Generators of $\Conw(\Sigma\times I)$}\label{s:genc}

In general one can show that all sets of knots~$\cal K$ for 
which~$\{\alpha(K)\;\vert\; K\in{\cal K}\}$ 
contains a set of representatives of conjugacy classes
of~$\pi_1(\Sigma)$ generate a dense subalgebra of the $z$-adic 
completion~$\varprojlim\,\Conw(\Sigma\times I)/z^n\Conw(\Sigma\times I)$
of~$\Conw(\Sigma\times I)$.
As we will see by the following example,
the set~$\cal K$ does not necessarily generate
the~$\Z[z]$-module~$\Conw(\Sigma\times I)$.
Choose a generator $s\in\pi_1(S^1\times I,*)$.
Denote the descending knot in~$S^1\times I^2$ belonging to~$s^n$ by~$K'_n$.
There are knots~$K_n''$ in~$S^1\times I^2$ that are a connected sum of~$K_n'$ and a Whitehead
knot in~$S^1\times I^2$ such that the equation
$K_n''=K'_n-zK'_{\sgn(n)(\vert n\vert +1)}K'_{-\sgn(n)}$ 
holds in~$\Conw(S^1\times I^2)$ (see Figure~\ref{f:whitehead} for~$n=3$).
Using Theorem~\ref{t:conwaypoly}
it is easy to see that the knots $K''_n$ ($n\in\Z$) 
do not generate~$\Conw(S^1\times I^2)$ as an algebra, 
but~$\{\alpha(K_n'')\;\vert\;n\in\Z\}=\pi_1(S^1\times I)$.

\begin{figure}[!h]
\centering
\setbox1=\hbox{\input{tex/contex}}
$\vcenter{\box1}$

\caption{$K_3''=K'_{3}-zK'_{4}K'_{-1}\in\Conw(\Sigma\times I$)}\label{f:whitehead}
\end{figure}

The knots $K_n''$ from the previous example are not descending.
For a set~$\Conjrep'$ consisting of cyclically reduced representatives of conjugacy classes
in~$\pi_1(\Sigma)$, where~$\Sigma$ is a decomposed surface, one can show 
by using Theorem~\ref{t:conwaypoly}
that the descending knots~$K_w$ ($w\in\Conjrep'$) 
generate~$\Q(z)\otimes_{\Z[z]}\Conw(\Sigma\times I)$ as a~$\Q(z)$-algebra, where~$\Q(z)$
is the quotient field of~$\Z[z]$. 
The knots $K_w$ ($w\in\Conjrep'$) do not necessarily generate $\Conw(\Sigma\times I)$ as
a~$\Z[z]$-algebra.
For example, let~$\Sigma$ be a 
disc with three holes 
decomposed such that $s(x_i)=i$, $d(x_i)=7-i$ ($i=1,2,3$). 
By  
Theorem~\ref{t:conwaypoly} and Proposition~\ref{p:emb}
there exists a basis of the~$\Z[z]$-module~$\Conquo(\Sigma\times I)$
containing the element~$K_{x_1{x_3^{-1}}x_2}^*
-K_{x_1x_2{x_3^{-1}}}^*$, where~$K^*$ denotes the mirror image of~$K$ with respect to
$\Sigma\times\{1/2\}$.
Let $\{x_2x_1x_3^{-1}, x_3^{-1}x_1x_2\}\subset \Conjrep'$.
If the knots~$K_w$ ($w\in\Conjrep'$) would generate~$\Conw(\Sigma\times I)$ 
we see by Theorem~\ref{t:conwaypoly} that we
could find a basis of $\Conquo(\Sigma\times I)$ containing the element
$K_{x_2x_1{x_3^{-1}}}-K_{{x_3^{-1}}x_1x_2}$.
But this is impossible because we can show by a computation that 

\begin{equation}\label{e:nogen}
K_{x_2x_1{x_3^{-1}}}-K_{{x_3^{-1}}x_1x_2}=(1+z^2)\left(K_{x_1{x_3^{-1}}x_2}^*
-K_{x_1x_2{x_3^{-1}}}^*\right)\in\Conw(\Sigma\times I).
\end{equation}

The following lemma says that the set~$\Conjrep$ chosen in Section~\ref{s:surf} has 
better properties.

\begin{lemma}\label{l:descgen}
The descending knots $K_w$ ($w\in\Conjrep$) generate
$\Conw(\Sigma\times I)$ as a $\cori$-algebra.
\end{lemma}
\begin{proof}
Let~$L$ be a diagram of an ordered based link. For the proof we use the following strategy:
we will never increase the number of crossings of~$L$, 
and make computations modulo diagrams with fewer crossings. 
This allows us to make crossing changes.
We will prove the lemma by induction, 
where the induction base is given by the following arguments  
for link diagrams without crossings, and the induction step is given by the same
arguments for link diagrams with crossings.

{\em Case $\partial\Sigma\not=\emptyset$, Step 1:} 
We apply an isotopy that does not increase the number of double points in 
$p_\Sigma(L)$
such that the projection of~$L$ to $\bigcup_{i=1}^k B_i$ consists of
parallel strands connecting points in $I_{s(x_i)}$ with points in $I_{d(x_i)}$.
Then we pass from $L$ to a diagram of a product of knots by crossing changes.

{\em Case $\partial\Sigma\not=\emptyset$, Step 2:} 
Each component
of~$L$ represents a unique word in the generators $x_i$ of $\pi_1(\Sigma)$.  
We claim that 
we can replace each of these components by a diagram of a 
based knot~$K$ representing a cyclically reduced word in~$\pi_1(\Sigma)$:
for words that are not cyclically reduced we
find a segment~$h$ of~$K$ corresponding to a cancellation such that~$p_\Sigma(h)$ connects two 
points~$P_1, P_2\in I_\nu$ for some~$\nu$ 
and all the parts of~$p_\Sigma(K)$ that go 
into~$B_0$ between~$P_1$ and~$P_2$ intersect~$p_\Sigma(h)$ 
(see Figure~\ref{f:cancellation}).

\begin{figure}[!h]
\centering
\setbox1=\hbox{\input{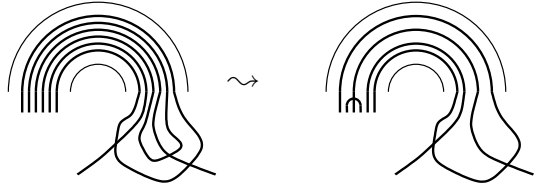}}
$\vcenter{\box1}$
\caption{Eliminating cancellations of generators of $\pi_1(\Sigma)$}\label{f:cancellation}
\end{figure}

The part~$h$ can first be moved to a suitable height by making crossing changes 
and then be pulled back by applying an isotopy
that does not increase the number of crossings as shown in Figure~\ref{f:cancellation}.
This will remove one cancellation from the word in the generators~$x_i, x_i^{-1}$
corresponding to the knot projection.
We continue by induction on the number of cancellations to pass to a cyclically reduced 
word.

{\em Case $\partial\Sigma\not=\emptyset$, Step 3:} 
Recall the definitions of~$\dot{w}$ and~$\hat{w}$ for~$w\in\pi_1(\Sigma)$ 
(Equations~(\ref{e:dotw}) and~(\ref{e:hatw})).
We continue with the modification of each based knot~$K$ in~$L$ by moving the basepoint
of~$K$ until $w=\hat{w}$ where $w=\alpha(K)$. 
By the Relation~$(Bas)$ this only contributes a factor~$\pm 1$.
We can further achieve that the projection~$A$ of the basepoint of~$K$ to~$\Sigma$ lies 
in~$I_{s(\dot{w})}$. 
Assume that there are points in~$p_\Sigma(K)\cap I_{s(\dot{w})}$ to the left 
of~$A$ in our pictures of~$\Sigma$.
Let~$\pe$ be the left neighbor of~$A$ in~$p_\Sigma(K)\cap I_{s(\dot{w})}$.
Notice that there may be points of~$p_{\Sigma}(L\setminus K)\cap I_{s(\dot{w})}$ 
between~$\pe$ and~$A$.
Let~$v$ be the element of~$\pi_1(\Sigma)$ that is represented by~$p_\Sigma(K)$ 
with basepoint~$\pe$ and with an orientation such that at~$\pe$ the curve~$p_\Sigma(K)$ 
enters the band. 
By the definition of the representative~$\hat{w}$ we have~$v\geq \dot{w}$.
Assume first~$v>\dot{w}$. We have~$v(1)=\dot{w}(1)$. 
Therefore, following the two strands of~$p_\Sigma(K)$ 
starting by entering a band at~$\pe$ and~$A$, 
we follow paths~$p$ and~$a$ respectively 
that pass through the same bands until we find a crossing~$c$
between them. 
Denote that part of~$L$ by~$T$ whose projection lies in a small neighborhood of the triangle
bounded by~$a$,~$p$, and~$I_{s(\dot{w})}$ with corners~$\pe$,~$A$, and~$c$.
By crossing changes we pass to a suitable height on~$T$ and then pull
back this part by an isotopy as shown in Figure~\ref{f:pullcback} until the points~$\pe$ 
and~$A$ are interchanged.
The projection of~$T$ 
arrives very close to~$I_{s(\dot{w})}$ in~$B_0$ such that no crossings with the 
remaining part of~$p_\Sigma(L)$ can appear.

\begin{figure}[!h]
\centering
\setbox1=\hbox{\input{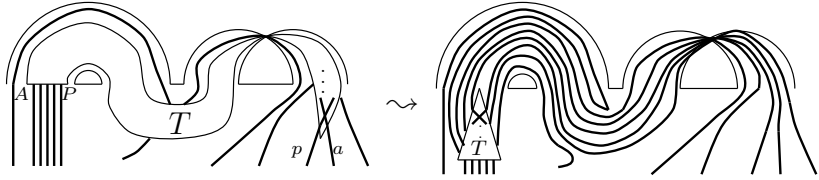}}
$\vcenter{\box1}$
\caption{Moving the basepoint to the left 
in~$p_\Sigma(K)\cap I_{s(\dot{w})}$}\label{f:pullcback}
\end{figure}

If $v=\dot{w}$, then we push the basepoint along~$K$ until its projection arrives at~$\pe$.
Therefore in any case this modification
does not increase the number of crossings of~$K$.
This implies that by induction we can assume that there are no points 
in~$p_\Sigma(K)\cap I_{s(\dot{w})}$ to the left of~$A$.
Then we pass from~$K$ to a descending knot by crossing changes. 
Since this works 
for all components of~$L$ we have proven the lemma in the 
case where~$\partial\Sigma\not=\emptyset$.

\smallskip

{\em Case $\Sigma=P^2$:} We choose an open disc $D\subset P^2$. Then
$P^2\setminus D=X$ is the M\"obius strip. 
It is easy to see that every link in~$P^2\times I$ can be represented by a link 
in~$X\times I$ and that a descending link in~$X\times I$ is isotopic to a descending
link in~$P^2\times I$. Therefore we conclude by using this lemma for~$\Sigma=X$.
\end{proof}

\medskip

Let $\Sigma$ be a surface with nonempty boundary.
Let $w\in\pi_1(\Sigma)\setminus\{e\}$.
In the free group~$\pi_1(\Sigma)$ the normalizer~$N(w)$
is isomorphic to~$\Z$. A generator~$z$ of~$N(w)$ is characterized by the property
that $w=z^k$ where~$k\in\Z$ and $\vert k\vert$ is as large as possible. 
If~$\sigma(z)=1$ and~$\sigmab(w)=0$, 
then by Proposition~\ref{p:ringiso} we have
$\overline{K}_w=0\in \Conw_0^\#(\Sigma\times I)$.
Nevertheless, we may have~$K_w\not=0\in\Conquo(\Sigma\times I)$. 
Hence it is not obvious how to calculate~$\nabla(K_w)$ by using only
the properties of~$\nabla$ given in Theorem~\ref{t:conwaypoly}.
Our next goal will be to see why the equation in Figure~\ref{f:desccable}
implies that~$\nabla(K_{(x_1x_2^{-1})^2})=
- z\oK_{x_1x_2^{-2}}\oK_{x_1}$.
In the rest of this
section we will reduce this question to 
a problem in~$\Conquo(X\times I)$, where~$X$ is the
M\"obius strip. 

\begin{figure}[!h]
\centering
\setbox1=\hbox{\input{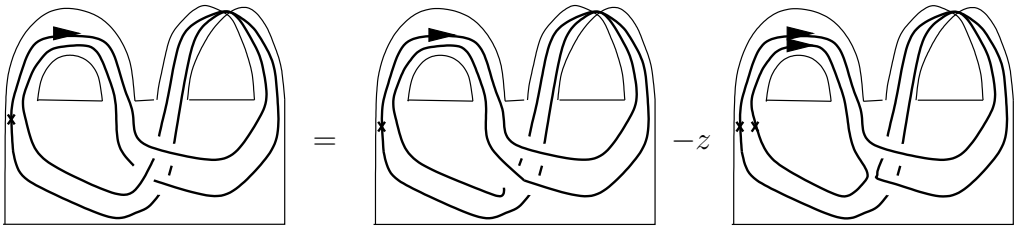}}
$\vcenter{\box1}$
\caption{What is $\nabla(K_{(x_1x_2^{-1})^2})$ ?}\label{f:desccable}
\end{figure}

Fix a generator~$s$ of~$\pi_1(X)$ and denote the descending knot~$K_{s^n}$ by~$K_n$.
For~$v\in\pi_1(\Sigma)$ with~$\sigmab(v)=1$ we have an embedding
$i_v:X\times I\longrightarrow \Sigma\times I$ such that~$i_v(K_1)=K_v$. 
Up to isotopy~$i_v$ is uniquely
determined by~$v$ (since~$\sigma(v)=1$ no framing is needed). 
By Proposition~\ref{p:emb} the embedding~$i_v$ induces a map
between skein-modules
${i_v}_*:\Conw(X\times I)\longrightarrow \Conw(\Sigma\times I)$. We can assume that~${i_v}_*$ 
is $\cori$-linear by making a suitable choice of the 
oriented subset~$U$ of~$X\times I$.
For~$w\in \pi_1(\Sigma)\setminus\{e\}$ there exists a unique generator~$u$ 
of~$N(w)$ such that~$w=u^n$ for~$n\in\N$. If~$\sigmab(u)=1$ and~$\sigmab(w)=0$, then
we define 
the {\em cabled descending knot} as~$\hat{K}_w=(i_u)_*(K_n)$. 
We define~$\hat{K}_w=K_w$ in the remaining cases.
Using this notation we can rewrite the equation in Figure~\ref{f:desccable} 
as~$K_{(x_1x_2^{-1})^2}=\hat{K}_{(x_1x_2^{-1})^2}-z K_{x_1x_2^{-2}} K_{x_1}
=\hat{K}_{(x_1x_2^{-1})^2}-z \hat{K}_{x_1x_2^{-2}} \hat{K}_{x_1}$.

\begin{lemma}\label{l:cabgen}
The 
knots $\hat{K}_w$ ($w\in\Conjrep$) generate
$\Conw(\Sigma\times I)$ as a $\cori$-algebra.
\end{lemma}
\begin{proof}
For $\Sigma=P^2$ there is nothing to prove. Assume that~$\partial\Sigma\not=\emptyset$.
Let~$L$ be a link in~$\Sigma\times I$.
By the first two steps of the proof of 
Lemma~\ref{l:descgen} we can assume that the connected components of~$L$
represent cyclically reduced words in~$\pi_1(\Sigma)$.
Applying the following Steps~1 and~2 to all components~$K$ of~$L$ will imply the
lemma by induction on the number of crossings.

\smallskip

{\em Step 1 (''Collect all starting points on the left''):}
Choose a maximal~$k>0$ such that~$\dot{w}=(a_1\ldots a_n)^k\in\pi_1(\Sigma)$ 
where~$w=\alpha(K)$
and~$a_i\in\{x_j^{\pm 1}\}$. 
Assume that $\sigmab(a_1\ldots a_n)=1$ (otherwise there is nothing to do).
Consider the~$k$ points
$A_i$ on~$K$ such that~$p_\Sigma(A_i)\in I_{s(a_1)}$ and~$p_\Sigma(K)$ with each of these 
points 
represents~$\hat{w}$. We want to move all the points~$p_\Sigma(A_i)$ to the left 
side in~$p_\Sigma(K)\cap I_{s(a_1)}$ 
without increasing the number of crossings of~$L$: if this is not already the case, then 
choose~$j$ such that the point~$P$ to the left of~$A:=p_\Sigma(A_j)$ 
in~$p_\Sigma(K)\cap I_{s(a_1)}$ is not one of the 
points~$p_\Sigma(A_i)$. Let~$(b_1\ldots b_n)^k$ be the element 
of~$\pi_1(\Sigma)$ represented by~$p_\Sigma(K)$ 
with respect to~$P$ and oriented such that it leaves~$B_0$ at~$P$.
Since we use this orientation of~$p_\Sigma(K)$ the word~$b_1\ldots b_n$ is obtained
from~$a_1\ldots a_n$ by a cyclic permutation.

The minimality of~$\dot{w}$ and the maximality of~$k$ implies 
$a_1\ldots a_n<b_1\ldots b_n$.
We have~$b_1=a_1$. Therefore, we can follow $p_\Sigma(K)$ on two strands along
$b_1\ldots b_m=a_1\ldots a_m$ ($m<n$) starting at~$P$ and~$A$ until we find a 
crossing~$c$
between these two strands. As in Step 3 of the proof of Lemma~\ref{l:descgen} we make crossing
changes and then pull
back this crossing 
and possibly some other parts of~$L\setminus K$ 
by an isotopy until the points~$P$ and~$A$ are interchanged (see Figure~\ref{f:pullcback}).
With this modification we do not increase the number of crossings of~$L$ and the 
point~$A$ 
has moved one step to the left in~$p_\Sigma(K)\cap I_{s(a_1)}$.
We will now verify that no other point~$A'=p_\Sigma(A_{j'})$ 
has moved to the right during this operation: 
assume that we reach~$A'$ after starting at~$P$ and traveling
along~$b_1\ldots b_{m'}$ with $0<m'\leq m$. 
Then for $x=a_1\ldots a_{m'}$ we either have~$a_1\ldots a_n=xyx$ or~$a_1\ldots a_n=xyx^{-1}$.
We can exclude the second case because this would imply~$x=a_1\ldots a_{m'}=b_1\ldots b_{m'}
=x^{-1}$ which is impossible. By the minimality of~$\dot{w}=(xyx)^k$ 
and the maximality of~$k$ we have~$xyx<yxx$ and~$xyx<xxy$.
From~$xyx<yxx$ we deduce~$xy<yx$. Then~$xyx<xxy$ and~$yx>xy$ implies~$\sigma(x)=1$.
Therefore, starting at~$P$ and~$A$
and traveling along~$x$, the first strand reaches~$I_{s(a_1)}$ to the right of the second one.
This means that by pulling back the crossing~$c$ we move the two points~$A$ and~$A'$ to 
the left
(see Figure~\ref{f:2basemove} for an example).
We can continue by induction until all points~$p_\Sigma(A_i)$ 
are on the left side of~$p_\Sigma(K)\cap I_{s(a_1)}$.

\begin{figure}[!h]
\centering
\setbox1=\hbox{\input{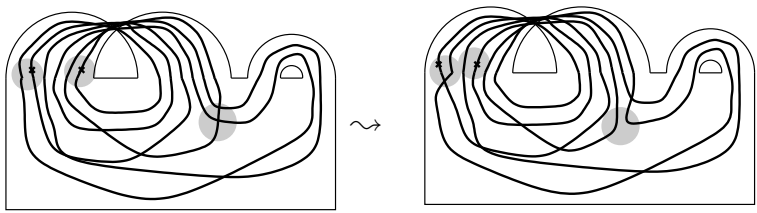}}
$\vcenter{\box1}$
\caption{Applying step~1 to the leftmost point of $p_\Sigma(K)\cap I_1$}\label{f:2basemove}
\end{figure}

\smallskip

{\em Step 2 (''move component into cabling position''):}
Successively for~$i=2,\ldots ,n$, we consider all parts~$s_\nu$ of~$K$ 
such that the projections~$p_\Sigma(s_\nu)$ connects
some interval~$I_j$ with~$I_{s(a_i)}$ inside of~$B_0$.
By crossing changes we move all the~$s_\nu$ to a suitable height and then
push the part~$T$ containing all crossings between these strands 
across the band at~$I_{s(a_i)}$. This isotopy can be chosen 
in the following way such that the
total number of crossings of~$L$ does not increase:
after being pushed across the band 
the part~$T$ arrives very close to~$I_{d(a_i)}\times I$
and the strands that now enter into~$I_{s(a_i)}$ without intersections 
in their projection to~$\Sigma$ can
be moved into a position such that they coincide 
with original parts of~$T$ except in 
neighborhoods of the crossings and except for orientations.

After having done the previous modification 
for $i=2,\ldots,n$ we make crossing changes
such that $K=i_{{a_1\ldots a_n}_*}(K')$, where $K'$ is a knot in $X\times I$.
Applying
the arguments in the proof of Lemma~\ref{l:descgen} for~$\Sigma=X$ to~$K'$ we see that~$K$ 
is equal to $i_{{a_1\ldots a_n}_*}(K_k)=\hat{K}_w$
modulo diagrams with fewer crossings.
\end{proof}
 
\medskip

In a diagram of a product of knots we can permute two components modulo diagrams with fewer
crossings. By Relation~$(Ord)$ this allows to express~$2K_w^2$ ($\sigma(w)=1$)
by diagrams with fewer crossings. 
This can be used in the calculation of~$\nabla(L)$ 
by working
with~$\Z[1/2][z]\otimes\Ring(\Sigma\times I)$. In order to prove that~$\nabla$ takes
values in~$\Z[z]\otimes\Ring(\Sigma\times I)$ we have to be more careful.
Since the presentation of~$\Conw_0(\Sigma\times I)$ by generators and relations 
(see the proof of Proposition~\ref{p:ringiso}) shows 
that~$K_w^2\not\equiv 0 \;\mod\;z$ we pass to a quotient 
of~$\Conw(\Sigma\times I)$.
Let~$J$ be the two-sided ideal of~$\Conw(\Sigma\times I)$ generated by
elements~$i_*(K_1^2)$ for all embeddings~$i:X\times I\longrightarrow\Sigma\times I$
as in Propositon~\ref{p:emb}. Then we have the following lemma.

\begin{lemma}\label{l:cabmodgen}
The $\cori$-module $\Conw(\Sigma\times I)/J$ 
is generated by the following set of links:\nopagebreak

$$
\left\{
\hat{K}_{w_1}\ldots \hat{K}_{w_n}\;\vert\;
n\in\N, w_i\in \Conjrep, w_i\leq w_{i+1},
(w_i=w_{i+1}\implies \sigmab(w_i)=0)
\right\}.
$$
\end{lemma}
\begin{proof}
The proof of Lemma~\ref{l:cabgen} shows
that the~$\cori$-module~$\Conw(\Sigma\times I)$ is generated by the 
links~$\hat{K}_{w_1}\ldots \hat{K}_{w_n}$ with~$n\in\N$, $w_i\in \Conjrep$
and~$w_i\leq w_{i+1}$.
Consider links~$L$ of the 
form~$L' K' K'' L''$ where~$L'$ and~$L''$ are ordered based links
and~$K'=K''=\hat{K}_w$ with~$\sigma(w)=1$.
For~$\sigma(w)=1$ we always have~$K_w=\hat{K}_w$.
As shown in the proof of Lemma~\ref{l:descgen} we may assume that in
a diagram of~$L$ the projections of the components~$K'$ and~$K''$ 
of~$L$ represent cyclically reduced words in~$\pi_1(\Sigma)$.
Furthermore, we may assume that the 
projections~$A'$ and~$A''$ of the basepoints of~$K'$
and~$K''$ lie in~$I_{s(\dot{w})}$.
When we follow the two strands of~$p_\Sigma(K'\cup K'')$ 
starting by entering a band at~$A'$ and~$A''$, 
we follow paths~$a'$ and~$a''$ respectively 
that pass through the same bands until we find a crossing~$c$
between them. We pull this crossing back as shown in Figure~\ref{f:pullcback}. 
We continue with the modification of~$L$
as in Step~2 of the proof of Lemma~\ref{l:cabgen}. We obtain that~$L$ is equal 
to~$L'i_*(K_1^2)L''$ modulo diagrams with fewer crossings, where 
we have~$\alpha(i_*(K_1))=w$ and where
we can
choose the embedding~$i:X\times I\longrightarrow \Sigma\times I$ such 
that~$i_*(K_1)$ is a knot descending with respect to its basepoint (but in general not a 
descending knot). The link~$L'i_*(K_1^2)L''$ lies in the ideal~$J$.
We obtain the lemma by induction on the number of crossings.
\end{proof}

\section{Generators of $\Conquo(X\times I)$}\label{s:gencx}

In this section
we will show by making computations in~$\Conquo(X\times I)$
that the knots~$\hat{K}_w$ with~$\hat{K}_w\not=K_w$
are~$0$ in~$\Conquo(\Sigma\times I)$. 
Verifying directly that

\begin{equation}\label{e:simpletorsion}
(z^2+ 4)i_*(K_1^2)=0\in\Conw(\Sigma\times I)
\end{equation}

gives a good
impression of the ideas used in this section. Readers that are mainly interested in the
case~$\Sigma=P^2$ can concentrate on the verification of Equation~(\ref{e:simpletorsion}).

Recall that we denote the descending knot~$K_{s^n}$ for a fixed choice of a 
generator~$s\in\pi_1(X)$ by~$K_n$.
In this section we represent ordered based links 
in~$X\times I$ by drawing only a part of them as shown in 
Figure~\ref{f:tanclos}. We say that we represent an ordered based 
link in~$L$ as the closure of a tangle.

\begin{figure}[!h]
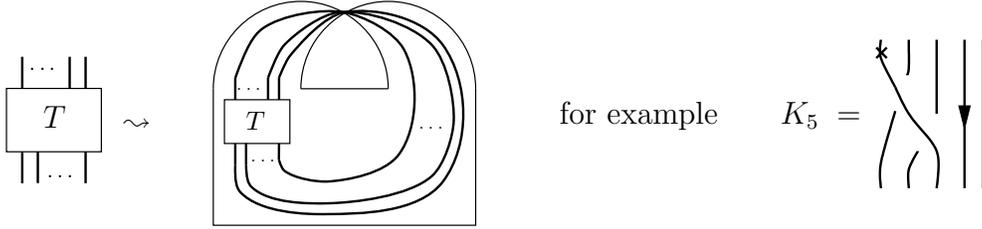

\centering
\setbox1=\hbox{\input{tex/clotan}}
$\vcenter{\box1}$
\qquad\mbox{for example}\qquad
\mbox{$K_5\;=$}\setbox1=\hbox{\input{tex/exspet3}}
$\vcenter{\box1}$
\caption{The closure of a tangle}\label{f:tanclos}
\end{figure}

The tangle~$T$ must have the same number~$n$ of endpoints at the top and 
at the bottom and the
strand at the $i$-th endpoint at the top 
(always by counting the endpoints from the left to the right side)
is directed downwards if and only if the $(n+1-i)$-th
strand at the bottom is directed downwards ($i=1,\ldots,n$). The ordered set of 
basepoints of~$T$ 
is in bijection with the components of the closure of~$T$.

The following lemmas are important for determining~$\Conquo(X\times I)$.
They also provide shortcuts in the computation of the Conway polynomial.

\begin{lemma}\label{l:K2k0}
For all $k,n,m\in\Z$ the following identities hold in $\Conquo(X\times I)$: 

$$
K_{2k}=0\qquad\mbox{and}\qquad K_{2n+1}K_{2m+1}=-K_{2m+1}K_{2n+1}.
$$
\end{lemma}

As a preparation for the proof of Lemma~\ref{l:K2k0} we first prove the following lemma.

\begin{lemma}\label{l:shortcuts}
Assume that $K_{2n}=0\in\Conquo(X\times I)$ for all~$n$ with $-k<n< k$. 
Then for~$i\in\Z$ with $-k<i<k$
the following
identities hold in $\Conquo(X\times I)$:

(1) If a link $L$ is a closure of a tangle~$T$ on~$2i$ strands as in 
Figure~\ref{f:stanclos} with $r=i$, then $L=0$.

(2) If a knot $K$ is a closure of a tangle on $l\in\{2i+1,2i+2\}$
strands as in Figure~\ref{f:stanclos} with $r=i+1$, $K$ is descending with respect 
to its basepoint, the basepoint lies on one of the first~$i+1$ strands of~$K$,
and~$K$ is homotopic to~$K_m$, then~$K=K_m=K_m^*$.

\begin{figure}[!h]
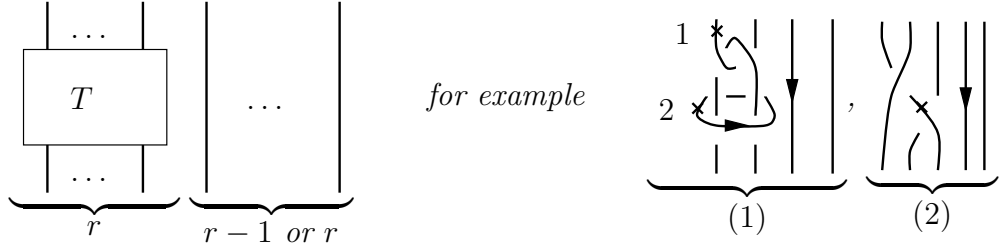

\centering
$\underbrace{\mbox{\setbox1=\hbox{\input{tex/spetan}}
$\vcenter{\box1}$}}_{\mbox{$r$}}$
$\underbrace{\mbox{\setbox1=\hbox{\input{tex/spetan2}}
$\vcenter{\box1}$}}_{\mbox{$r-1$ or $r$}}$
\qquad\mbox{for example}\qquad
$\underbrace{\mbox{\setbox1=\hbox{\input{tex/exspet2}}
$\vcenter{\box1}$}}_{\mbox{$(1)$}}$,
$\underbrace{\mbox{\setbox1=\hbox{\input{tex/exspet}}
$\vcenter{\box1}$}}_{\mbox{$(2)$}}$ 
\caption{Special tangles on an even or odd number of strands}\label{f:stanclos}
\end{figure}
\end{lemma}
\begin{proof}[Proof of Lemma~\ref{l:shortcuts}]
(1) Since $2i$ is even and the diagram of~$L$ is as shown in 
Figure~\ref{f:stanclos} we can pass 
from~$L$ to a product $\pm K_{2i_1}\ldots K_{2i_r}$ with~$-i<i_\nu<i$
by crossing changes and by Relation~$(Bas)$.
If a crossing
of~$L$ is spliced, then we obtain again a closure of a tangle as in Figure~\ref{f:stanclos}.
By induction on the number of crossings of~$L$ we can expand~$L$ 
as a linear combination over~$\cori$ of monomials~$K_{2j_1}\ldots K_{2j_m}$ 
with~$-i<j_i<i$.
Under the hypotheses of the lemma all these monomials 
equal~$0$ in $\Conquo(X\times I)$.

(2) Let $K$ be as in Part~($2$) of the lemma.
When we push the basepoint of~$K$ along the generator~$s$ of~$\pi_1(X)$ and~$l$ is odd, then
Relation~$(Bas)$ contributes the sign~$(-1)^{1\cdot(l-1)}=1$. 
If~$l$ is even, then we will see in the following that we 
will push the basepoint along~$s$ an even 
number of times, because the basepoint is on one of the first~$i+1$ strands of~$K$. This will
not give a sign contribution either.
If the basepoint
passes through a crossing, then by changing this crossing 
the knot will become descending with respect to the new basepoint.
If the crossing is spliced, then we obtain a product of two
knots that are descending with respect to their basepoints. One of these
knots is isotopic to~$K_0=0$ or satisfies the conditions from Part~(1) of the lemma. 
Therefore 
we do not change~$K$ if we push the basepoint along~$K$ and change crossings until
the basepoint is close to the left boundary of~$X$ in our pictures. This implies 
$K=K_m$ if $K$ is homotopic to $K_m$ (the number~$m$ 
satisfies $-l\leq m\leq l$ and $l\equiv m\;\mod\;2$). 
There exists a diagram of $K_{m}$ as a closure of
a tangle as in Figure~\ref{f:stanclos} on $m$ strands having $[(m-1)/2]$~crossings, 
such that by changing arbitrary crossings and by 
splicing one arbitrary crossing we obtain a product of two knots (see 
Figure~\ref{f:tanclos} for $m=5$). 
Using this we
obtain $K_{m}=K_{m}^*$ by similar arguments as above.
\end{proof}

\begin{proof}[Proof of Lemma~\ref{l:K2k0}]
We prove the lemma by induction. We have $K_0=0$.
Assume that the lemma is true for all $k, n, m$ with~$\vert k\vert \leq \ell$ and 
$\vert n\vert+\vert m\vert< \ell$.
We have to prove the lemma for $\vert k\vert -1 =\ell=\vert n\vert +\vert m\vert$.
Since inversion of the orientations of all components of a link induces a 
$\cori$-linear involution of~$\Conw(M)$, we can assume without loss of generality
that~$k>0$ and~$n\geq 0$. 
We will prove the following two equations:

\begin{eqnarray}\label{e:u2n}
2K_{2k} & = & -z\sum_{i=0}^{k-1} K_{2i+1}K_{2(k-i)-1},\label{e:even}\\
K_{2n+1}K_{2m+1}+K_{2m+1}K_{2n+1} & = & zp(z)K_{2k}\qquad\mbox{for 
some $p(z)\in\Z[z]$.}\label{e:odd}
\end{eqnarray}

Equations~(\ref{e:even}) and~(\ref{e:odd}) imply that $(2+z^2q(z))K_{2k}=0$ for some 
$q(z)\in\Z[z]$. Therefore $K_{2k}=0\in\Conquo(X\times I)$. 
This together with Equation~(\ref{e:odd}) implies also
\mbox{$K_{2n+1}K_{2m+1}=-K_{2m+1}K_{2n+1}
\in\Conquo(X\times I)$}.

\smallskip

We prove Equation~(\ref{e:even}) for $k=3$ in Figure~\ref{f:bppusheven}
and explain afterwards why this equation holds for arbitrary~$k$. 

\begin{figure}[!h]
\centering
\setbox1=\hbox{\input{tex/K2kg0}}
$\vcenter{\box1}$
\caption{A refinement of $2 \overK_{2k}=0\in\Conw_0(\Sigma\times I)$}\label{f:bppusheven}
\end{figure}

The first equality in Figure~\ref{f:bppusheven} follows in general from Part~(2)
of Lemma~\ref{l:shortcuts} for~$l=2k$ and from Relation~$(Bas)$.
The new knot is descending with respect to its basepoint.
The second equality follows by applying the skein relation successively~$k$ times to crossings
that appear when the new basepoint is pulled to the left.
When all crossings are changed, then we obtain again~$K_{2k}$.
When one of the~$k$ crossings is spliced, then we obtain a product of two 
descending knots 
each of which satisfies the conditions of Part~(2) of Lemma~\ref{l:shortcuts}
for odd~$l$. This gives us the~$k$ product links $K_iK_{2k-i}$ where~$i$ runs 
successively over the numbers~$1,3,\ldots,2k-1$. Looking at signs 
we obtain Equation~(\ref{e:even}).

\smallskip

Proof of Equation~(\ref{e:odd}):
The link $K_{2n+1}K_{2m+1}$ is the closure of a tangle as shown in 
Figure~\ref{f:ko1ko2}. In this diagram the crossings between the two components of
$K_{2n+1}K_{2m+1}$ are numbered from~$1$ to~$x:=\vert 2m+1\vert$. The position of
the basepoints on~$K_{2n+1}K_{2m+1}$ is not important.

\begin{figure}[!h]
\centering
\setbox1=\hbox{\input{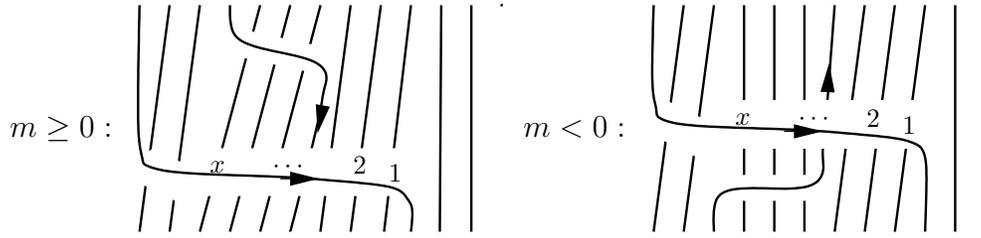}}
$\vcenter{\box1}$
\caption{The links $K_{5}K_{7}$ and $K_5K_{-7}$}\label{f:ko1ko2}
\end{figure}

Let~$D_i(n,m)$ ($i=1,\ldots,x=\vert 2m+1\vert$) be the 
knot obtained by changing the crossings $1,2,\ldots,i-1$ and by 
splicing the crossing~$i$ in Figure~\ref{f:ko1ko2} (see Figure~\ref{f:Didiag}). 
For later use some crossings of the diagram of~$D_i(n,m)$ are marked by the symbol~$*$
in Figure~\ref{f:Didiag}.

\begin{figure}[!h]
\centering
\setbox1=\hbox{\input{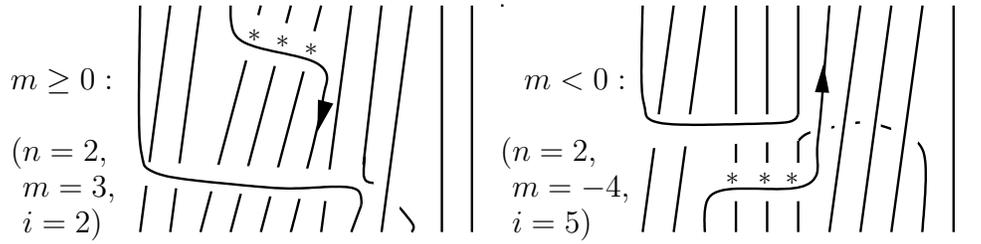}}
$\vcenter{\box1}$
\caption{The diagram $D_i(n,m)$}\label{f:Didiag}
\end{figure}

Since we can permute the two components of~$K_{2n+1}K_{2m+1}$ by changing all
crossings between them we obtain by
the skein-relation and by Relation~$(Ord)$ that

\begin{equation}\label{e:disum}
K_{2n+1}K_{2m+1}=-K_{2m+1}K_{2n+1}+z\sum_{i=1}^{\vert 2m+1\vert} \epsilon_i D_i(n,m)
\end{equation}

for some basepoint on~$D_i(n,m)$ and some signs~$\epsilon_i\in\{\pm 1\}$. 
Equation~(\ref{e:odd}) follows from Equation~(\ref{e:disum}),
the induction hypothesis and the subsequent lemma.
\end{proof}

Let $J_k$ be the two-sided ideal of $\Conw(\Sigma\times I)$ generated by 
knots~$K_{2n}$ with $-k<n< k$. 
Recall the definition of the diagram~$D_i(n,m)$ from the 
proof of Lemma~\ref{l:K2k0} (see Figure~\ref{f:Didiag}).
Then the following holds.

\begin{lemma}\label{l:Di}
For $n\in\N$, $m\in\Z$ and $i\leq \vert 2m+1\vert$ we have
$D_i(n,m)=p(z)K_{2k}+a$, where $p(z)\in\Z[z]$, $k=n+\vert m\vert +1$, and~$a$ is some
element of~$J_k$. 
\end{lemma}
\begin{proof}
A binary tree with root $L=D_i(n,m)$ is given by 
the following recursive description:

1. We start at the top left of the diagram of~$L$ and follow the orientation of the strand
until we come back to the point where we started 
or until we travel along an undercrossing strand of
a crossing where we did not travel along the overcrossing strand of that crossing
before.

2.a If we come back to the starting point on~$L$ in step~1, then 
the tree for~$L$ consists of a leaf labeled by~$L$.

2.b If we reach a crossing first as an undercrossing in step 1, then we
apply Relations~$(Bas)$ and~$(Ord)$ and then a skein relation to 
this crossing. We obtain $L=\epsilon_1 L_1+\epsilon_2 z L_2$ for 
some~$\epsilon_j\in\{\pm 1\}$.
The tree for~$L$ consists of a root vertex connected to a tree for~$L_1$ and
to a tree for~$L_2$.

\smallskip

The knot~$D_i(n,m)$ is equal to a weighted sum 
of the labels~$L'$ of the leaves in this tree with weights~$\pm z^j$, where~$j$ is
equal to the number of crossings spliced on the path in the tree leading from the
root to~$L'$.
The leaves are labeled by links~$L'$ that are products~$K L''$, where~$K$ is a descending knot.
We will further examine the paths in the tree leading from the root to a leaf.
There is a unique path in the tree where no crossing is spliced and 
the label of the leaf at the end of this path is $K_{2n+2m+2}$.

If a crossing of~$D_i(n,m)$ is spliced, then the first one has to be one of the
crossings marked by a~$*$ in Figure~\ref{f:Didiag}. Let $D^s$ be the resulting link diagram
(see Figure~\ref{f:Dijdiag}). 

\begin{figure}[!h]
\centering
\setbox1=\hbox{\input{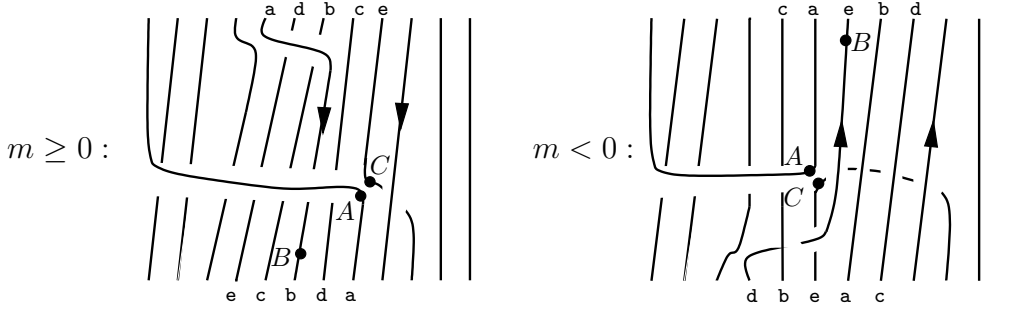}}
$\vcenter{\box1}$
\caption{The diagram $D^s$ obtained from $D_i(n,m)$}\label{f:Dijdiag}
\end{figure}

For $m\geq 0$ (resp.\ $m<0$)
let~$B$ be a point on~$D^s$ near the~$1+n+\vert m\vert$-th endpoint at the bottom 
(resp.\ at the top) of~$D^s$.
Traveling along~$D^s$ as in step~1 of the description of the tree for~$D^s$
we first pass through a strand of the crossing that was
spliced by passing from $K_{2n+1}K_{2m+1}$ to $D_i(n,m)$. 
We label a point 
on this strand by~$A$ and a point
on the other strand
of this spliced crossing by~$C$ (see Figure~\ref{f:Dijdiag}).
When we continue to travel along $D^s$ we will arrive first at the point~$B$
and then at the point~$C$ on~$D^s$ without the need of a
further modification of $D^s$.
Forgetting the component~$K'$ of~$D^s$ containing the points~$A, B, C$, we obtain a 
knot~$K''$. 
The upper and lower boundary points of the diagram of~$K''$ consist of a 
subset of the corresponding boundary points of~$K_{2m+1}$ in Figure~\ref{f:ko1ko2},
the $\delta_{1\sgn(m)}+\vert m\vert$-th upper and lower boundary point 
does not belong to~$K''$ because~$B$ lies on~$K'$, and
the strands of~$K''$ to the right of the point~$B$ do not intersect.
This implies that~$K''$ is the closure of a tangle on~$2r$ strands as in 
Figure~\ref{f:stanclos}.
The link~$D^s$ is not the product of~$K'$ 
and~$K''$, but the point~$C$ lies between the~$r$-th and $r+1$-st strand of~$K''$ and
when we travel along~$K'$ from~$C$ back to the starting point we 
cross~$K''$ only in strands to the right of~$C$.
This implies that no matter how~$D^s$ is modified in the passage from~$D^s$ to a 
leaf~$L$, the component of~$L$ containing the points~$A, B, C$ will always be homotopic
to~$K_{2j}$ for 
some~$j\in\{n-\vert m\vert,\ldots,n+\vert m\vert+1\}\subseteq\{-k+1,\ldots k\}$.
We have~$K_{2j}\in J_k$ for~$\vert j\vert<k$. 
If the label~$L$ of the leaf of the tree is of the form $L=K_{2k}L'$ for 
some link $L'$, then $L=K_{2k}$. 
\end{proof}

\medskip

Combining Lemma~\ref{l:K2k0} with
Lemma~\ref{l:cabmodgen} we obtain the following proposition.

\begin{prop}\label{p:conquogen}
The $\cori$-module $\Conquo(\Sigma\times I)$ is generated by the following set of links:

$$
\gen=\left\{
\begin{array}{l|l}
 & 
n\in\N, 
w_i\in \Conjrep^*, w_i\leq w_{i+1},\\
\raisebox{1.5ex}[0cm][0cm]{$K_{w_1}\ldots K_{w_n}$}
& (w_i=w_{i+1}\implies \sigmab(w_i)=0),
(w_i=e\implies i=n=1)
\end{array}
\right\}.
$$
\end{prop}
\begin{proof}
The cabled descending knots $\hat{K}_{w}$ with $\hat{K}_{w}\not=K_{w}$ are indexed
by $w\in\Conjrep\setminus(\Conjrep^*\cup\{e\})$. 
By definition they are of the 
form~$\hat{K}_{w}={i_v}_*(K_{2n})$ for some~$v\in\pi_1(\Sigma)$ with~$\sigmab(v)=1$ 
and~$n>0$. 
By Lemma~\ref{l:K2k0} we
have $K_{2n}=0\in\Conquo(X\times I)$ implying $\hat{K}_{w}=0\in\Conquo(\Sigma\times I)$.
Also by Lemma~\ref{l:K2k0}, we have~$2K_1^2=0\in\Conquo(X\times I)$ 
implying~$i_*(K_1^2)=0\in\Conquo(\Sigma\times I)$ for an arbitrary 
embedding~$i:X\times I\longrightarrow\Sigma\times I$.
This shows that the ideal~$J$ of Lemma~\ref{l:cabmodgen} is contained 
in~$\Tor_{\Z[z]}(\Conw(\Sigma\times I))$.
It is easy to see that for all $3$-manifolds~$M$ we 
have~$zK_eK_w=0\in \Conw(M)$. For non-orientable~$M$ we have~$K_e=0$ and~$e\not\in\Conjrep^*$.
The three arguments from above together with Lemma~\ref{l:cabmodgen} imply the proposition.
\end{proof}

%
%

\section{The weight system of the Conway polynomial}\label{s:wC}

We used ordered based links 
to define the Conway
polynomial~$\nabla$
taking values in~$\cori\otimes\Ring(\Sigma\times I)$.
In this section we will construct a map~$W_\glz^\ob$ from
a space of ordered based chord diagrams to~$\Ring(\Sigma\times I)\otimes\Q$. 
In the next section
we will use the map~$W_\glz^\ob$ together with a universal Vassiliev
invariant~$Z_{\Sigma\times I}^\ob$ to
show that~$\nabla$ is well-defined. 

Let $G$ be a group. Let $\Gamma={S^1}^{\amalg l}$ 
be a disjoint union of~$l$ oriented circles.
A {\em $G$-\labeled{} ordered based chord diagram}~$D$ with support~$\Gamma$
consists of a finite set $S=A\amalg B\amalg C$ of mutually distinct points
on~$\Gamma$, such that

\begin{itemize}
\item on each circle lies exactly one element of $B$ called {\em basepoint},
the set~$B$ is linearly ordered,

\item $C$ is partitioned into subsets of cardinality two called {\em chords}, and

\item to each point~$p$ in~$A$ there is assigned an element of the group~$G$ called {\em label}
of~$p$.
\end{itemize}

Usually we call~$G$-\labeled{} ordered based chord diagrams
simply chord diagrams.
We consider two chord diagrams~$D$ and~$D'$ with support~$\Gamma$
and~$\Gamma'$ respectively as being equal,
if there exists a homeomorphism
between~$\Gamma$ and~$\Gamma'$
that preserves all additional data.
Define
the {\em degree}~$\deg(D)$ of a chord diagram~$D$
as the number of its chords.

We represent a chord diagram graphically as follows:
the circles of the $1$-manifold~$\Gamma$ are oriented counterclockwise in the pictures.
We connect the two endpoints of a chord by a thin line.
The labels of a chord diagram are represented by
marking the points of~$A$ by a dot on a circle and by
writing the labels close to these marked points.
The basepoints are marked by the symbol~$\times$ on the circle.
We draw the basepoints from the left to the right in increasing
order or we label the basepoints by elements of an ordered set.
An example of a picture of a $\Z$-\labeled{} chord diagram~$D$
of degree~$5$ is shown in Figure~\ref{f:Dex}.

\begin{figure}[!h]
\centering 
\setbox1=\hbox{\input{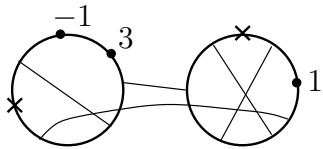}}
$\vcenter{\box1}$
\caption{A picture of a
$\Z$-\labeled{} chord diagram} \label{f:Dex}
\end{figure}

In the following definitions we define vector spaces by generators and
graphical relations.
We use the convention that all
diagrams in a graphical relation coincide everywhere except
for the parts we show, and that all configurations of the hidden
parts are possible.

\begin{defi}\label{d:A}
Let $G$ be a group, and $\sigmab:G\longrightarrow\Z/2$ a homomorphism
of groups.
Define~$\A^\ob(G,\sigma)$ to be the graded
$\Q$-vector space generated by $G$-\labeled{} ordered based
chord diagrams~$D$ modulo the Relations~(4T), ($\sigma$-Nat), (Rep), (Bas), and~(Ord) below. 
\end{defi}

The four-term relation (4T):\nopagebreak

$$
\Picture{
\thicklines
\put(-0.8,1){\vector(0,-1){2}}
\put(0,1){\vector(0,-1){2}}
\put(0.8,1){\vector(0,-1){2}}
\thinlines
\put(-0.8,-0.3){\line(1,0){0.8}}
\put(0,0.4){\line(1,0){0.8}}
}\quad + \quad
\Picture{
\thicklines
\put(-0.8,1){\vector(0,-1){2}}
\put(0,1){\vector(0,-1){2}}
\put(0.8,1){\vector(0,-1){2}}
\thinlines
\put(-0.8,-0.3){\line(1,0){1.6}}
\put(0,0.4){\line(1,0){0.8}}
} \quad = \quad
\Picture{
\thicklines
\put(-0.8,1){\vector(0,-1){2}}
\put(0,1){\vector(0,-1){2}}
\put(0.8,1){\vector(0,-1){2}}
\thinlines
\put(-0.8,0.4){\line(1,0){0.8}}
\put(0,-0.3){\line(1,0){0.8}}
} \quad + \quad
\Picture{
\thicklines
\put(-0.8,1){\vector(0,-1){2}}
\put(0,1){\vector(0,-1){2}}
\put(0.8,1){\vector(0,-1){2}}
\thinlines
\put(-0.8,0.4){\line(1,0){1.6}}
\put(0,-0.3){\line(1,0){0.8}}
}
$$

\bigskip\bigskip

The Relation~($\sigma$-Nat):
$\quad\begin{picture}(2,1)(-1,-0)
\thicklines
\put(-1,1){\vector(1,0){2}}
\put(-1,-0.5){\vector(1,0){2}}
\thinlines
\put(-0.4,-0.5){\line(0,1){1.5}}
\put(0.3,-0.5){\circle*{0.15}}
\put(0.3,1){\circle*{0.15}}
\put(0.1,0.5){\makebox(0.4,0.5){$s$}}
\put(0.1,-1){\makebox(0.4,0.5){$s$}}
\end{picture}\ = \ (-1)^{\sigma(s)} \
\begin{picture}(2,1.1)(-1,-0.1)
\thicklines
\put(-1,1){\vector(1,0){2}}
\put(-1,-0.5){\vector(1,0){2}}
\thinlines
\put(0.3,-0.5){\line(0,1){1.5}}
\put(-0.4,-0.5){\circle*{0.15}}
\put(-0.4,1){\circle*{0.15}}
\put(-0.6,0.5){\makebox(0.4,0.5){$s$}}
\put(-0.6,-1){\makebox(0.4,0.5){$s$}}
\end{picture}$

\bigskip\smallskip

The Relations~(Rep):
$\quad\begin{picture}(2,0.5)(-1,-0)
\thicklines
\put(-1,0.25){\vector(1,0){2}}
\put(-0.5,0.25){\circle*{0.15}}
\put(0.4,0.25){\circle*{0.15}}
\put(-0.6,-0.3){\hbox{\hss$a$\hss}}
\put(0.3,-0.3){\hbox{\hss$b$\hss}}
\end{picture}\ = \
\begin{picture}(2,0.5)(-1,0)
\thicklines
\put(-1,0.25){\vector(1,0){2}}
\put(0,0.25){\circle*{0.15}}
\put(-0.25,-0.3){\hbox{\hss$ab$\hss}}
\end{picture}\quad , \quad
\begin{picture}(2,0.5)(-1,0)
\thicklines
\put(-1,0.25){\vector(1,0){2}}
\put(0,0.25){\circle*{0.15}}
\put(-0.4,-0.45){\makebox(0.8,0.5){$e$}}
\end{picture}\ = \
\begin{picture}(2,0.5)(-1,0)
\thicklines
\put(-1,0.25){\vector(1,0){2}}
\end{picture}\ ,$

\medskip

where~$e$ is the neutral element of~$G$ and $ab$ denotes the product
of~$a$ and~$b$ in the group~$G$.

\smallskip

The Relations~$(Bas)$: \nopagebreak

\begin{center}
\setbox1=\hbox{\input{tex/bas}}
$\vcenter{\box1}$
\end{center}\nopagebreak

where the element $b\in G$ is defined as the product of
the labels on the hidden part of the shown component 
(the order of the multiplication is not important). 

\smallskip

The Relation~$(Ord)$: \nopagebreak

\begin{center}
\setbox1=\hbox{\input{tex/relord}}
$\vcenter{\box1}$
\end{center}\nopagebreak

meaning that the change of the order of two neighbored
circles gives the sign contribution $(-1)^{\sigmab(a)\sigmab(b)}$,
where~$a$ and~$b$ are the products of the labels on the two circles.

\begin{defi}\label{d:ab}
Let $\Ab^\ob(G,\sigma)$ be the quotient of $\A^\ob(G,\sigma)$ by the 
framing
independence relation~(FI):
$\ 
\begin{picture}(2,0.5)(-1,0.4)
\thicklines
\put(-1,0.25){\vector(1,0){2}}
\thinlines
\qbezier[70](-0.6,0.25)(-0.6,0.75)(-0.1,0.75)
\qbezier[70](0.4,0.25)(0.4,0.75)(-0.1,0.75)
\end{picture} \quad = \quad 0.
$
\end{defi}

Extending the order of the basepoints of chord diagrams $D_1$, $D_2$ to the disjoint
union~$D_1\amalg D_2$ 
by requiring that
all basepoints of~$D_1$ are smaller than those of~$D_2$ 
induces a multiplication on~$\A(G,\sigma)$ turning it into a 
graded ring.  
If $G=\pi_1(M)$ and $\sigma:\pi_1(M)\longrightarrow\Z/2$ is the first Stiefel-Whitney class 
of~$M$, then one has an isomorphism 
of rings $\Conw_0(M)\otimes \Q\cong\A_0(G,\sigma)$, where $\A_0(G,\sigma)$ is
the degree-$0$ part of $\A(G,\sigma)$. 

Let~$D$ be a picture of a \obG{} such that the basepoint is the 
highest point on each
circle and all labels lie on a horizontal line. An example is 
shown in Figure~\ref{obGpic}.

\begin{figure}[!h]
\centering
\setbox1=\hbox{\input{tex/Wglobex1}}
$\vcenter{\box1}$
\caption{A special picture of a \obG}\label{obGpic}
\end{figure}

We used the order of the oriented circles and the basepoints to draw
Figure~\ref{obGpic}. Now we forget this data. We
replace each chord as shown in Figure~\ref{f:Wgl2}.

\begin{figure}[!h]
\centering
\setbox1=\hbox{\input{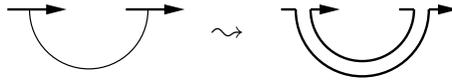}}
$\vcenter{\box1}$
\caption{Replacing a chord}\label{f:Wgl2}
\end{figure}

The result is an immersion of \labeled{} circles into the
plane. By a homotopy this immersion can be rearranged to a {\em
standard embedding}, meaning that 
all labels lie on a horizontal line, no parts of
the circles lie below this line, and the projections of the circles to
the horizontal line are disjoint. 
As an example, after having replaced the chords of the diagram in 
Figure~\ref{obGpic} we can rearrange the immersion to a standard
embedding as shown in Figure~\ref{stdemb}.

\begin{figure}[!h]
\centering
\setbox1=\hbox{\input{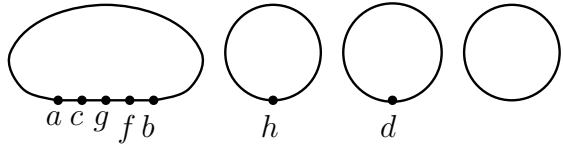}}
$\vcenter{\box1}$
\caption{A standard embedding of \labeled{} circles}\label{stdemb}
\end{figure}

During the homotopy we count the number $m$ of times the order
of the projections of
two labels $s, s'$ with $\sigma(s)=\sigma(s')=1$
to the horizontal line changes.
For the example in Figures~\ref{obGpic} and~\ref{stdemb} we have 
$(-1)^m=(-1)^{\sigmab(b)\sigmab(cfg)+
\sigmab(d)\sigmab(fgh)+\sigma(f)\sigma(g)}$.
The chosen order of the circles
and the chosen highest points of the circles in
a standard embedding of labeled circles
determine an ordered based chord diagram (without chords)
that we call $D'$ for a moment. 

\begin{prop}\label{p:Wglob}
With the notation from above, the linear map\nopagebreak

$$
W_\gl^\ob:\A^\ob(G,\sigma)\longrightarrow \A_0^\ob(G,\sigma)
$$\nopagebreak

that sends a $G$-\labeled{} ordered based chord diagram~$D$ to
$(-1)^m D'$ is well-defined.
\end{prop}
\begin{proof}
Changing the order of two circles or pushing a label around a circle
in a standard embedding of labeled circles as in
Figure~\ref{stdemb} gives the same sign contribution to~$(-1)^m$ and
to~$D'$. So~$W_\gl^\ob$ is well-defined for \obG s. We have to
verify that~$W_\gl^\ob$ respects the defining relations
of~$\A^\ob(G,\sigma)$:

The map $W_\gl^\ob$ respects the relation
shown in Figure~\ref{f:2T}, where any order of the shown parts of the
diagram is possible, because one can slide one thickened chord
as in Figure~\ref{f:Wgl2} along the other one (see~\cite{BNG}, Section~3). 
This relation implies the Relation~(4T).

\begin{figure}[!h]
\centering
\setbox1=\hbox{\input{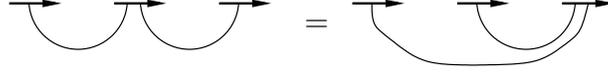}}
$\vcenter{\box1}$
\caption{A two term relation}\label{f:2T}
\end{figure}

We can easily verify the compatibility of~$W_\gl^\ob$ with the
Relations~(Rep),~(Bas), and~(Ord).
The compatibility of $W_\gl^\ob$ with the Relation~($\sigma$-Nat) follows
from Figure~\ref{f:sNatc} which shows an equation between
elements of $\A_0^\ob(G,\sigma)$ defined by
immersed \labeled{} circles as in the definition of~$W_\gl^\ob$. 

\begin{figure}[!h]
\centering
\setbox1=\hbox{\input{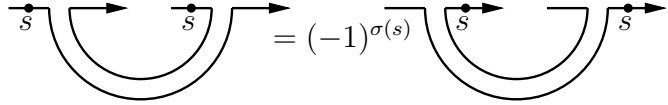}}
$\vcenter{\box1}$
\caption{The compatibility of $W_\gl^\ob$ with the
Relation~($\sigma$-Nat)}\label{f:sNatc}
\end{figure}

The formula in Figure~\ref{f:sNatc} holds because when we push
the labels~$s$ to their new positions, then
in the projection to the horizontal line they 
are commuted with the same labels from the remaining part of
the diagram and the commutation between the~$s$ on the right side with
the~$s$ on the left side gives the
sign contribution $(-1)^{\sigma(s)}$.
\end{proof}

\medskip

Let $p:\A_0(\pi_1(\Sigma),\sigma)\longrightarrow\Ring(\Sigma\times I)\otimes \Q$ 
be the homomorphism
of rings induced by mapping a circle with a single mark~$g$ to~$\overK_g$.
Replacing the chords in a Relation~(FI) as in Figure~\ref{f:Wgl2} produces a circle without
labels together with at least one additional component. 
This implies the following corollary.

\begin{coro}\label{c:Wb}
For non-orientable surfaces~$\Sigma$ the map $p\circ W_\gl^\ob$ descends to a 
map\nopagebreak

$$
W_\glz^\ob:\Ab(\pi_1(\Sigma),\sigma)\longrightarrow \Ring(\Sigma\times I)\otimes \Q.
$$
\end{coro}

The reason for the notation $W_\gl^\ob$ and~$W_\glz^\ob$ is given in
Section~1.17 of \cite{Lie}, where~$W_\gl^\ob$ and
$W_\glz^\ob$ are related to weight systems associated 
with the Lie superalgebra~$\gl(V)$ of endomorphisms of 
a $\Z/2$-graded vector space~$V$. In the case of~$W_\glz^\ob$ the
superdimension of~$V$ is~$0$. We will make use of this relation in Section~\ref{s:squareproof}.


The map~$W_\gl^\ob$ can be described recursively as follows.
For chord diagrams~$D$ of degree~$0$ we have~$W_\gl^\ob(D)=D$.
If~$D$ is a \obG{} with
at least one chord~$c$, then one can use Relations~$(Bas)$ and~$(Ord)$ such that~$D$ looks
like one of the two diagrams to the left of
the symbol~$\leadsto$ in Figure~\ref{f:Wglrec}.

\begin{figure}[!h]
\centering
\setbox1=\hbox{\input{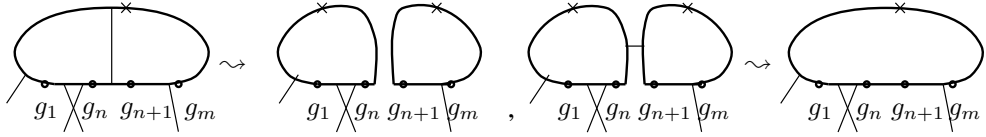}}
$\vcenter{\box1}$
\caption{A recursive description of $W_\gl^\ob$}\label{f:Wglrec}
\end{figure}

In Figure~\ref{f:Wglrec} only the chord~$c$ of $D$ is replaced as shown in
Figure~\ref{f:Wgl2}. We 
obtain a \obG{}~$D(c)$ with~$\deg(D(c))=\deg(D)-1$ and with

\begin{equation}\label{e:recWglob}
W_\gl^\ob(D)=W_\gl^\ob(D(c)).
\end{equation}

We can express the recursion formula for $W_\gl^\ob$ from Figure~\ref{f:Wglrec} without
using embeddings by the local replacement rules 
involving basepoints as shown in Figure~\ref{f:Wglrec2}.

\begin{figure}[!h]
\centering
\setbox1=\hbox{\input{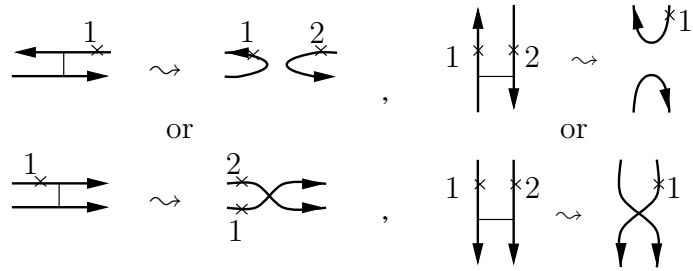}}
$\vcenter{\box1}$
\caption{Descriptions of $W_\gl^\ob$ and $W_\glz^\ob$ 
by local replacement rules}\label{f:Wglrec2}
\end{figure}

\section{Proof of Theorems~\ref{t:conwaypoly} and~\ref{t:specase}}\label{s:eC}

Let $\Sigma=\bigcup_{i=0}^k B_i$ be a decomposed surface or let~$\Sigma=P^2$. 
Let~$L$ be a link in~$\Sigma\times I$.
Recall the definition of 
the universal Vassiliev invariant~$Z_{\Sigma\times I}(L)$
from~\cite{Lie}.
For a decomposed surface~$\Sigma$ 
this definition is roughly 
as follows: after applying an isotopy we can assume that~$L\cap B_i\times I$ 
is in a standard position for all~$i>0$. The 
tangle~$L\cap B_0\times I$ is turned into a non-associative tangle~$T$ by choosing a 
bracketing on the ordered set $L\cap (\partial B_0)\times I$ satisfying some conditions.
Then $Z_{\Sigma\times I}(L)$ is obtained from 
the universal Vassiliev invariant~$Z(T)$ 
of the non-associative tangle~$T$ (see~\cite{LeM}, \cite{BN2}) 
by gluing labeled intervals to those pairs of boundary points of diagrams in~$Z(T)$ whose
corresponding boundary points of~$T$ are connected by an interval 
in~$L\cap (\Sigma\setminus B_0)\times I$.

Let~$\L^\ob(M)$ be the set of isotopy
classes of ordered based links in~$M$. Let~$\widehat{\A}^\ob(M)$ be
the graded completion of $\Ab^\ob(M)=\Ab^\ob(\pi_1(M),\sigma)$. Define a map

\begin{equation}
Z_{\Sigma\times I}^\ob:\L^\ob(\Sigma\times I)\longrightarrow
\widehat{\A}^\ob(\Sigma\times I)
\end{equation}

by equipping the chord diagrams in the series 
$Z_{\Sigma\times I}(L)$ with an
order and with basepoints according to the order and the basepoints of~$L$.
The map~$Z^\ob_{\Sigma\times I}$ is well-defined because chord endpoints
commute with the basepoints by the Relation~$(Bas)$ (see~\cite{Li2}).
Let~$\hRing(M):=(\Ring(M)\otimes\Q)[[h]]$. 
Define 

\begin{equation}\label{e:defWhat}
\widehat{W}_\glz^\ob:\widehat{\A}^\ob(G,\sigma)\longrightarrow
\hRing(\Sigma\times I)
\end{equation}

by extending $\widehat{W}_\glz^\ob(D)=W_\glz^\ob(D)h^{\deg\, D}$ for any \obG{}~$D$ to formal
power series.
Finally, define the invariant~$\Ch$ of ordered based links by

$$
\Ch=\widehat{W}_\glz^\ob\circ Z^\ob_{\Sigma\times I}.
$$

Turn $\hRing(M)$ into a $\cori$-module 
by~$z\cdot a:=(e^{h/2}-e^{-h/2})a$ for $a\in \hRing(M)$.

\begin{lemma}\label{l:Ch}
The invariant $\Ch$ of ordered based links induces a $\cori$-linear 
map

\begin{eqnarray}
& & \Ch:\Conquo(\Sigma\times I)\longrightarrow\hRing(\Sigma\times I) \mbox{ satisfying}\nonumber\\
& & \Ch(L)\equiv \overL\, \mod\, h\quad\mbox{for all ordered
based links~$L$.}\label{e:equivtL}
\end{eqnarray}
\end{lemma}
\begin{proof}
Let
$(L_+,L_-,L_{\vert\vert})$ be a skein triple of ordered based links.
Let $B_1=\picBo$ and $B_2=\picBt{}$. Define the composition of diagrams 
$D_1, D_2$ by placing~$D_1$ onto the top of~$D_2$.
In the recursive description of~$W_\gl^\ob$ (see Figure~\ref{f:Wglrec2}) a
chord $B_i\circ {\setlength{\unitlength}{15pt}\picChord}$
between two parallel parts of the
oriented circles is replaced by $B_{3-i}\circ{\setlength{\unitlength}{15pt}\picX}$.
By the explicit description of $Z_{\Sigma\times I}^\ob$
and by Equation~(\ref{e:recWglob}) we have for some~$i\in\{1,2\}$:

\begin{eqnarray*}
& & \Ch(L_+)-\Ch(L_-)\\
& = & \widehat{W}_\glz^\ob\left(B_i\circ\exp\left({\picChord\, /2}\right)\circ\picX
\right)
-\widehat{W}_\glz^\ob
\left(B_i\circ\exp\left({-\picChord\, /2}\right)\circ\picX\right)\\
& = & \sum_{n=0}^\infty \frac{1}{4^{n}(2n+1)!}
\widehat{W}_\glz^\ob\left(B_i\circ\picChord^{2n+1}\circ\picX\right)\\
& = & \sum_{n=0}^\infty \frac{h^{2n+1}}{4^{n}(2n+1)!}
\widehat{W}_\glz^\ob\left(B_{3-i}\right) = \left(e^{h/2}-e^{-h/2}\right)\Ch(L_{\vert\vert}) 
= z\cdot \Ch(L_{\vert\vert}).
\end{eqnarray*}

It is obvious by Relations~$(Ord)$ and~$(Bas)$ for ordered based links and chord diagrams 
that~$\Ch$ is well-defined 
on~$\Conw(\Sigma\times I)$. Since~$\hRing(\Sigma\times I)$ has no  
non-trivial torsion elements as a~$\cori$-module, $\Ch$ descends to~$\Conquo(\Sigma\times I)$.
Equation~(\ref{e:equivtL}) follows because for ordered based links~$L$ 
whose $i$-th component is homotopic to~$w_i\in\pi_1(\Sigma\times I)$ 
we have~$Z_{\Sigma\times I}^\ob(L)\equiv D_L\;\mod\;h$ where
$D_L\in\Ab(\Sigma\times I)_0$ 
is a product of~$r$ oriented circles with a single label~$w_i$ and 
$W_{\glz}^\ob(D_L)=\pi(D_L)=\overL\in\Ring(\Sigma\times I)\otimes \Q$.
\end{proof}

\medskip

It remains to make a change of parameters to show that the Conway
polynomial is well-defined.

\medskip

\begin{proof}[Proof of Theorem~\ref{t:conwaypoly}]
(1) The uniqueness of~$\nabla$ follows because the $\cori$-module
$\cori\otimes \Ring(\Sigma\times I)$ has a trivial torsion submodule and
the condition $\nabla(L)=\overL$ for descending links~$L$ with $\theta(L)\not=0$ 
prescribes the image of the generators~$\gen$ of~$\Conquo(\Sigma\times I)$ 
(see Proposition~\ref{p:conquogen}).
Let us prove the existence of~$\nabla$.
Proposition~\ref{p:ringiso} implies that the set~$\genb=\{\overL\;\vert\; L\in\gen\}$
is a basis of the~$\cor$-module $\Ring(\Sigma\times I)$.
Define a $\cori$-linear map $\kappa:\cori\otimes_\Z\Ring(\Sigma\times I)\longrightarrow
\Conquo(\Sigma\times I)$ on basis elements by 

$$
\kappa(1\otimes \overK_{w_1}\ldots \overK_{w_n}):=K_{w_1}\ldots K_{w_n}.
$$

The map~$\kappa$ is surjective. 
Consider the composition 

$$
\begin{CD}
\cori\otimes_\Z\Ring(\Sigma\times I) @>\kappa>> 
\Conquo(\Sigma\times I)@>\Ch{}>> {\hRing}(\Sigma\times I).
\end{CD}
$$

By Lemma~\ref{l:Ch} we have~$\Ch(\kappa(\overline{L}))\equiv\overline{L}\;\mod\;h$
for all ordered based links~$L$.
This implies that~$\genb$ is mapped injectively to the basis~$\Ch(\kappa(\genb))\,\mod\,h$ of the
$\Q$-vector space~$\hRing(\Sigma\times I)/(h)$.
Considering inductively the sets~$\Ch(\kappa(\genb))\,\mod\,h^n$
we see that $\Ch(\kappa(\genb))$
is linearly independent over~$\Q[z]$ (where $z\cdot a=
(e^{h/2}-e^{-h/2})a$).
This implies that~$\Ch\circ \kappa:\cori\otimes\Ring(\Sigma\times I)
\longrightarrow\hRing(\Sigma\times I)$ 
is injective. Hence~$\kappa$ is an isomorphism and the map 
$\kappa^{-1}:\Conquo(\Sigma\times I)\longrightarrow\cori\otimes\Ring(\Sigma\times I)$ 
has the property $\kappa^{-1}(L)=\overL$ for all
$\Conjrep$-descending links~$L$ with~$\theta(\overL)\not=0$.
We obtain~$\nabla$ as the
composition of the canonical projection~$\Conw(\Sigma\times I)
\longrightarrow\Conquo(\Sigma\times I)$ with~$\kappa^{-1}$.
\end{proof}

\medskip

\begin{proof}[Proof of Theorem~\ref{t:specase}]
(1) 
Let $\Sigma$ be the M\"obius strip~$X$. Let $L\subset X\times I$ be a~$\Conjrep$-descending
link. If~$\theta\left(\overL\right)=0$, then~$L=0\in\Conquo(X\times I)$ by Lemma~\ref{l:K2k0}. 
Theorem~\ref{t:conwaypoly} implies~$\nabla(L)=0$ because~$\nabla$ factors through
$\Conquo(X\times I)$.

It follows easily from Lemma~\ref{l:K2k0} and Theorem~\ref{t:conwaypoly} that 
$\nabla(L_1L_2)=\nabla(L_1)\nabla(L_2)$ for all $\Conjrep$-descending links
$L_1, L_2\subset X\times I$. Lemma~\ref{l:descgen} implies this formula for all
ordered based links~$L_1, L_2\subset X\times I$.

The case $\Sigma=P^2$ of the theorem can be deduced from the case~$\Sigma=X$. 
For~$\Sigma=I^2$ and~$\Sigma=S^1\times I$ Part~(1) of the theorem is trivial.

\smallskip
 
(2) It is enough to show $L=(-1)^{\vert L\vert_0} L^*\in\Conquo(\Sigma\times I)$ 
for~$\Sigma=X$.
Let $L$ be a diagram of an ordered based link in $X\times I$.
Let~$k_1(L)$ be the 
number of crossings of~$L$ and let~$k_2(L)$ be the minimal number of crossing
changes needed to pass from~$L$ to a diagram of an~$\Conjrep$-descending link. 
We prove the theorem by induction on the lexicographical order on $(k_1(L),k_2(L))$. 

Assume
that $k_2(L)=0$. Then $L=K_{i_1}\ldots K_{i_r}$ is a product of descending knots. By 
Lemma~\ref{l:K2k0} we have $L=0=L^*$ if $\vert L\vert_0>0$. Assume that $\vert L\vert_0=0$.
Then Lemma~\ref{l:K2k0} and Lemma~\ref{l:shortcuts} imply 

\begin{eqnarray*}
K_{i_1}\ldots K_{i_r} & = & (-1)^{r(r-1)/2}K_{i_r}\ldots K_{i_1}=(-1)^{r(r-1)/2}
K_{i_r}^*\ldots K_{i_1}^*\\
 & = & \left(K_{i_1}\ldots K_{i_r}\right)^*=(-1)^{\vert L\vert_0}
\left(K_{i_1}\ldots K_{i_r}\right)^*.
\end{eqnarray*}

If $k_1(L)=0$, then we also have $k_2(L)=0$ and we are back in the first case.
Now let $k_1(L)>0$ and $k_2(L)>0$. 
Choose a crossing of $L$ such that $k_2(L_1)=k_2(L)-1$, where~$L_1$ is obtained 
from~$L$ by changing this crossing and let~$L_2$ be obtained by splicing the
crossing. Then for an $\epsilon\in\{\pm 1\}$ we have by induction

\begin{eqnarray*}
L & = & L_1+\epsilon zL_2=(-1)^{\vert L_1\vert_0} L_1^*+\epsilon z (-1)^{\vert L_2\vert_0}
L_2^*\\
& = & (-1)^{\vert L\vert_0}\left(L_1^*-\epsilon z L_2^*\right)=(-1)^{\vert L\vert_0}L^*.
\end{eqnarray*}

In this computation we used that
for a skein triple $(L_+,L_-,L_{\vert\vert})$ of links we have 
$\vert L_+\vert_0=\vert L_-\vert_0=\vert L_{\vert\vert}\vert_0\pm 1$ and
$(L_-^*,L_+^*,L_{\vert\vert}^*)$ is also a skein triple.  

%
\end{proof}

\section{Coverings and the Conway polynomial}\label{s:squareproof}

Let $\Ab(G,\sigma)$ be the vector space generated by $G$-labeled chord diagrams 
(without order and basepoints) modulo the Relations
(4T), ($\sigma$-Nat), (Rep) and (FI) (see Section~\ref{s:wC}). 
For a $3$-manifold~$M$ we define $\Ab(M):=\Ab(\pi_1(M), \sigmab)$ where~$\sigmab$ is the
orientation character of~$M$.
The universal Vassiliev 
invariant~$Z_{\Sigma\times I}$ (see~\cite{Lie})
takes values in a completion of~$\Ab(\Sigma\times I)$. 
For a labeled diagram~$D$ define

$$
\sigmab(D)=\sum \sigmab(g)\in\Z/2
$$ 

where
the sum runs over all labels~$g$ on the skeleton of~$D$.

Let us recall a result of~\cite{Li2} for the special case of
the $2$-fold covering
$p:S^2\times I\longrightarrow P^2\times I$.
For a chord diagram~$D$
labeled by elements of~$\pi_1(P^2,*)\cong\Z/2$ 
we omit the points labeled by the neutral element in our pictures.
The points of~$D$ 
labeled by the non-trivial element are simply represented by a dot (without label).
There exists a map $p^*:\Ab(P^2\times I)\longrightarrow\Ab(S^2\times I)$
defined by the replacement rules shown in Equation~(\ref{e:defpstar}).

\begin{equation}\label{e:defpstar}
p^*\left( \mbox{\setbox1=\hbox{\begin{picture}(0,0)%
\epsfig{file=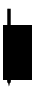}%
\end{picture}%
\setlength{\unitlength}{4144sp}%
\begingroup\makeatletter\ifx\SetFigFont\undefined%
\gdef\SetFigFont#1#2#3#4#5{%
  \reset@font\fontsize{#1}{#2pt}%
  \fontfamily{#3}\fontseries{#4}\fontshape{#5}%
  \selectfont}%
\fi\endgroup%
\begin{picture}(60,381)(409,41)
\end{picture}
}
$\vcenter{\box1}$}\right)=\mbox{\setbox1=\hbox{\begin{picture}(0,0)%
\epsfig{file=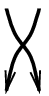}%
\end{picture}%
\setlength{\unitlength}{4144sp}%
\begingroup\makeatletter\ifx\SetFigFont\undefined%
\gdef\SetFigFont#1#2#3#4#5{%
  \reset@font\fontsize{#1}{#2pt}%
  \fontfamily{#3}\fontseries{#4}\fontshape{#5}%
  \selectfont}%
\fi\endgroup%
\begin{picture}(194,419)(421,0)
\end{picture}
}
$\vcenter{\box1}$}\quad , \quad 
p^*\left( \mbox{\setbox1=\hbox{\begin{picture}(0,0)%
\epsfig{file=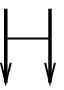}%
\end{picture}%
\setlength{\unitlength}{4144sp}%
\begingroup\makeatletter\ifx\SetFigFont\undefined%
\gdef\SetFigFont#1#2#3#4#5{%
  \reset@font\fontsize{#1}{#2pt}%
  \fontfamily{#3}\fontseries{#4}\fontshape{#5}%
  \selectfont}%
\fi\endgroup%
\begin{picture}(246,381)(418,41)
\end{picture}
}
$\vcenter{\box1}$}\right)=\mbox{\setbox1=\hbox{\begin{picture}(0,0)%
\epsfig{file=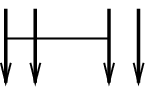}%
\end{picture}%
\setlength{\unitlength}{4144sp}%
\begingroup\makeatletter\ifx\SetFigFont\undefined%
\gdef\SetFigFont#1#2#3#4#5{%
  \reset@font\fontsize{#1}{#2pt}%
  \fontfamily{#3}\fontseries{#4}\fontshape{#5}%
  \selectfont}%
\fi\endgroup%
\begin{picture}(651,381)(418,41)
\end{picture}
}
$\vcenter{\box1}$}-\mbox{\setbox1=\hbox{\begin{picture}(0,0)%
\epsfig{file=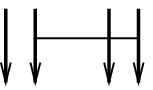}%
\end{picture}%
\setlength{\unitlength}{4144sp}%
\begingroup\makeatletter\ifx\SetFigFont\undefined%
\gdef\SetFigFont#1#2#3#4#5{%
  \reset@font\fontsize{#1}{#2pt}%
  \fontfamily{#3}\fontseries{#4}\fontshape{#5}%
  \selectfont}%
\fi\endgroup%
\begin{picture}(651,380)(418,42)
\end{picture}
}
$\vcenter{\box1}$}
\end{equation}

For a link~$L\subset P^2\times I$ we have 

\begin{equation}\label{e:rp2cover}
p^*\circ Z_{P^2 \times I}(L)=Z_{S^2\times I}(p^{-1}(L))
\end{equation}

(see \cite{Li2}, Theorem~4).
Links in $S^2\times I$ are in a natural bijection with links in~$S^3$ or~$\R^3$
and~$Z_{S^2\times I}$ corresponds to the usual Kontsevich integral under this bijection.

Comparing~$\Ab(P^2\times I)$ with~$\Ab^\ob(P^2\times I)$ we remark huge
differences: elements
of odd degree in $\Ab(P^2\times I)$ are~$0$
(see~\cite{Lie}, proof of Proposition~1.2) whereas this is not true 
for~$\Ab^\ob(P^2\times I)$ (otherwise the recursive description 
in Figure~\ref{f:Wglrec2} of the non-trivial
map~$W_\glz^\ob$  could not exist). 
However,
for ordered based
diagrams~$D$ on a single circle with $\sigmab(D)=1$, order and basepoints are
superfluous and the difference between $\Ab(P^2\times I)$ 
and~$\Ab^\ob(P^2\times I)$ vanishes for these special diagrams.

As for usual chord diagrams there are many descriptions of  
spaces isomorphic to $\Ab(G,\sigma)$. 
By $G$-labeled trivalent diagrams we shall mean usual trivalent diagrams~$D$ (see~\cite{BN1})
together with a distinguished subset of points on edges of~$D$ that are 
equipped with local orientations and labeled by 
elements of~$G$. For points on the skeleton of~$D$ we assume that the local orientation
coincides with the orientation of the skeleton. In particular, labeled chord diagrams
are labeled trivalent diagrams.
The space~$\Ab(G,\sigma)$ is isomorphic to the $\Q$-vector space 
generated by trivalent diagrams
modulo Relations~(FI), (Rep) (see Section~\ref{s:wC}),
(STU), (IHX), (AS) (see~\cite{BN1}) and modulo 
the Relations~(Comult) and (inv) shown below.

\begin{eqnarray*}
\mbox{(Comult)}  &  & \mbox{\setbox1=\hbox{\input{tex/comult}}
$\vcenter{\box1}$}\\
\mbox{(Inv1)} & & \mbox{\setbox1=\hbox{\input{tex/inv}}
$\vcenter{\box1}$}
\end{eqnarray*}

This can be proven along the same lines as in~\cite{BN1} (with a slight difference in the
case where an internal trivalent vertex connects to the skeleton of the diagram via two 
edges). A similar description by generators and
relations exists for~$\Ab^\ob(G,\sigma)$.
The space~$\Ab(G,\sigma)$ is a coalgebra with
comultiplication~$\Delta$ given by

\begin{equation}\label{e:comult}
\Delta(D)= \sum_{D=D_1\cup D_2} D_1\otimes D_2,
\end{equation}

where the sum runs over all trivalent diagrams $D_1, D_2\subset D$, such that
$D\setminus \Gamma=(D_1\setminus \Gamma)\amalg (D_2\setminus \Gamma)$.
Besides the comultiplication~$\Delta$, there exists a 
map

\begin{equation}\label{e:comultob}
\Delta^\ob:\Ab(G,\sigma)\longrightarrow\Ab^\ob(G,\sigma)\otimes\Ab^\ob(G,\sigma)
\end{equation}
 
defined by the same formula as in Equation~(\ref{e:comult}), where this time
$D_1$ and $D_2$ are equipped with order and basepoints that are arbitrarily chosen
for the first diagram and copied to the second.

Now consider $\Z/2$-labeled trivalent diagrams~$D$ with skeleton~$\Gamma$
where~$D\setminus \Gamma$ is connected and all labeled points lie close to univalent vertices
of~$D\setminus\Gamma$.
Define diagrams~$p_1(D)$ and~$p_2(D)$ 
by the local replacement rules shown in Equations~(\ref{e:defp1p2a}) 
and~(\ref{e:defp1p2b}).

\begin{eqnarray}
p_i\left(\mbox{\setbox1=\hbox{\begin{picture}(0,0)%
\epsfig{file=arrowdot.ps}%
\end{picture}%
\setlength{\unitlength}{4144sp}%
\begingroup\makeatletter\ifx\SetFigFont\undefined%
\gdef\SetFigFont#1#2#3#4#5{%
  \reset@font\fontsize{#1}{#2pt}%
  \fontfamily{#3}\fontseries{#4}\fontshape{#5}%
  \selectfont}%
\fi\endgroup%
\begin{picture}(60,381)(409,41)
\end{picture}
}
$\vcenter{\box1}$}\right)=\mbox{\setbox1=\hbox{\begin{picture}(0,0)%
\epsfig{file=cross.ps}%
\end{picture}%
\setlength{\unitlength}{4144sp}%
\begingroup\makeatletter\ifx\SetFigFont\undefined%
\gdef\SetFigFont#1#2#3#4#5{%
  \reset@font\fontsize{#1}{#2pt}%
  \fontfamily{#3}\fontseries{#4}\fontshape{#5}%
  \selectfont}%
\fi\endgroup%
\begin{picture}(194,419)(421,0)
\end{picture}
}
$\vcenter{\box1}$} & , & 
p_i\left(\mbox{\setbox1=\hbox{\begin{picture}(0,0)%
\epsfig{file=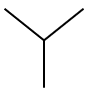}%
\end{picture}%
\setlength{\unitlength}{4144sp}%
\begingroup\makeatletter\ifx\SetFigFont\undefined%
\gdef\SetFigFont#1#2#3#4#5{%
  \reset@font\fontsize{#1}{#2pt}%
  \fontfamily{#3}\fontseries{#4}\fontshape{#5}%
  \selectfont}%
\fi\endgroup%
\begin{picture}(384,384)(439,17)
\end{picture}
}
$\vcenter{\box1}$}\right)=\mbox{\setbox1=\hbox{\begin{picture}(0,0)%
\epsfig{file=y.ps}%
\end{picture}%
\setlength{\unitlength}{4144sp}%
\begingroup\makeatletter\ifx\SetFigFont\undefined%
\gdef\SetFigFont#1#2#3#4#5{%
  \reset@font\fontsize{#1}{#2pt}%
  \fontfamily{#3}\fontseries{#4}\fontshape{#5}%
  \selectfont}%
\fi\endgroup%
\begin{picture}(384,384)(439,17)
\end{picture}
}
$\vcenter{\box1}$},\label{e:defp1p2a}\\
p_1\left(\mbox{\setbox1=\hbox{\begin{picture}(0,0)%
\epsfig{file=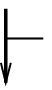}%
\end{picture}%
\setlength{\unitlength}{4144sp}%
\begingroup\makeatletter\ifx\SetFigFont\undefined%
\gdef\SetFigFont#1#2#3#4#5{%
  \reset@font\fontsize{#1}{#2pt}%
  \fontfamily{#3}\fontseries{#4}\fontshape{#5}%
  \selectfont}%
\fi\endgroup%
\begin{picture}(205,381)(418,41)
\end{picture}
}
$\vcenter{\box1}$}\right)\, = \,  
p_2\left(\mbox{\setbox1=\hbox{\begin{picture}(0,0)%
\epsfig{file=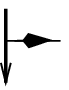}%
\end{picture}%
\setlength{\unitlength}{4144sp}%
\begingroup\makeatletter\ifx\SetFigFont\undefined%
\gdef\SetFigFont#1#2#3#4#5{%
  \reset@font\fontsize{#1}{#2pt}%
  \fontfamily{#3}\fontseries{#4}\fontshape{#5}%
  \selectfont}%
\fi\endgroup%
\begin{picture}(284,381)(418,41)
\end{picture}
}
$\vcenter{\box1}$}\right)\, = \, \mbox{\setbox1=\hbox{\begin{picture}(0,0)%
\epsfig{file=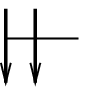}%
\end{picture}%
\setlength{\unitlength}{4144sp}%
\begingroup\makeatletter\ifx\SetFigFont\undefined%
\gdef\SetFigFont#1#2#3#4#5{%
  \reset@font\fontsize{#1}{#2pt}%
  \fontfamily{#3}\fontseries{#4}\fontshape{#5}%
  \selectfont}%
\fi\endgroup%
\begin{picture}(366,381)(418,41)
\end{picture}
}
$\vcenter{\box1}$} & , &
p_1\left(\mbox{\setbox1=\hbox{\begin{picture}(0,0)%
\epsfig{file=halfchd.ps}%
\end{picture}%
\setlength{\unitlength}{4144sp}%
\begingroup\makeatletter\ifx\SetFigFont\undefined%
\gdef\SetFigFont#1#2#3#4#5{%
  \reset@font\fontsize{#1}{#2pt}%
  \fontfamily{#3}\fontseries{#4}\fontshape{#5}%
  \selectfont}%
\fi\endgroup%
\begin{picture}(284,381)(418,41)
\end{picture}
}
$\vcenter{\box1}$}\right)\, = \,  
p_2\left(\mbox{\setbox1=\hbox{\begin{picture}(0,0)%
\epsfig{file=halfch.ps}%
\end{picture}%
\setlength{\unitlength}{4144sp}%
\begingroup\makeatletter\ifx\SetFigFont\undefined%
\gdef\SetFigFont#1#2#3#4#5{%
  \reset@font\fontsize{#1}{#2pt}%
  \fontfamily{#3}\fontseries{#4}\fontshape{#5}%
  \selectfont}%
\fi\endgroup%
\begin{picture}(205,381)(418,41)
\end{picture}
}
$\vcenter{\box1}$}\right) \, = \, \mbox{\setbox1=\hbox{\begin{picture}(0,0)%
\epsfig{file=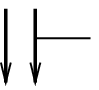}%
\end{picture}%
\setlength{\unitlength}{4144sp}%
\begingroup\makeatletter\ifx\SetFigFont\undefined%
\gdef\SetFigFont#1#2#3#4#5{%
  \reset@font\fontsize{#1}{#2pt}%
  \fontfamily{#3}\fontseries{#4}\fontshape{#5}%
  \selectfont}%
\fi\endgroup%
\begin{picture}(421,380)(418,42)
\end{picture}
}
$\vcenter{\box1}$}\label{e:defp1p2b}
\end{eqnarray}

It is easy to see that 
the map~$p^*$ satisfies Equation~(\ref{e:pstartrival}) 
for diagrams~$D$ as above.

\begin{equation}\label{e:pstartrival}
p^*(D)=p_1(D)+(-1)^{\deg\;D} p_2(D).
\end{equation}

For diagrams~$D$ where each of the~$c>1$ connected components of~$D\setminus \Gamma$ satisfies
the conditions above, $p^*(D)$ can be expressed as a sum with~$2^c$ terms by extending
Equation~(\ref{e:pstartrival})
multilinearly over these connected components. 
Using Relations (Comult) and~(Inv1), we can apply
this formula whenever~$D\setminus \Gamma$ is a forest.

We call a connected component of a diagram~$D\setminus \Gamma$ a wheel 
with (resp.\ without) a dot, if it
looks like shown on the left (resp.\ right) side of Figure~\ref{f:wheels}. 

\begin{figure}[!h]
\centering 
\setbox1=\hbox{\input{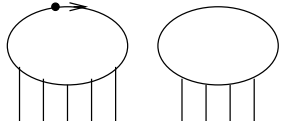}}
$\vcenter{\box1}$
\caption{A wheel with and without a dot} \label{f:wheels}
\end{figure}

As a consequence of Equation~(\ref{e:pstartrival}) we can construct many elements in $\Ker(p^*)$.
An example is given in the following lemma.

\begin{lemma}\label{l:pstarwheeldot}
Let $D$ be a trivalent diagram such that one component of~$D\setminus \Gamma$ is
a wheel with a dot.
Then~$p^*(D)=0$.
\end{lemma}
\begin{proof}

Applying the $(STU)$-relation to the wheel with a dot we obtain the two
diagrams $D_1$ and $D_2$ shown in Figure~\ref{f:D1D2}.

\begin{figure}[!h]
\centering 
\setbox1=\hbox{\input{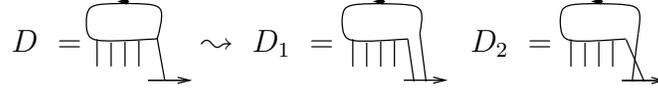}}
$\vcenter{\box1}$
\caption{$D=D_1-D_2$ by the (STU)-relation} \label{f:D1D2}
\end{figure}

By Equation~(\ref{e:pstartrival}) we then have

\begin{eqnarray*}
& & p^*(D) = p^*(D_1-D_2)\\
& = &  p_1(D_1)+(-1)^{\deg(D_1)}p_2(D_1)-p_1(D_2)-(-1)^{\deg(D_2)}p_2(D_2)=0
\end{eqnarray*}

because $p_i(D_1)=p_i(D_2)$.
\end{proof}


Consider an interval~$J$ on the skeleton~$\Gamma={S^1}^{\amalg\ell}$ 
of a diagram~$D$. Cutting~$D$ at the endpoints of~$J$ we obtain a diagram~$D'$
with skeleton~${S^1}^{\amalg\ell-1}\amalg J\amalg (S^1\setminus J)$. 
When $D'$ is the union of diagrams~$D_1$ and~$D_2$ with skeletons~$J$ and
${S^1}^{\amalg\ell-1}\amalg (S^1\setminus J)$, respectively, then we say that
$D_1$ is an {\em isolated part} of~$D$.
The map~$W_\glz^\ob:\Ab^\ob(P^2\times I)\longrightarrow R(P^2\times I)\otimes\Q$ 
has the following 
special property.

\begin{lemma}\label{l:Wnablais0}
Let $D$ be an ordered based $\Z/2$-labeled diagram such that
$\sigma(D)=0$, or 
the number of circles of the skeleton of~$D$ plus the degree of~$D$ is even, 
or $D$ contains an isolated part of odd degree.
Then $W_\glz^\ob(D)=0$.
\end{lemma}

The proof of the lemma follows easily by expanding~$W_\glz^\ob(D)$ as a linear combination of
multiples of~$\overK_e=0$ and~$\overK_s^2=0\in R(P^2\times I)\otimes\Q$.

Some properties of~$W_\glz^\ob$ for ordered based $\Z/2$-labeled diagrams can 
more conveniently be proven using the Lie superalgebra~$\gl(1\vert 1)$.
We identify $\gl(1\vert 1)$ with $2\times 2$-matrices by using a homogeneous basis
of the defining $(1\vert 1)$-dimensional representation of~$\gl(1\vert 1)$.
The rough idea, how the usual description of the weight system associated to~$\gl(1\vert 1)$
(see~\cite{Vai}) extends to
$\Z/2$-labeled diagrams is by mapping a dot on the skeleton to the morphism~$\tau$
in Equation~(\ref{e:dotonarrow}).

\begin{equation}\label{e:dotonarrow}
\mbox{\setbox1=\hbox{\begin{picture}(0,0)%
\epsfig{file=arrowdot.ps}%
\end{picture}%
\setlength{\unitlength}{4144sp}%
\begingroup\makeatletter\ifx\SetFigFont\undefined%
\gdef\SetFigFont#1#2#3#4#5{%
  \reset@font\fontsize{#1}{#2pt}%
  \fontfamily{#3}\fontseries{#4}\fontshape{#5}%
  \selectfont}%
\fi\endgroup%
\begin{picture}(60,381)(409,41)
\end{picture}
}
$\vcenter{\box1}$}\mapsto
\left(\begin{array}{ll}
0& 1\\
1  & 0
\end{array}\right)=\tau
\end{equation}

(compare~\cite{Lie}, 
Section~1.17; recall from there how labels can influence signs). 
Given an ordered based~$\Z/2$-labeled diagram~$D$ with~$\sigma(D)=1$
we cut the first circle of the skeleton of~$D$ at its basepoint. 
Then the construction of weight systems maps the resulting diagram 
to a morphism of the form 
$\lambda\cdot\tau$. It is easy to show 
that~$W_\glz^\ob(D)=\lambda \overK_s$. In particular, 
the number~$\lambda=:W^s(D)$
depends only on~$D\in\Ab^\ob(P^2\times I)$.

%

Define~$W^e:\Ab(S^2\times I)\longrightarrow\Q$ 
by $W^e(D)=\lambda$ whenever $W_\glz^\ob(D)=\lambda \overK_e$.
The following
properties of $W=W^e$ for diagrams without labels are proven in~\cite{Vai}
\footnote{The diagram on the right side of Equation~(\ref{e:Wn2}) is a trivalent
diagram in the sense of this paper (that means it must not 
contain a connected component that does not connect to the skeleton 
of the diagram).}.
They extend to~$W=W^s$ in a straightforward way.

\begin{eqnarray}
W\left(\mbox{\setbox1=\hbox{\begin{picture}(0,0)%
\epsfig{file=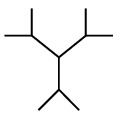}%
\end{picture}%
\setlength{\unitlength}{4144sp}%
\begingroup\makeatletter\ifx\SetFigFont\undefined%
\gdef\SetFigFont#1#2#3#4#5{%
  \reset@font\fontsize{#1}{#2pt}%
  \fontfamily{#3}\fontseries{#4}\fontshape{#5}%
  \selectfont}%
\fi\endgroup%
\begin{picture}(519,488)(439,-87)
\end{picture}
}
$\vcenter{\box1}$}\right)=
W\left(\mbox{\setbox1=\hbox{\begin{picture}(0,0)%
\epsfig{file=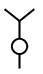}%
\end{picture}%
\setlength{\unitlength}{4144sp}%
\begingroup\makeatletter\ifx\SetFigFont\undefined%
\gdef\SetFigFont#1#2#3#4#5{%
  \reset@font\fontsize{#1}{#2pt}%
  \fontfamily{#3}\fontseries{#4}\fontshape{#5}%
  \selectfont}%
\fi\endgroup%
\begin{picture}(159,294)(439,107)
\end{picture}
}
$\vcenter{\box1}$}\right)=0 & &\label{e:Wn1}\\
W\left(\mbox{\setbox1=\hbox{\begin{picture}(0,0)%
\epsfig{file=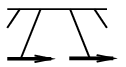}%
\end{picture}%
\setlength{\unitlength}{4144sp}%
\begingroup\makeatletter\ifx\SetFigFont\undefined%
\gdef\SetFigFont#1#2#3#4#5{%
  \reset@font\fontsize{#1}{#2pt}%
  \fontfamily{#3}\fontseries{#4}\fontshape{#5}%
  \selectfont}%
\fi\endgroup%
\begin{picture}(528,261)(457,140)
\end{picture}
}
$\vcenter{\box1}$}\right)=
W\left( \mbox{\setbox1=\hbox{\begin{picture}(0,0)%
\epsfig{file=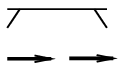}%
\end{picture}%
\setlength{\unitlength}{4144sp}%
\begingroup\makeatletter\ifx\SetFigFont\undefined%
\gdef\SetFigFont#1#2#3#4#5{%
  \reset@font\fontsize{#1}{#2pt}%
  \fontfamily{#3}\fontseries{#4}\fontshape{#5}%
  \selectfont}%
\fi\endgroup%
\begin{picture}(528,261)(457,140)
\end{picture}
}
$\vcenter{\box1}$}\right)& &\label{e:Wn2}\\
W\left(\mbox{\setbox1=\hbox{\begin{picture}(0,0)%
\epsfig{file=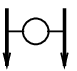}%
\end{picture}%
\setlength{\unitlength}{4144sp}%
\begingroup\makeatletter\ifx\SetFigFont\undefined%
\gdef\SetFigFont#1#2#3#4#5{%
  \reset@font\fontsize{#1}{#2pt}%
  \fontfamily{#3}\fontseries{#4}\fontshape{#5}%
  \selectfont}%
\fi\endgroup%
\begin{picture}(307,307)(455,78)
\end{picture}
}
$\vcenter{\box1}$}\right)=-2
W\left( \mbox{\setbox1=\hbox{\begin{picture}(0,0)%
\epsfig{file=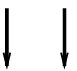}%
\end{picture}%
\setlength{\unitlength}{4144sp}%
\begingroup\makeatletter\ifx\SetFigFont\undefined%
\gdef\SetFigFont#1#2#3#4#5{%
  \reset@font\fontsize{#1}{#2pt}%
  \fontfamily{#3}\fontseries{#4}\fontshape{#5}%
  \selectfont}%
\fi\endgroup%
\begin{picture}(307,307)(455,78)
\end{picture}
}
$\vcenter{\box1}$}\right)& &\label{e:Wn3}
\end{eqnarray}

Equations~(\ref{e:Wn1})--(\ref{e:Wn3}) for~$W=W^e, W^s$, 
Equations~(\ref{e:comult}), (\ref{e:comultob}), 
(\ref{e:pstartrival}), and the (STU)-relation  
imply the following lemma.

\begin{lemma}\label{l:moreW}
Equations~(\ref{e:Wn1})--(\ref{e:Wn3}) also hold for

$$
W=W^e\circ p^*\quad\mbox{and}\quad 
W=(W^s\otimes W^s)\circ\Delta^\ob.$$
\end{lemma}

Now we are ready to prove the following important lemma.

\begin{lemma}\label{l:deltapstar}
Let $D$ be a $\Z/2$-labeled diagram with skeleton~$S^1$ and $\sigma(D)=1$. Then

$$(W^s\otimes W^s)\circ\Delta^\ob(D)=W^e(p^*(D)).
$$
\end{lemma}
\begin{proof}
The proof is devided into four steps.

(1) Assume that a component of~$D\setminus \Gamma$ is
a wheel with a dot. Then~$W^s(D)=0$
(by Equation~(\ref{e:Wn2}) and Lemma~\ref{l:Wnablais0} it is sufficient to verify this
for a degree-$2$ wheel with a dot. This case follows from a straightforward
compuation).

\smallskip

(2) We show now that when the lemma holds for~$D$,
then it also holds for all diagrams~$D'$ with~$D\subset D'$ 
where~$C=D'\setminus D$ is
connected and~$C$ is not a tree. 
If~$C$ is a wheel with a dot then we have

\begin{equation}\label{e:Wis0}
(W^s\otimes W^s)\circ\Delta^\ob(D')=0\quad\mbox{ and }\quad W^e(p^*(D'))=0
\end{equation}

by Part~(1) of this proof and Lemma~\ref{l:pstarwheeldot}.
If~$C$ is not a wheel, then Equation~(\ref{e:Wis0})
follows from Lemma~\ref{l:deltapstar} and Equation~(\ref{e:Wn1}).
If~$C$ is a wheel of odd degree, then Equation~(\ref{e:Wis0})
is implied by the $(STU)$-relation, the previous case (where~$C$ was not a wheel), and
Lemma~\ref{l:Wnablais0}.
Finally, if~$C$ is a wheel of even degree without a dot, 
then we have $W_\glz^\ob(D')=-2W_\glz^\ob(D)$ 
by Equations~(\ref{e:Wn2}) and~(\ref{e:Wn3}) which implies

$$
(W^s\otimes W^s)\circ\Delta^\ob(D')
=-4(W^s\otimes W^s)\circ\Delta^\ob(D)=
-4W^e(p^*(D))=W^e(p^*(D')).
$$

\smallskip

(3) Now assume that the lemma holds for a diagram~$D$ and consider
$D\subset D'$ such that $C=D'\setminus D$ is a tree.
We call the tree~$C$ a comb if we cannot apply Equation~(\ref{e:Wn1}) to the 
part~$C$ of the diagram~$D´$.
By Lemma~\ref{l:moreW} and Equation~(\ref{e:Wn2}) 
we only need to consider combs of degrees~$1$, $2$, $3$, and~$4$.
By the $(STU)$-relation and Part~(2) of this proof,
we may arrange the univalent vertices 
of~$C$ on the skeleton~$S^1$ in any order we want to.
We then apply Relation~(Comult) to reduce the configuration of labels we need to consider.
In the following we investigate the remaining possibilities.

{$\deg\, C= 1:$} 
By Relation~$(FI)$ and Lemma~1.9 of~\cite{Lie} there is nothing to prove in this case
(alternatively, there is a direct proof similar to the case {$\deg\, C= 3$}).

{$\deg\, C=2:$} We only need to consider the case shown in 
Equation~(\ref{e:deg2}), 
where we apply Relations~$(STU)$ and~$(AS)$ to reduce this case to Part~$(2)$ of this proof.

\begin{equation}\label{e:deg2}
\mbox{\setbox1=\hbox{\input{tex/deg2}}
$\vcenter{\box1}$}
\end{equation}

{$\deg\, C=3:$} All possible configurations of $C\subset D'$ 
can be reduced the case shown in 
Equation~(\ref{e:deg3}), where we prove that $p^*(D')=0$.

\begin{equation}\label{e:deg3}
p^*\left(\mbox{\setbox1=\hbox{\begin{picture}(0,0)%
\epsfig{file=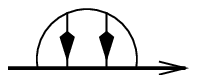}%
\end{picture}%
\setlength{\unitlength}{4144sp}%
\begingroup\makeatletter\ifx\SetFigFont\undefined%
\gdef\SetFigFont#1#2#3#4#5{%
  \reset@font\fontsize{#1}{#2pt}%
  \fontfamily{#3}\fontseries{#4}\fontshape{#5}%
  \selectfont}%
\fi\endgroup%
\begin{picture}(854,301)(429,7)
\end{picture}
}
$\vcenter{\box1}$}\right)
=\mbox{\setbox1=\hbox{\begin{picture}(0,0)%
\epsfig{file=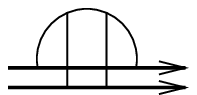}%
\end{picture}%
\setlength{\unitlength}{4144sp}%
\begingroup\makeatletter\ifx\SetFigFont\undefined%
\gdef\SetFigFont#1#2#3#4#5{%
  \reset@font\fontsize{#1}{#2pt}%
  \fontfamily{#3}\fontseries{#4}\fontshape{#5}%
  \selectfont}%
\fi\endgroup%
\begin{picture}(854,391)(429,-38)
\end{picture}
}
$\vcenter{\box1}$}
-\mbox{\setbox1=\hbox{\begin{picture}(0,0)%
\epsfig{file=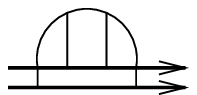}%
\end{picture}%
\setlength{\unitlength}{4144sp}%
\begingroup\makeatletter\ifx\SetFigFont\undefined%
\gdef\SetFigFont#1#2#3#4#5{%
  \reset@font\fontsize{#1}{#2pt}%
  \fontfamily{#3}\fontseries{#4}\fontshape{#5}%
  \selectfont}%
\fi\endgroup%
\begin{picture}(854,391)(429,-38)
\end{picture}
}
$\vcenter{\box1}$}
=0
\end{equation}

By Lemma~\ref{l:Wnablais0} we have 
$(W^s\otimes W^s)\circ\Delta^\ob(D')=0$ because~$\Delta^\ob(D')$ 
is a linear combination of elements~$D_1\otimes D_2$ such that~$D_1$ or~$D_2$ has
an isolated part of odd degree.

{$\deg\, C=4:$} It is sufficient to consider a comb as shown on the left side
of Equation~(\ref{e:deg4}). 

\begin{equation}\label{e:deg4}
\mbox{\setbox1=\hbox{\input{tex/deg4}}
$\vcenter{\box1}$}
\end{equation}

By the $(STU)$-relation and Part~$(2)$
of this proof, we can equivalently consider the linear combination of two diagrams shown
on the right side of Equation~(\ref{e:deg4}). 
But by Relation~$(AS)$ this element equals~$0$.

\smallskip

(4) Conclusion:
For diagrams of degree~$0$ the lemma is obvious. Using 
Parts~(2) and~(3) of this proof the lemma
follows by induction.
\end{proof}

It would be interesting to know to what extent 
Theorem~\ref{t:Csquare} can be generalized to links. A 
counterexample to Lemma~\ref{l:deltapstar} for chord diagrams on 
more than one circle is shown in Figure~\ref{f:counterex}. 

\begin{figure}[!h]
\centering 
\setbox1=\hbox{\input{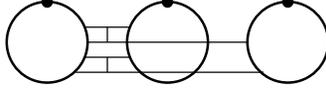}}
$\vcenter{\box1}$
\caption{A diagram $D$ with $(W^s\otimes W^s)\circ\Delta^\ob(D)=0$ 
and $W^e(p^*(D))\not=0$} \label{f:counterex}
\end{figure}

Lemma~\ref{l:deltapstar} is the main ingredient in the following proof of
Theorem~\ref{t:Csquare}.

\begin{proof}[Proof of Theorem~\ref{t:Csquare}]
Define $\Q[[h]]$-valued maps~$\What^e$, $\What^s$ by 
$\What^e(D)=W^e(D) h^{\deg\, D}$ and~$\What^s(D)=W^s(D) h^{\deg\, D}$.
Let~$K^\ob\subset P^2\times I$ be a based knot with $\sigma(K)=1$ and let~$K$ 
be~$K^\ob$ without its basepoint.
Then we have

\begin{eqnarray}
& & \What^e\circ Z_{S^2\times I}(p^{-1}(K))  = \What^e\circ p^*(Z_{P^2\times I}(K))
\label{e:nablahat1}\\
& = & \left(\What^s\otimes \What^s\right)\circ\Delta^\ob(Z_{P^2\times I}(K)) = 
\left(\What^s\circ Z^\ob_{P^2\times I}(K^\ob)\right)^2\label{e:nabla2hat},
\end{eqnarray}

where the first equality follows from Equation~(\ref{e:rp2cover}),
the second equality follows from Lemma~\ref{l:deltapstar}, 
and the last one follows because~$Z_{P^2\times I}(K)$ is group-like.

It remains to verify that Equations~(\ref{e:nablahat1}), (\ref{e:nabla2hat}) are 
compatible with the replacement
of parameters in the proof of Theorem~\ref{t:conwaypoly}:
by Proposition~\ref{p:conquogen} there exist $f,g\in\Z[z]$ such that

\begin{equation}\label{e:twonab}
\nabla(p^{-1}(K))=g(z)\overK_e\quad\mbox{and}\quad\nabla(K^\ob)=f(z)\overK_s.
\end{equation}

By Lemma~\ref{l:Ch} we have

\begin{eqnarray}
& &\What^e\circ Z_{S^2\times I}(p^{-1}(K))=
g(e^{h/2}-e^{-h/2})\What^e\circ Z_{S^2\times I}(K_e)\\
& \mbox{and} &
\What^s\circ Z_{P^2\times I}(K^\ob)=f(e^{h/2}-e^{-h/2})\What^s\circ Z_{P^2\times I}(K_s).
\end{eqnarray}

As a simple consequence of the construction of~$Z_{P^2\times I}$ we have

\begin{equation}\label{e:nablahat5}
\left(\What^s\circ Z_{P^2\times I}(K_s)\right)^2=\What^e\circ Z_{S^2\times I}(K_e).
\end{equation}

Equations~(\ref{e:nablahat1})--(\ref{e:nablahat5}) imply 
$g(z)=f(z)^2$ what we wanted to prove.
\end{proof}

Lemma~\ref{l:deltapstar} 
can be generalized from~$P^2$ to the M\"obius strip~$X$ using the weight system
$W_\glz^\ob\;\mod\;(\overK_i\overK_j\;\vert\; i,j\in\Z)$. A 
generalization of Theorem~\ref{t:Csquare} to this case is not straightforward because
Equation~(\ref{e:twonab}) becomes more complicated and
for ordered based links~$L$ in the ideal $(K_iK_j\;\vert\;i,j\in\Z)\subset\Conw(X\times I)$ 
we may have

\begin{equation}\label{e:not0mod}
\What_\gl^\ob\circ Z_{X\times I} (L)\not \equiv 0\;\mod\; (\overK_i\overK_j\;\vert\; i,j\in\Z)
\end{equation}

(for example for $L=K_{-1} K_1 K_3$).

\section*{Appendix}

In the following we apply the methods of this paper to determine the Homfly and Kauffman
skein modules of cylinders over oriented surfaces with boundary. 
For the Homfly skein module we rediscover the main result of~\cite{Pr1}. Our choice of
representatives of conjugacy classes in~$\pi_1(\Sigma)$ is different to the choice made
in~\cite{Pr1}. Our choice has the advantage that
it allows us to determine the structure of the $2$-variable Kauffman skein module 
with similar methods.
This result seems to be new.

\subsection*{The Homfly skein module}\label{s:aH}

The {\em Homfly skein module} of an oriented $3$-manifold is
generated over the ring
$\Z[x^{\pm 1},y^{\pm 1}]$ by isotopy classes of oriented links modulo the skein relation
of the Homfly polynomial

\begin{equation}\label{e:skeinH}
x \picXps-x^{-1} \picXms = y\picIIs,
\end{equation}

where the links $(\picXps,\picXms,\picIIs)$ differ only locally as shown in the diagrams.
For technical reasons we include the empty link~$\emptyset$ in the definition of
the Homfly skein module and relate it to the trivial knot~$O$ by the 
equation~$(x-x^{-1}) \emptyset  = y O$.

In difference to the Conway skein module, the $3$-manifold must be oriented for
the definition of the Homfly skein module, and 
we need no order or basepoints on~$L$.
Let~$\Sigma$ be an oriented decomposed surface in the sense of Section~\ref{s:surf}. 
Define the ordered set~$\Conjrep$ of representatives of conjugacy classes of elements 
in~$\pi_1(\Sigma)$ as in Section~\ref{s:surf}. 
Then $\Conjrep$-descending links are defined by forgetting the order and basepoints of
the links in Definition~\ref{d:desclink}.
Define~$\Conjrep^\dagger=\Conjrep\setminus\{e \}$, 
where $e\in\pi_1(\Sigma)$ is the neutral element.
For a link~$L\subset\Sigma\times I$
we define~$L^\dagger$ as~$L$ without
the components that are homotopic to~$e$. Let~$n(L)$ be the number of components
of~$L\setminus L^\dagger$.
Let~$\Z[t_{\Conjrep^\dagger}]$ (resp.\ $\Z[t_{\Conjrep}]$) 
be a polynomial ring with indeterminates~$t_w$
in one-to-one correspondence with
elements~$w\in\Conjrep^\dagger$ (resp.\ $w\in\Conjrep$).
For a knot~$K\subset\Sigma\times I$ 
the element~$t_K\in\Z[t_{\Conjrep}]$ is defined as $t_K=t_w$, where $K$ is homotopic 
to~$w$. For a link~$L$ the element~$t_L$ is defined as the product of elements~$t_K$
where~$K$ runs over all
components of~$L$. 
By~$t_\emptyset$ we mean~$1\in\Z[t_{\Conjrep}]$.
Then we have the following theorem.

\begin{theorem}\label{t:Homfly}
There exists a unique invariant 

$$
H(L)\in\Z[x^{\pm 1}, y^{\pm 1}]\otimes_\Z \Z[t_{\Conjrep^\dagger}]
$$

of links~$L$ in~$\Sigma\times I$ that depends only on the class of~$L$ in the 
Homfly skein module
and satisfies

\begin{equation}
H\left(L\right)= t_{L^\dagger} \left(\frac{x-x^{-1}}{y}\right)^{n(L)}\
\mbox{whenever $L$ is $\Conjrep$-descending.}\label{e:Hdesc}
\end{equation}
\end{theorem}

The link invariant~$H$ is called the Homfly polynomial.
As in the proof of Lemma~\ref{l:descgen} we see that 
Equations~(\ref{e:skeinH}) and~(\ref{e:Hdesc}) are
sufficient to calculate~$H(L)$ for every 
link\footnote{This argument is much simpler than the corresponding proof for the 
Conway polynomial in Sections~\ref{s:genc} and~\ref{s:gencx}}.
This implies the
uniqueness stated in Theorem~\ref{t:Homfly}.
Under the condition that~$H$ factors through the Homfly skein module
it can easily be shown that 
Equation~(\ref{e:Hdesc}) is
implied by the same equation for links~$L$ with~$n(L)=0$. 
This implies that the Homfly skein module of~$\Sigma\times I$ is a quotient 
of~$\Z[x^{\pm 1}, y^{\pm 1}]\otimes_\Z \Z[t_{\Conjrep^\dagger}]$.
We give a sketch of the proof of the existence of~$H$ in the rest of this
section.
This will imply that the Homfly skein module of~$\Sigma\times I$ is isomorphic 
to~$\Z[x^{\pm 1}, y^{\pm 1}]\otimes_\Z \Z[t_{\Conjrep^\dagger}]$.

For a group~$G$ 
define~$\A(G)=\A^\ob(G,0)$ and~$\Ab(G)=\Ab^\ob(G,0)$
(see Definitions~\ref{d:A} and~\ref{d:ab}). In these definitions the order and the basepoints 
on chord diagrams and relations~$(Ord)$ and~$(Bas)$ are superfluous because~$\sigma=0$. We
identify the 
degree-$0$ part~$\A_0(G)$ of~$\A(G)$ with the isomorphic polynomial 
ring~$\Q[t_{\Conjrep}]$.
The analogue of Corollary~\ref{c:Wb} is not true for~$\sigma=0$ 
because~$t_e t_w\not=0$ for~$w\in G$.
There are two 
equivalent ways of solving this problem. One way is used in the following, and
the second way is used in Section~\ref{s:aF} for the Kauffman polynomial.
For a $G$-labeled{} chord diagram~$D$, we define the element~$\iota(D)$ by
replacing each chord as shown in Figure~\ref{f:defiota}.

\begin{figure}[h]
$$
\setlength{\unitlength}{15pt}
\begin{picture}(1.6,1)(-0.8,-0.1)
\thicklines
\put(-0.7,1){\vector(0,-1){2}}
\put(0.7,-1){\vector(0,1){2}}
\thinlines
\put(-0.7,0){\line(1,0){1.4}}
\end{picture}\ \leadsto \
\begin{picture}(1.6,1)(-0.8,-0.1)
\thicklines
\put(-0.7,1){\vector(0,-1){2}}
\put(0.7,-1){\vector(0,1){2}}
\thinlines
\put(-0.7,0){\line(1,0){1.4}}
\end{picture} -\frac{1}{2}\left(\
\begin{picture}(1.6,1)(-0.8,-0.1)
\thicklines
\put(-0.7,1){\vector(0,-1){2}}
\put(0.7,-1){\vector(0,1){2}}
\thinlines
\qbezier[80](-0.7,0.4)(-0.3,0.4)(-0.3,0)
\qbezier[80](-0.7,-0.4)(-0.3,-0.4)(-0.3,0)
\end{picture} +
\begin{picture}(1.6,1)(-0.8,-0.1)
\thicklines
\put(-0.7,1){\vector(0,-1){2}}
\put(0.7,-1){\vector(0,1){2}}
\thinlines
\qbezier[80](0.7,0.4)(0.3,0.4)(0.3,0)
\qbezier[80](0.7,-0.4)(0.3,-0.4)(0.3,0)
\end{picture}\ \right)
$$\vspace{8pt}
\caption{The deframing map $\iota$}\label{f:defiota}
\end{figure}

The definition determines a linear map
$\iota:\Ab(G)\longrightarrow \A(G)$, such that 
$p\circ\iota=\id$ where
$p:\A(G)\longrightarrow \Ab(G)$ denotes the canonical 
projection (compare~\cite{BN1}, Exercise~3.16).
The map~$\iota$ induces a map $\iota^*:\A(G)^*\longrightarrow\Ab(G)^*$ 
called deframing projection. 
For~$\sigma=0$ we denote the map~$W_\gl^\ob$ 
(see Proposition~\ref{p:Wglob}) also by~$W_\gl$.
Define

\begin{equation}
\Wgl:\Ab(G)\longrightarrow \Q[t_{\Conjrep}]\ , \quad
\Wgl=W_\gl\circ\iota.
\end{equation}

The map $\Wgl$ is called the weight system of the Homfly
polynomial. 
Proceeding as in Section~\ref{s:eC}, we extend $\Wgl$ to the completion
$\widehat{\A}(G)$ of~$\Ab(G)$
by 

\begin{equation}
\widehat{W}_\gl(D)=\Wgl(D)h^{\deg\,D}\in \Q[t_{\Conjrep}][[h]].
\end{equation}
 
For
a link $L\subset\Sigma\times I$ let $\Hh(L)\in \Q[t_{\Conjrep}][[h]]$ be given
by 
$\Hh(L)=\widehat{W}_\gl\circ Z_{\Sigma\times I}(L)$ where $Z_{\Sigma\times I}$ 
denotes the universal Vassiliev invariant of links in~$\Sigma\times I$.
Denote the inclusion 
map~$\Q[t_{\Conjrep}][[h]]\longrightarrow \Q[t_{\Conjrep}][[h,h^{-1}]]$ by~$i$.
We turn~$\Q[t_{\Conjrep}][[h,h^{-1}]]$ into a $\Z[x^{\pm 1}, y^{\pm 1}]$-module by

\begin{equation}
x\cdot a=e^{t_e h/2} a\quad\mbox{and}\quad y\cdot a=(e^{h/2}-e^{-h/2})a.
\end{equation}

With the notation from above we have the following lemma.

\begin{lemma}\label{l:eHomfly}
The link invariant~$i\circ \Hh$ induces a~$\Z[x^{\pm 1}, y^{\pm 1}]$-linear map 
from the Homfly skein module to $\Q[t_{\Conjrep}][[h,h^{-1}]]$ and satisfies
%
$$\Hh(L)\equiv t_L\,\mod\,h\
\mbox{for every link $L$.}$$
%
%
%
%
\end{lemma}

The proof of Lemma~\ref{l:eHomfly} and the completion of the proof of
Theorem~\ref{t:Homfly} are similar to the proofs in Section~\ref{s:eC}.

\subsection*{The Kauffman skein module}\label{s:aF}

The {\em Kauffman skein module} of an oriented $3$-manifold is
generated over the ring
$\Z[x^{\pm 1},y^{\pm 1}]$ by isotopy classes of framed oriented 
links modulo the relation that reversing the orientation of a component induces the identity
map of the Kauffman skein module
and modulo the skein relations
of the Kauffman polynomial

\begin{eqnarray}
\picuXpl-\picuXml & = & y\left(\picuIIl-\picuInvl\right)\label{e:Kauffmanskeinrel}, \\
\picKt & = & x\picKI,\label{e:kau2}
\end{eqnarray}

where the framed links $(\picuXp,\picuXm,\picuII, \picuInv)$ and $(\picKts,\picKIs)$
differ only locally as shown by the 
diagrams. In diagrams of framed links the framing is assumed to be the 
so-called blackboard framing 
(the framing pointing to the reader).
For technical reasons we include the empty link~$\emptyset$ in the definition of
the Kauffman skein module and relate it to the trivial knot with $0$-framing~$O$ by the 
equation~$(x-x^{-1}+1)\emptyset=y O$.

Let~$\Sigma$ be an oriented decomposed surface in the sense of Section~\ref{s:surf}. 
Define the ordered set~$\Conjrep$ of representatives of conjugacy classes of elements 
in~$\pi_1(\Sigma)$ as in Section~\ref{s:surf}. 
Let $\Conjrep_\pm=\{\min\{a, a^{-1}\}\,\vert\, a\in\Conjrep\}$. 
The set~$\Conjrep_\pm$ is in one-to-one correspondence with
homotopy classes of non-oriented knots.
In this section, 
$\Conjrep$-descending links are defined by forgetting the order and basepoints
of
the links in Definition~\ref{d:desclink} and equipping them with arbitrary framing.

For a diagram~$L$ of a framed oriented link in $\Sigma\times I$, 
define the {\em writhe} of~$L$ as~$\w(L)=k_+-k_-$, 
where~$k_+$ is the number of positive crossings
and~$k_-$ is the number of negative crossings in the diagram.
The writhe is an isotopy invariant of framed oriented links and of framed non-oriented knots.

Define~$\Conjrep^\dagger_\pm=\Conjrep_\pm\setminus\{e \}$, 
where $e\in\pi_1(\Sigma)$ is the neutral element.
For a framed link~$L$ we define~$L^\dagger$, $n(L)$ and $t_L\in\Z[t_{\Conjrep_\pm}]$ 
in the same way as in Section~\ref{s:aH}. 
Then we have the following theorem.

\begin{theorem}\label{t:Kauffman}
There exists a unique invariant 

$$
F(L)\in\Z[x^{\pm 1}, y^{\pm 1}]\otimes_\Z \Z[t_{\Conjrep^\dagger_\pm}]
$$

of framed links~$L$ in~$\Sigma\times I$ that depends only on the class
of~$L$ in the Kauffman skein module
and satisfies

$$
F(L)=t_{L^\dagger}\left(\frac{x-x^{-1}}{y}+1\right)^{n(L)}x^{\w(L)}
\ \mbox{whenever $L$ is~$\Conjrep_\pm$-descending.}$$
\end{theorem}

The link invariant $F$ is called the Kauffman polynomial.
The polynomial $x^{-w(L)} F(L)$ is defined for each framed oriented link~$L$ and does not
depend on the framing of~$L$, hence is an isotopy invariant of oriented links.
As in the proof of Lemma~\ref{l:descgen} we see that the conditions in 
Theorem~\ref{t:Kauffman} are
sufficient to calculate~$F(L)$ for every 
link.
This implies the
uniqueness stated in Theorem~\ref{t:Kauffman}.
We give a sketch of the proof of the existence of~$F$ in the rest of this
section.
This will also imply that the Kauffman skein module of~$\Sigma\times I$ is isomorphic 
to~$\Z[x^{\pm 1}, y^{\pm 1}]\otimes_\Z \Z[t_{\Conjrep_\pm^\dagger}]$.

We say that a circle is locally oriented if it is decomposed into a
finite number of oriented intervals.
Let $C_{\rm loc}(G)$ be the $\Q$-vector space generated by
disjoint unions of locally oriented circles 
with a finite number of distinct points on
oriented parts of the circles labeled{} by elements of~$G$ modulo
homeomorphisms of these diagrams and the
following relations:

\medskip

(Rep):
$\quad\begin{picture}(2,0.5)(-1,-0)
\thicklines
\put(-1,0.25){\vector(1,0){2}}
\put(-0.5,0.25){\circle*{0.15}}
\put(0.4,0.25){\circle*{0.15}}
\put(-0.6,-0.3){\hbox{\hss$a$\hss}}
\put(0.3,-0.3){\hbox{\hss$b$\hss}}
\end{picture}\ = \
\begin{picture}(2,0.5)(-1,0)
\thicklines
\put(-1,0.25){\vector(1,0){2}}
\put(0,0.25){\circle*{0.15}}
\put(-0.25,-0.3){\hbox{\hss$ab$\hss}}
\end{picture}\quad , \quad
\begin{picture}(2,0.5)(-1,0)
\thicklines
\put(-1,0.25){\vector(1,0){2}}
\put(0,0.25){\circle*{0.15}}
\put(-0.4,-0.45){\makebox(0.8,0.5){$e$}}
\end{picture}\ = \
\begin{picture}(2,0.5)(-1,0)
\thicklines
\put(-1,0.25){\vector(1,0){2}}
\end{picture}\ ,$

\medskip

(Ori):
$\quad\begin{picture}(4,0.5)(-2,-0)
\thicklines
\put(-2,0.25){\vector(1,0){4}}
\put(-2,0.25){\vector(1,0){0.75}}
\put(0,0.25){\vector(-1,0){0.75}}
\put(0,0.25){\circle*{0.15}}
\put(-0.2,-0.3){\hbox{\hss$g$\hss}}
\thinlines
\put(-1,0.35){\line(0,-1){0.2}}
\put(1,0.35){\line(0,-1){0.2}}
\end{picture}\ = \
\begin{picture}(4,0.5)(-2,0)
\thicklines
\put(-2,0.25){\vector(1,0){4}}
\put(0,0.25){\circle*{0.15}}
\put(-0.2,-0.3){\hbox{\hss$g^{-1}$\hss}}
\end{picture}$.

\bigskip

Here and in Figure~\ref{f:defWosp} 
we represent elements of $C_{\rm loc}(G)$ graphically by formal linear
combinations of parts of pictures 
of immersions of labeled{}, locally oriented circles. The points where
the local orientation changes are marked by the symbol~$\vert$.

Recall the definition of the vector space~$\A(G)=\A^\ob(G,0)$
(see Definition~\ref{d:A}).
Define a map $\beta_1:\A(G)\longrightarrow C_{\rm loc}(G)$ 
by replacing each chord as shown in
Figure~\ref{f:defWosp} (it is easy to see that~$\beta_1$ is well-defined).

\begin{figure}[!h]
$$
\begin{picture}(1.2,1)(-0.6,-0.1)
\thicklines
\put(-0.5,1){\vector(0,-1){2}}
\put(0.5,1){\vector(0,-1){2}}
\thinlines
\put(-0.5,0){\line(1,0){1}}
\end{picture}\ \leadsto \
\begin{picture}(1.2,1)(-0.6,-0.1)
\thicklines
\put(-0.5,1){\vector(1,-2){1}}
\put(0.5,1){\vector(-1,-2){1}}
\end{picture}
\ - \
\begin{picture}(1.2,1)(-0.6,-0.1)
\thicklines
\qbezier[60](-0.5,1)(-0.5,0.5)(0,0.5)
\qbezier[60](0.5,1)(0.5,0.5)(0,0.5)
\put(0.354,0.646){\vector(-1,-1){0.1}}
\put(-0.354,0.646){\vector(1,-1){0.1}}
\put(0,0.4){\line(0,1){0.2}}
\qbezier[60](-0.5,-1)(-0.5,-0.5)(0,-0.5)
\qbezier[60](0.5,-1)(0.5,-0.5)(0,-0.5)
\put(0.5,-0.9){\vector(0,-1){0.1}}
\put(-0.5,-0.9){\vector(0,-1){0.1}}
\put(0,-0.4){\line(0,-1){0.2}}
\end{picture}
$$\vspace{8pt}
\caption{The map $\beta_1$}\label{f:defWosp}
\end{figure}
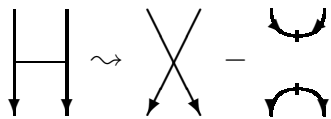

Using Relations $(Rep)$ and $(Ori)$ we can replace each element of
$C_{\rm loc}(G)$ by a linear combination of
disjoint unions of oriented circles labeled{} by a
single element of~$G$. Again by Relations~$(Rep)$ and~$(Ori)$, oriented circles with
labels~$g$, $hgh^{-1}$, and $g^{-1}$ define the same element of~$C_{\rm loc}(G)$.
A map $\beta_2:C_{\rm loc}(G)\longrightarrow \Q[t_{\Conjrep_\pm}]$ is defined
by replacing oriented circles with a label~$g\in\Conjrep_\pm$ by~$t_g$.
It is easy to see that~$\beta_2$ is well-defined.
We define the linear map

\begin{equation}
W_\osp:\A(G)\longrightarrow \Q[t_{\Conjrep_\pm}] ,\quad
W_\osp=\beta_2\circ\beta_1.
\end{equation}

The map~$W_\osp$ is called the weight system
of the Kauffman polynomial. 
Proceeding as in Section~\ref{s:eC}, we extend $W_\osp$ to the completion
$\widetilde{\A}(G)$ of~$\A(G)$
by 

$$
\widetilde{W}_\osp(D)=W_\osp(D)h^{\deg\,D}\in \Q[t_{\Conjrep_\pm}][[h]].$$

For
a framed oriented link $L\subset\Sigma\times I$ 
let~$\Fh(L)\in \Q[t_{\Conjrep}][[h]]$ be given 
by

$$
\Fh(L)=\widetilde{W}_\osp\circ Z^f_{\Sigma\times I}(L)$$ 

where~$Z^f_{\Sigma\times I}$ 
denotes the universal Vassiliev invariant of framed oriented
links in~$\Sigma\times I$.
\footnote{The invariant~$Z^f_{\Sigma\times I}$ is defined by the same formulas 
as~$Z_{\Sigma\times I}$ (see~\cite{Lie}, \cite{LeM}, \cite{BN2}).}
As in Section~\ref{s:aH} we denote the inclusion map~$\Q[t_{\Conjrep_\pm}][[h]]\longrightarrow 
\Q[t_{\Conjrep_\pm}][[h,h^{-1}]]$ by~$i$. This time, we 
turn~$\Q[t_{\Conjrep_\pm}][[h,h^{-1}]]$ into a $\Z[x^{\pm 1}, y^{\pm 
1}]$-module by 

\begin{equation}\label{e:xymodF}
x\cdot a=e^{(t_e-1) h/2} a\quad\mbox{and}\quad y\cdot a=(e^{h/2}-e^{-h/2})a.
\end{equation}

With the notation from above we have the following lemma (compare~\cite{LM1}).

\begin{lemma}\label{l:Kauffman}
The invariant~$i\circ \Fh$ of framed oriented links 
induces a~$\Z[x^{\pm 1}, y^{\pm 1}]$-linear map 
from the Kauffman skein module to $\Q[t_{\Conjrep_\pm}][[h,h^{-1}]]$ and satisfies

$$\Fh(L)  \equiv   t_L\,\mod\,h\
\mbox{for every link $L$.}
$$
\end{lemma}
{\bf Sketch of proof:}
Using $W_\osp(
\begin{picture}(1.1,0.5)(-0.55,-0.05)
\thicklines
\put(-0.5,0){\vector(1,0){1}}
\thinlines
\put(0,0){\line(0,1){0.5}}
\end{picture}
)=-W_\osp(
\begin{picture}(1.1,0.5)(-0.55,-0.05)
\thicklines
\put(0.5,0){\vector(-1,0){1}}
\thinlines
\put(0,0){\line(0,1){0.5}}
\end{picture}
)
$ 
and the analogous dependence of~$Z_{\Sigma\times I}^f$ on orientations 
we see that~$\widetilde{F}(L)$ does not
depend on the orientation of~$L$.

Using
$W_\osp(
\begin{picture}(1.2,0.6)(-0.6,0.1)
\setlength{\unitlength}{11pt}
\thicklines
\put(-1,0.25){\vector(1,0){2}}
\thinlines
\qbezier[70](-0.5,0.25)(-0.5,0.65)(-0.1,0.65)
\qbezier[70](0.3,0.25)(0.3,0.65)(-0.1,0.65)
\end{picture})=(t_e-1)W_\osp(
\begin{picture}(1.2,0.6)(-0.6,0.1)
\setlength{\unitlength}{11pt}
\thicklines
\put(-1,0.25){\vector(1,0){2}}
\end{picture})$ and the explicit description of~$Z_{\Sigma\times I}^f$ 
it is easy to that

\begin{equation}\label{e:Ffrdep}
\Fh\left(\picKt\right)=e^{(t_e-1)h/2}\widetilde{F}\left(\picKI\right).
\end{equation}

The property~$\Fh(L)\equiv t_L\,\mod\,h$ follows directly from the analogous property of the
degre-$0$ part of $\tilde{Z}^f_{\Sigma\times I}$ and the definition of~$W_\osp$.
Therefore it remains only to show that 
Equation~(\ref{e:Fkappaskein}) holds with~$\kappa=1$, $\tilde{y}=e^{h/2}-e^{-h/2}$,
and $(L_+,L_-,L_{\vert\vert},L_=)=\left(\picXps, \picXms, \picIIs, \picuInvmed\right)$,
where the link~$L_=$ has arbitrary orientation.

\begin{equation}\label{e:Fkappaskein}
\widetilde{F}(L_+)-\widetilde{F}(L_-)=
\tilde{y}
\left(\widetilde{F}(L_{\vert\vert})-\kappa \widetilde{F}(L_=)\right)
\end{equation}

For the following computation we extend $\widetilde{W}_\osp$ to
locally oriented labeled chord diagrams in the unique way that
respects the Relation~(Ori). We do not
indicate the local orientation in the picture if it is not of
importance. 

\begin{eqnarray*}
& & \widetilde{F}(L_+)-\widetilde{F}(L_-)=
\widetilde{W}_\osp\left(\exp\left(\picChord /2\right)\circ\picux\right)-
\widetilde{W}_\osp\left(\exp\left(-\picChord /2\right)\circ\picux\right)=\\
& & \sum_{n=0}^\infty\widetilde{W}_\osp\left(\left(
\left(\picux-\picuinf\right)^n
-\left(\picuinf-\picux\right)^n
\right)\circ\picux\right)
\frac{h^n}{2^n n!}=\\
& & \sum_{n=0}^\infty\widetilde{W}_\osp\left(
\left(\picux^n-\picu2+\left(\picu2-\picuinf\right)^n
-\left(-\picux\right)^n+\left(-\picu2\right)^n-\left(\picuinf-\picu2\right)^n
\right)\circ\picux \right)\frac{h^n}{2^n n!}=\\
& & \widetilde{y}\widetilde{F}(L_{\vert\vert})
+\sum_{n=0}^\infty\widetilde{W}_\osp\left(\left(
-\picu2+\left(\picu2-\picuinf\right)^n
+\left(-\picu2\right)^n-\left(\picuinf-\picu2\right)^n\right)\circ\picux
\right)\frac{h^n}{2^n n!}=\\
& & \widetilde{y}\widetilde{F}(L_{\vert\vert})
+\frac{1}{t_e}\sum_{n=0}^\infty\widetilde{W}_\osp\left(
\left(-\picu2+\left(\picu2-\picuinf\right)^n
+\left(-\picu2\right)^n-\left(\picuinf-\picu2\right)^n\right)\circ\picuinf
\right)\frac{h^n}{2^n n!}=\\
& & \widetilde{y}\widetilde{F}(L_{\vert\vert})
+\frac{-e^{h/2}+e^{(1-t_e)h/2}+e^{-h/2}-e^{(t_e-1)h/2}}{t_e}
\widetilde{W}_\osp\left(\picuinf\right)=\\
& &  \widetilde{y}\left(\widetilde{F}(L_{\vert\vert})
-\frac{[t_e-1]_{e^{h/2}}+1}{t_e}
\widetilde{W}_\osp\left(\picuinf\right)\right),
\end{eqnarray*}

where in the last equation we used the notation 

$$
[t_e-1]_{e^{h/2}}=\frac{e^{(t_e-1)h/2}-e^{-(t_e-1)h/2}}{e^{h/2}-e^{-h/2}}.
$$

Let~$O$ be the trivial knot with $0$-framing.
It is easy to see that we have

$$
\widetilde{W}_\osp\left(\mbox{\setbox1=\hbox{\begin{picture}(0,0)%
\epsfig{file=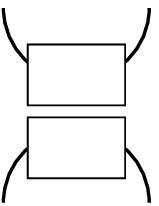}%
\end{picture}%
\setlength{\unitlength}{0.00083300in}%
\begingroup\makeatletter\ifx\SetFigFont\undefined%
\gdef\SetFigFont#1#2#3#4#5{%
  \reset@font\fontsize{#1}{#2pt}%
  \fontfamily{#3}\fontseries{#4}\fontshape{#5}%
  \selectfont}%
\fi\endgroup%
\begin{picture}(732,964)(2368,-1893)
\put(2501,-1353){\makebox(0,0)[lb]{\smash{\SetFigFont{12}{14.4}{\familydefault}{\mddefault}{\updefault}$\ \nu^{1/2}$}}}
\put(2501,-1703){\makebox(0,0)[lb]{\smash{\SetFigFont{12}{14.4}{\familydefault}{\mddefault}{\updefault}$\ \nu^{1/2}$}}}
\end{picture}
}
$\vcenter{\box1}$}
\right)=\frac{\widetilde{F}(O)}{t_e}
\widetilde{W}_\osp\left(\picuInvl \right),
$$

where~$\nu=Z_f(O)$. 
We have $\widetilde{F}(O)\equiv t_e\,\mod\,h$ which implies that~$\widetilde{F}(O)$
is invertible in~$\Q[t_{\Conjrep_\pm}][t_e^{-1}][[h]]$. The computations above show that
Equation~(\ref{e:Fkappaskein}) holds with
 
\begin{equation}\label{e:kappa}
\kappa=([t_e-1]_{e^{h/2}}+1)/\widetilde{F}(O)\in
\Q[t_{\Conjrep_\pm}][t_e^{-1}][[h]].
\end{equation}

\begin{figure}[!h]
\centering
\setbox1=\hbox{\input{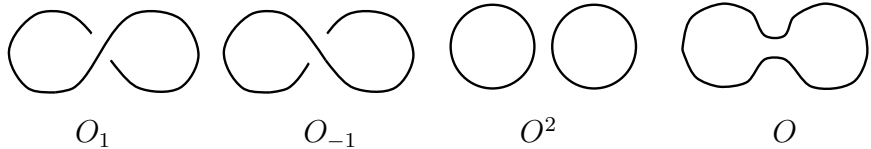}}
$\vcenter{\box1}$
\caption{Link diagrams used to determine~$\widetilde{F}(O)$}\label{f:deteFO}
\end{figure}

Applying Equation~(\ref{e:Fkappaskein}) to the link diagrams of
Figure~\ref{f:deteFO} and using Equation~(\ref{e:Ffrdep}) we obtain

\begin{equation}\label{e:FO}
\widetilde{F}(O^2)=([t_e-1]_{e^{h/2}}+\kappa)\widetilde{F}(O).
\end{equation}

It is easy to see that~$\widetilde{F}(O^2)=\widetilde{F}(O)^2$.
With~$\kappa$ as in Equation~(\ref{e:kappa}) Equation~(\ref{e:FO}) has a
unique solution
satisfying~$\widetilde{F}(O)\equiv t_e\,\mod\,h$,
namely~$\widetilde{F}(O)=[t_e-1]_{e^{h/2}}+1$. This
implies~$\kappa=1$ which completes the proof.
$\Box$

\medskip

Theorem~\ref{t:Kauffman} follows from Lemma~\ref{l:Kauffman} in a similar way as
Theorem~\ref{t:conwaypoly} follows from Lemma~\ref{l:Ch} (see Section~\ref{s:eC}).




}

\begin{thebibliography}{BGRT}

\bibitem[Ale]{Ale} J.\ W.\ Alexander, {\em Topological invariants of knots and links},
Trans.\ Amer.\ Math.\ Soc.\ 30 (1928), 275--306.

\bibitem[AMR]{AMR} J.\ E.\ Andersen, J.\ Mattes, N.\ Reshetikhin, {\em
Quantization of the algebra of chord diagrams}, 
Math.\ Proc.\ Camb.\ Philos.\ Soc.\ 124, No.\ 3 (1998), 451--467.

\bibitem[Bla]{Bla} R.\ C.\ Blanchfield, {\em Intersection theory of manifolds with
operators with applications to knot theory}, Ann.\ of Math.\ 65, No.\ 2 (1957), 340--356.

\bibitem[BN1]{BN1}
D.\ Bar-Natan, {\em On the Vassiliev knot invariants}, Topology 34
(1995), 423--472.

\bibitem[BN2]{BN2}
D.\ Bar--Natan, {\em Non--associative tangles}, Geometric topology
proceedings of the Georgia International Topology Conference (W.\
H.\ Kazez ed.),  139--183,
Amer.\ Math.\ Soc.\ and international Press, Providence (1997).



\bibitem[BNG]{BNG} D.\ Bar-Natan and S.\ Garoufalidis, {\em On the Melvin-Morton-Rozansky conjecture},
Invent.\ Math.\ 125 (1996), 103--133.





\bibitem[Con]{Con} J.\ H.\ Conway, {\em An enumeration of knots and
links, and some of their algebraic properties}, Computational
Problems in Abstract Algebra, Pergamon (1970), 329--358.





\bibitem[HaK]{KaH} R.\ Hartley, A.\ Kawauchi, {\em Polynomials of amphicheiral knots}, 
Math.\ Ann.\ 243 (1979), 63--70.

\bibitem[HOM]{HOM} P.\ Freyd, J.\ Hoste, W.\ B.\ R.\ Lickorish,
K.\ Millet, A.\ Ocneanu and D.\ Yetter,
{\em A new polynomial invariant of knots and links},
Bull.\ Amer.\ Math.\ Soc. 12 (1985), 239--246.

\bibitem[Jon]{Jon} V.\ F.\ R.\ Jones, {\em A polynomial invariant of links via von Neumann
algebras}, Bull.\ Amer.\ Math.\ Soc.\ 12 (1985), 103--111.  




\bibitem[Kau]{Ka2} L.\ H.\ Kauffman, {\em An invariant of regular isotopy}, Trans.\ Am.\
Math.\ Soc.\ 318, No.\ 2 (1990), 417--471.


\bibitem[Kaw]{Kaw} A.\ Kawauchi, {\em A survey of knot theory}, Birkh\"auser, 1996.




\bibitem[Lie]{Lie} J.\ Lieberum, {\em Invariants de Vassiliev pour les 
entrelacs dans~$S^3$ et dans
les vari\'{e}t\'{e}s de dimension trois}, Th\`{e}se de Doctorat de l'Universit\'{e}
Louis Pasteur (Strasbourg~I) (1998).

\bibitem[Li2]{Li2} J.\ Lieberum, {\em Universal Vassiliev invariants of links
in $3$-manifolds}, in preparation.

\bibitem[Liv]{Liv} G.\ R.\ Livesay, {\em Involutions with two fixed points on the three-sphere},
Ann.\ of Math.\ 78, No.\ 2 (1963), 582--593.

\bibitem[LM1]{LM1} T.\ Q.\ T.\ Le and J.\ Murakami, {\em Kontsevich integral for the Kauffman
polynomial}, Nagoya Math. J., 142 (1996), 39--65.

\bibitem[LM2]{LeM} T.\ Q.\ T.\ Le and J.\ Murakami, {\em The universal 
Vassiliev-Kontsevich invariant
for framed oriented links}, Comp. Math.\ 102 (1996), 41--64.



\bibitem[Pr1]{Pr1} J.\ H.\ Przytycki, {\em Skein module of links in a handlebody}, Ohio State
Univ.\ Math.\ Res.\ Inst.\ Publ.\ 1 (1992), 315--342.

\bibitem[Pr2]{Pr2} J.\ H.\ Przytycki, {\em Algebraic topology based on knots : an introduction},
Knots~'96 (Tokyo), 279--297, World Sci.\ Publishing, River Edge, NJ (1997).

\bibitem[Rub]{Rub} J.\ H.\ Rubinstein, {\em Heegard splittings and a theorem of Livesay},
Proc.\ of the Amer.\ Math.\ Soc.\ 60 (1976), 317--320.

\bibitem[Smi]{Smi} P.\ A.\ Smith, {\em Fixed points of periodic transformations},
Appendix~B in S.\ Lefschetz, {\em Algebraic topology}, Amer.\ Math.\ Soc.\ Colloquium 
publications~27, New York (1942).



\bibitem[Vai]{Vai} A.\ Vaintrob, {\em Melvin-Morton conjecture 
and primitive Feynman diagrams}, University of Utah preprint, May 
1996, see q-alg/9605028. 

\bibitem[Vog]{Vog} P.\ Vogel, {\em Invariants de Vassiliev des n\oe uds}, S\'{e}minaire Bourbaki 769
(1993), 1--17, Ast\'{e}risque 216 (1993), 213--232.


\bibitem[Vas]{Vas} V.\ A.\ Vassiliev, {\em Cohomology of knot spaces}, Theory of Singularities
and its Applications (V.\ I.\ Arnold ed.), 23--69, Amer.\ Math.\ Soc., Providence (1990).
\end{thebibliography}
\end{document}